\documentclass[11pt]{article}

\usepackage{ifthen}

\makeatletter
\renewcommand\@cite[2]{%
	#1\ifthenelse{\boolean{@tempswa}}
	{,#2}{}}
\def\@biblabel#1{}
\makeatother

\selectfont

\usepackage{etoolbox}
\AtBeginEnvironment{thebibliography}{\linespread{1}\selectfont}

\usepackage[greek,english]{babel}
\usepackage[utf8]{inputenc}
\usepackage[T1]{fontenc}
\usepackage{lmodern}
\usepackage{geometry}
\usepackage{makecell}
\usepackage{ragged2e}
\usepackage{booktabs}
\usepackage{multirow}
\usepackage[counterclockwise]{rotating}
\usepackage{amssymb,amsmath,latexsym,color,textcomp,marvosym}
\usepackage{ulem,url}
\usepackage[counterclockwise]{rotating}
\usepackage{dsfont}
\usepackage{mathrsfs}
\usepackage{comment}
\usepackage{tabulary}
\usepackage{mathtools}
\usepackage{appendix}
\usepackage{tikz}
\usepackage{float}

\newcommand{\btheta}{\boldsymbol{\theta}}

\newcommand{\brho}{\boldsymbol{\rho}}
\newcommand{\bvartheta}{\boldsymbol{\vartheta}}

\newcommand{\bA}{\boldsymbol{A}}
\newcommand{\ba}{\boldsymbol{a}}

\newcommand{\bx}{\boldsymbol{x}}

\newcommand{\bu}{\boldsymbol{u}}
\newcommand{\bX}{\boldsymbol{X}}

\newcommand{\bs}{\boldsymbol{s}}

\newcommand{\bnu}{\boldsymbol{\nu}}
\newcommand{\bmu}{\boldsymbol{\mu}}

\newcommand{\scal}[1]{\langle #1 \rangle}

\renewcommand{\thefootnote}{\raise.8ex\hbox{\scriptsize{\arabic{footnote}}}}
\small\normalsize

\newcommand{\sign}{\text{sign}}
\newcommand{\tg}{\text{tg}}

\usepackage{etoolbox}
\patchcmd{\thebibliography}{\section*{\refname}}{}{}{}

\def \argmin{\mathop{\hbox{\rm arg min}}}

\def\1g{1\hskip -3pt \mbox{l}}

\newtheorem{theo}{Theorem}
\newtheorem{lem}{Lemma}
\newtheorem{prop}{Proposition}

\newtheorem{rem}{Remark}
\newtheorem{exemple}{Example}

\newtheorem{assumption}{Assumption}

\numberwithin{equation}{section}
\numberwithin{theo}{section}
\numberwithin{lem}{section}
\numberwithin{cor}{section}
\numberwithin{exemple}{section}
\numberwithin{rem}{section}
\numberwithin{prop}{section}
\numberwithin{definition}{section}

\allowdisplaybreaks

\makeatletter
\def\@bibitem#1{\item[]%
	\if@filesw\immediate\write\@auxout{\string \bibcite {#1}{\the\value{\@listctr }}}\fi\ignorespaces}
\makeatletter

\newcommand\blfootnote[1]{%
	\begingroup
	\renewcommand\thefootnote{}\footnote{#1}%
	\addtocounter{footnote}{-1}%
	\endgroup
}

\tikzstyle{bag} = [text width=3em, text centered]
\tikzstyle{bbag} = [text width=5em, text centered]
\tikzstyle{Bbag} = [text width=23em, align = right]

\begin{document}

	
	\begin{center}
		
		{ \Large Conditional Moments of Noncausal Alpha-Stable Processes\\ and the Prediction of Bubble Crash Odds}
		
		
		\vspace*{0.3cm}
		
		{\sc Sébastien Fries}\blfootnote{\textit{Address for correspondence:} Sébastien Fries, 1105 De Boelelaan, 1081HV, Amsterdam, the Netherlands. 
		s.f.fries@vu.nl}\\ \textit{Vrije Universiteit Amsterdam}

		
		\vspace*{0.5cm}
		
		\textbf{Abstract}
	\end{center}
	
	{\justify
		\small
		Noncausal, or anticipative, heavy-tailed processes generate trajectories featuring locally explosive episodes akin to speculative bubbles in financial time series data.
		For $(X_t)$ a two-sided infinite $\alpha$-stable moving average (MA), conditional moments up to integer order four are shown to exist provided $(X_t)$ is anticipative enough, despite the process featuring infinite marginal variance. 
		Formulae of these moments at any forecast horizon under any admissible parameterisation are provided.
		Under the assumption of errors with regularly varying tails, closed-form formulae of the predictive distribution during explosive bubble episodes are obtained  and expressions of the ex ante crash odds at any horizon are available.
		It is found that the noncausal autoregression of order 1 (AR(1)) with AR coefficient $\rho$ and tail exponent $\alpha$ generates bubbles whose survival distributions are geometric with parameter $\rho^{\alpha}$. 
		This property extends to bubbles with arbitrarily-shaped collapse after the peak, provided the inflation phase is noncausal AR(1)-like.
		It appears that mixed causal-noncausal processes generate explosive episodes with dynamics \textit{à la} Blanchard and Watson (1982) which could reconcile rational bubbles with  tail exponents greater than 1.
		Applications of the conditional moments to bubble modelling by noncausal processes are discussed and
		the use of the closed-form crash odds is illustrated on the Nasdaq and S\&P500 series.
		
	}
	
	\vspace*{0.3cm}
	
	{\small
		\noindent \textit{Keywords:} Noncausal process, Conditional moments, Speculative bubble,  Crashes, Prediction, Rational expectation\\
	}
	
	
	
	
	\vspace*{2cm}
	
	\addtocontents{toc}{\protect\setcounter{tocdepth}{-1}}
	
	\section{Introduction}
	
\begin{figure}[h!]
		\centering
		\includegraphics[scale=0.4]{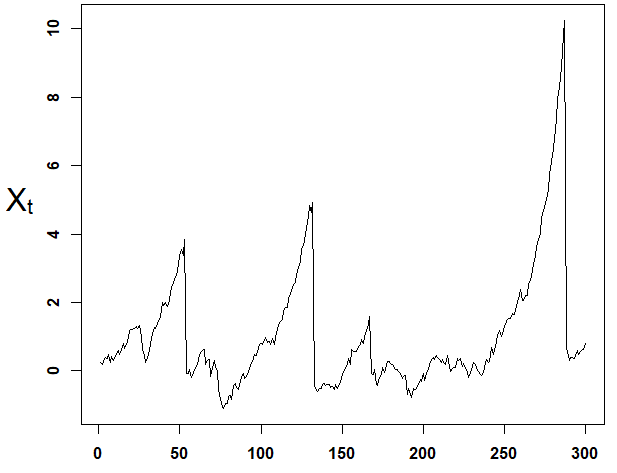}
		\caption{\footnotesize Sample path of an elementary bubble-generating noncausal process: the noncausal AR(1), strictly stationary solution of $X_t = \rho X_{t+1} + \varepsilon_t$, $\rho=0.95$, with $\alpha$-stable errors.
		} 
		\label{fig:traj_ar1}
	\end{figure}

	\noindent Dynamic models often admit solution processes for which the current value of the variable is a function of future values of an independent error process. Such solutions, called \textit{anticipative} or \textit{noncausal}, have attracted increasing attention in the financial and econometric literatures. In particular, noncausal processes have been found convenient for modelling locally explosive phenomena in financial time series such as speculative bubbles, while featuring heavy-tailed marginals and conditional heteroscedastic effects 
	[\cite{bns19, cav17,fri17,gj18,gz17,hs19}, Hecq et al. (2016, 2017a,b), \cite{hen15}] (see also \cite{che17,lns12}, Lanne and Saikkonen (2011, 2013)). 
	Figure \ref{fig:traj_ar1} depicts a typical simulated path of an elementary noncausal process, the $\alpha$-stable noncausal AR(1), featuring multiple bubbles.
	Noncausal processes, shown to be suitable candidates for bubble components in rational expectation price models [\cite{gjm16}], may offer a possibility to forecast the future trajectories of bubbles and to infer the odds of crashes.
	This would enable for instance risk managers to assess large downside risks during prolonged bull markets and the regulator to adjust requirements and restrictions to ensure resilience of the financial system.\\
	\indent However, the limited knowledge about the predictive distribution of noncausal processes, especially during explosive bubble events, is impeding the ability to forecast them, thus limiting their use in practical applications.
	Taking notice of the absence of closed-form formulae for conditional moments and the predictive density except in special cases, two simulation- and sample-based methods have been proposed in the noncausal literature to approximate the conditional distribution of noncausal processes [\cite{lls12,gj16}].
	While offering flexible alternatives for forecasting noncausal processes beyond the special cases, \cite{vh20} find that these methods can become computationally intense for larger prediction horizons and that accurately capturing the dynamics during explosive episodes may prove challenging 
	[see also \cite{ghj19}]. 
	Partial results have been obtained by \cite{gz17} on the conditional moments of noncausal AR(1) processes driven by independent and identically distributed (i.i.d.) $\alpha$-stable errors, which have been extended to mixed causal-noncausal AR processes with single ill-located root by \cite{fri17}. 
	Despite stable noncausal processes featuring infinite marginal variance, their conditional moments may exist up to integer order four. 
	In special cases, expressions of the conditional expectation and variance have been obtained, and revealed that noncausal processes can feature GARCH type effects in calendar time despite such effects not being explicitly included in the modelling. 
	Provided the expressions of the conditional moments are derived, this suggests that point forecasts of noncausal processes based on their conditional expectation, variance, skewness and kurtosis could be formulated -as opposed to other predictors specifically introduced to circumvent the infinite variance of $\alpha$-stable processes, such as minimum $L^\alpha${-dispersion} or maximum {covariation}  (see \cite{kar13} and the references therein).\\
	\indent The aim of this paper is to provide practical analytical results to compute the conditional moments of $\alpha$-stable noncausal processes and to compute the crash odds of bubbles that such processes generate.
	First, the paper extends the literature on the conditional moments $\mathbb{E}[X_2^p|X_1]$ of arbitrary bivariate $\alpha$-stable random vectors $(X_1,X_2)$ [Cioczek-Georges and Taqqu (1995a,b, 1998), \cite{har91}, Samorodnitsky and Taqqu (1994) (\cite{st94} hereafter)] by providing formulae for the orders $p=2,3,4$. 
	We then apply these results to derive a complete characterisation of the conditional moments $\mathbb{E}[X_{t+h}^p|X_t]$, $p=1,2,3,4$, $h\ge1$, for $(X_t)$ an infinite \textit{two-sided} moving average process driven by i.i.d. $\alpha$-stable errors
	\vspace*{-0.2cm}
	\begin{align}\label{def:2sidedma} 
	X_t & = \sum_{k\in\mathbb{Z}} a_k \varepsilon_{t+k}, 
	\end{align}
	
	\vspace*{-0.cm}
	
	\noindent where $(a_k)$ is a non-random coefficients sequence satisfying mild conditions for $(X_t)$ to be well defined and strictly stationary.
	Second, the conditional distribution of noncausal processes during explosive bubble episodes is analysed.
	Provided the errors have probability tails similar to that of $\alpha$-stable distributions, in the sense that they also feature power-law tails, we obtain closed-form formulae valid during explosive episodes.
	These expressions provide illuminating interpretations on the dynamics of the bubbles that such models generate and a practical way to quantify the crash odds.
	Implications and parallels with the literature on rational expectation bubble models in the line of \cite{bw82} are discussed.\\
	\indent The rest of the paper is organised as follows. Section \ref{sec:alpha-vector} recalls properties of bivariate stable distributions and provides our results on the conditional moments up to order four of arbitrary bivariate $\alpha$-stable vectors.
	Applying these results to models of the form \eqref{def:2sidedma}, Section \ref{sec:processes} proposes a sufficient condition on the coefficients $(a_k)$ for the existence of conditional moments, characterises their expressions when they exist, and discusses several examples and methodological aspects.
	Section \ref{sec:crash_odds} derives closed-form formulae for the predictive distribution of noncausal processes during explosive bubble episodes.
	Section \ref{sec:applications} 
	proposes an application of the crash odds formulae on the Nasdaq and S\&P500 series. 
	Proofs and complementary results are collected in a Supplementary File.
	
	\section{Conditional moments of bivariate $\alpha$-stable vectors}
	\label{sec:alpha-vector}
	
	We begin by recalling some properties of bivariate stable vectors $(X_1,X_2)$
	and then propose new expressions for their higher-order conditional power moments $\mathbb{E}[X_2^p|X_1]$. 
	These expressions will apply to $(X_t,X_{t+h})$ when considering $\alpha$-stable noncausal processes in the next section.
	Letting $\alpha\in(0,2)$, a random vector $\bX=(X_1,X_2)$ is said to be an $\alpha$-stable random vector in $\mathbb{R}^2$ (see Theorem 2.3.1 in \cite{st94}) if there exists a unique pair $(\Gamma,\bmu^0)$, where $\Gamma$ is a finite measure on the Euclidean unit sphere $S_2$ and $\bmu^0$ a vector in $\mathbb{R}^2$, such that, for any $\bu\in\mathbb{R}^2$, the characteristic function of $\bX$ writes
	\begin{align}\label{def:char_fun_vec}
	\mathbb{E}\Big[e^{i\scal{\bu,\bX}}\Big] = \exp\bigg\{-\int_{S_2}|\scal{\bu,\bs}|^\alpha \bigg(1-i\,\sign(\scal{\bu,\bs})w(\alpha,\scal{\bu,\bs})\bigg)\Gamma(d\bs)+i\,\scal{\bu,\bmu^0}\bigg\},
	\end{align}
	where $\scal{\cdot,\cdot}$ is the canonical inner product, $w(\alpha,s) = \tg\left(\frac{\pi \alpha}{2}\right)$, if $\alpha \ne 1$, and $w(1,s) = - \frac{2}{\pi}\ln|s|$ otherwise, for $s\in\mathbb{R}$. The measure $\Gamma$ and the vector $\bmu^0$ are respectively called the \textit{spectral measure} and the \textit{shift vector} of $\bX$. The pair $(\Gamma,\bmu^0)$ is said to be the \textit{spectral representation} of $\bX$. 
    The  spectral measure $\Gamma$ of a stable vector $\bX$ in particular  completely characterises the tail dependence between its components: it holds that $\mathbb{P}(\bX/\|\bX\|\in A \,\, | \,\,\|\bX\| > x) \rightarrow \Gamma(A)/\Gamma(S_2)$ as $x\rightarrow+\infty$, for any continuity set $A\subset S_2$ (Theorem 4.4.8 in \cite{st94}), where $\|\cdot\|$ denotes the Euclidean norm.
    Intuitively, the more mass $\Gamma$ attributes to some points of the unit sphere $S_2$, the more likely $(X_1,X_2)$ is to be colinear to these points when it is large in norm. 
	The  counterpart of \eqref{def:char_fun_vec} for a univariate $\alpha$-stable variable reads $
	\mathbb{E}[e^{iuX}] = \exp\Big\{-\sigma^{\alpha}|u|^{\alpha}\Big(1-i\beta \,\sign(u)w(\alpha,u)\Big)+iu\mu\Big\}$, for some asymmetry $\beta\in[-1,1]$, scale $\sigma>0$ and location $\mu\in\mathbb{R}$. 
	
	Stable distributions are known to have very few moments. 
	However, the distribution of one component conditionally on the other can have more moments according to the degree of dependence between them. If the spectral measure $\Gamma$ of an $\alpha$-stable random vector $\bX=(X_1,X_2)$ satisfies
	\vspace*{-0.1cm}
	\begin{equation}\label{eq:nu_cond}
	\int_{S_2} |s_1|^{-\nu}\Gamma(ds)<+\infty, \quad \text{for some} \quad \nu\ge0,   
	\end{equation}
	then, $\mathbb{E}\big[|X_2|^\gamma\big|X_1=x\big]<+\infty$ for almost every $x$ if $0\le\gamma<\min(\alpha+\nu,2\alpha+1)<5$ (see Theorem 5.1.3 in \cite{st94} for details), entailing that conditional moments up to integer order four may exist although $\mathbb{E}\big[|X_2|^\alpha]=+\infty$.
	The conditional expectation of arbitrary $\alpha$-stable bivariate vectors has been studied in details and its expression is recalled
	 in Theorem \ref{theo:cond_expect_stable_general_alpnot1}
	below, while the conditional variance received attention most exclusively in the symmetric $\alpha$-Stable case. \\ 
	\hspace*{0.5cm} We provide and prove new formulae for the conditional power moments of order 2, 3, and 4 of arbitrary (not necessarily symmetric) $\alpha$-stable bivariate vectors $(X_1,X_2)$. 
	The second order moment in the case $\alpha=1$, which requires special treatment when not restricting to symmetric stable distributions, is also considered.
	In the rest of this section, we assume without loss of generality that the shift vector $\bmu^0=(\mu_1^0,\mu_2^0)$ is zero. This can be done without loss of generality because, assuming the conditional moment of order $p$ exists, $\mathbb{E}\big[X_2^p\big|X_1=x\big]  = \mathbb{E}\big[(X_2-\mu_2^0+\mu_2^0)^p\big|X_1-\mu_1^0=x-\mu_1^0\big] = \sum_{j=0}^p \text{C}_p^j (\mu_2^0)^{p-j}\mathbb{E}\big[\tilde{X}_2^j\big|\tilde{X}_1=\tilde{x}\big]$
	where $\tilde{x}=x-\mu_1^0$, and $(\tilde{X}_1,\tilde{X}_2)=(X_1-\mu_1^0,X_2-\mu_2^0)$ has the same spectral measure as $(X_1,X_2)$ and zero shift parameter. 
	We first consider the case $\alpha\ne1$ and introduce useful constants and functions which generalise existing quantities in the literature. For $p\in\{1,2,3,4\}$, when they exist, define
		\begin{equation}\label{def:sbkl}
	\begin{aligned}
	\sigma_1^\alpha  = \int_{S_2}|s_1|^\alpha\Gamma(d\bs), \hspace*{0.2cm} \beta_1  = \dfrac{\int_{S_2}s_1^{<\alpha>}\Gamma(d\bs)}{\sigma_1^\alpha}, 
	\hspace*{0.2cm} \kappa_p  = \dfrac{\int_{S_2} (\frac{s_2}{s_1})^p |s_1|^\alpha\Gamma(d\bs)}{\sigma_1^\alpha}, \hspace*{0.2cm} \lambda_p  = \dfrac{\int_{S_2}(\frac{s_2}{s_1})^p s_1^{<\alpha>} \Gamma(d\bs)}{\sigma_1^\alpha},
	\end{aligned}
	\end{equation}
	where $y^{<r>}=\sign(y)|y|^r$ for any $y,r\in\mathbb{R}$. For any $n\in\mathbb{N}$, $\btheta=(\theta_{1},\theta_{2})\in\mathbb{R}^2$, $x\in\mathbb{R}$, define $\mathcal{H}$ as
	\begin{align}\label{def:hcal_res}
	\mathcal{H}(n,\btheta;x) & =  \int_0^{+\infty} e^{-\sigma_1^\alpha u^\alpha} u^{n(\alpha-1)} \Big(\theta_{1}\cos(ux-a\beta_1\sigma_1^\alpha u^\alpha) + \theta_{2}\sin(ux-a\beta_1\sigma_1^\alpha u^\alpha)\Big)du.
	\end{align}
	The quantities $\sigma_1$ and $\beta_1$ denote the scale and asymmetry parameters of the marginal distribution of $X_1$, whereas the constants $\kappa_p$'s and $\lambda_p$'s generalise standard dependence measures invoked in the literature.
	Noticeably, $\kappa_1=\int_{S_2}s_2 s_1^{<\alpha-1>}\Gamma(d\bs)/\sigma_1^\alpha$ corresponds to the normalised \textit{covariation} between $X_2$ and $X_1$. 
	This dependence measure was been introduced by \cite{mil78} and \cite{cm81} to replace the ill-defined covariance between two symmetric $\alpha$-stable random variables, and has been a popular tool to formulate point forecasts of infinite variance $\alpha$-stable processes [see \cite{kar13} and the references therein].
	The new constants $\kappa_p$ and $\lambda_p$, $p\ge 2$ introduced here, which intervene in the expressions of the higher order conditional moments of $(X_1,X_2)$, can be seen as extending this dependence measure to higher powers of $X_1$ and $X_2$ in the asymmetric case. 
	The new family of functions $\mathcal{H}$ introduced contains functions related to the marginal density of the stable random variable $X_1\sim\mathcal{S}(\alpha,\beta_1,\sigma_1,0)$, $\alpha\ne1$: $f_{X_1}(x) := \frac{1}{\pi}\mathcal{H}\big(0,(1,0);x\big) = \frac{1}{\pi} \int_0^{+\infty} e^{-\sigma_1^\alpha u^\alpha} \cos(ux-a\beta_1\sigma_1^\alpha u^\alpha)du$.
	The following result recalls the expression of the conditional expectation in the case $\alpha\ne1$.
	\begin{theo}[Theorem 5.2.2, \cite{st94}]\label{theo:cond_expect_stable_general_alpnot1}
		Let $(X_1,X_2)$ be an $\alpha$-stable random vector with spectral representation $(\Gamma,\boldsymbol{0})$. 
		For $\alpha\in(0,2)\setminus \{1\}$ and letting $\Gamma$ satisfy \eqref{eq:nu_cond} with $\nu > 1-\alpha$ if $\alpha \in (0,1)$,
		\begin{align}\label{eq:1cm_gen_recall_alpnot1}
		& \mathbb{E}\Big[X_{2}\Big|X_1=x\Big] = \kappa_1 x  + \dfrac{a(\lambda_1-\beta_1\kappa_1)}{1+a^2\beta_1^2} \Bigg[a \beta_1x + \dfrac{1-xH(x)}{\pi f_{X_1}(x)}\Bigg], \hspace*{5cm}
		\end{align}
		where $a=\tg(\pi\alpha/2)$, $\beta_1$, $\kappa_1$ and $\lambda_1$ are as in \eqref{def:sbkl} 
		and $H(\,\cdot\,):=\mathcal{H}\big(0,(0,1);\,\cdot\,\big)$.
	\end{theo}
	We now state our result in the case $\alpha\ne1$ for the conditional moments of order two, three and four.
	\begin{theo}\label{theo:cond_moments_stable_general_alpnot1}
		Let $(X_1,X_2)$ be an $\alpha$-stable random vector with spectral representation $(\Gamma,\boldsymbol{0})$.
		
		\vspace*{-1.cm}
		
		\begin{align}
		& \text{For } \alpha\in(1/2,2)\setminus \{1\} \text{ and } \Gamma \text{ satisfying } \eqref{eq:nu_cond} \text{ with } \nu > 2-\alpha, \nonumber\\
		&\mathbb{E}\Big[X_{2}^2\Big|X_1=x\Big] = \kappa_2 x^2  + \frac{ax(\lambda_2-\beta_1\kappa_2)}{1+(a\beta_1)^2} \Bigg[a\beta_1x + \dfrac{1-xH(x)}{\pi f_{X_1}(x)}\Bigg]  \label{eq:scm_gen}\\
		& \hspace{4cm} - \dfrac{\alpha^2\sigma_1^{2\alpha}}{\pi f_{X_1}(x)}\mathcal{H}\Big(2,\btheta_1;x\Big). \nonumber\\
		& \text{For } \alpha\in(1,2) \text{ and } \Gamma \text{ satisfying } \eqref{eq:nu_cond} \text{ with } \nu > 3-\alpha,  \nonumber\\
		& \mathbb{E}\Big[X_{2}^3\Big|X_1=x\Big] = \kappa_3x^3  + \dfrac{ax^2(\lambda_3-\beta_1\kappa_3)}{1+(a\beta_1)^2} \bigg[a\beta_1x + \dfrac{1-xH(x)}{\pi f_{X_1}(x)}\bigg]\label{eq:tcm_gen}\\
		& \hspace{4cm} -\dfrac{\alpha^2\sigma_1^{2\alpha}}{2 \pi f_{X_1}(x)} \bigg[x\mathcal{H}\Big(2,\btheta_2;x\Big) + \alpha\sigma_1^\alpha\mathcal{H}\Big(3,\btheta_3;x\Big)\bigg].\nonumber\\
		& \text{For } \alpha\in(3/2,2) \text{ and } \Gamma \text{ satisfying } \eqref{eq:nu_cond} \text{ with } \nu > 4-\alpha,  \nonumber\\
		& \mathbb{E}\Big[X_{2}^4\Big|X_1=x\Big]  = \kappa_4 x^4  + \dfrac{ax^3(\lambda_4-\beta_1\kappa_4)}{1+(a\beta_1)^2} \bigg[a\beta_1x+\dfrac{1-xH(x)}{\pi f_{X_1}(x)}\bigg]\label{eq:fcm_gen}\\
		& \hspace{4cm} -\dfrac{\alpha^2\sigma_1^{2\alpha}}{\pi f_{X_1}(x)} \bigg[\dfrac{x^2}{2}\mathcal{H}\Big(2,\btheta_4;x\Big)+\dfrac{\alpha x\sigma_1^{\alpha}}{6}\mathcal{H}\Big(3,\btheta_5;x\Big)+\dfrac{\alpha^2\sigma_1^{2\alpha}}{3}\mathcal{H}\Big(4,\btheta_6;x\Big)\bigg].\nonumber
		\end{align}
		Here,  $a = \text{tg}\,(\pi \alpha/2)$, 
		$H(\,\cdot\,)=\mathcal{H}\big(0,(0,1);\,\cdot\,\big)$, $\btheta_1=(\theta_{11},\theta_{12})$ in \eqref{eq:scm_gen} is given by
		\vspace*{-0.2cm}
		\begin{align*}
		\theta_{11} & = \kappa_1^2-a^2\lambda_1^2 + a^2\beta_1\lambda_2 - \kappa_2,
		& \theta_{12} & = a(\lambda_2+\beta_1\kappa_2)-2a\lambda_1\kappa_1,
		\end{align*}
		
		\vspace*{-0.2cm}
		
		\noindent and the remaining $\btheta_i$'s in \eqref{eq:tcm_gen}-\eqref{eq:fcm_gen}, which depend only on $\alpha$, $\beta_1$, and the $\kappa_p$'s and $\lambda_p$'s in \eqref{def:sbkl}, are given in \eqref{eq:th2}-\eqref{def:nu_fourthcm2} in the Supplementary File. 
	\end{theo}
	\textit{Proof.} See Sections \ref{sec:prelim_elem}, \ref{sec:lemc2} and \ref{sec:prooftheo22} in the Supplementary File.\\
	Let us now turn to the case $\alpha=1$. 
	The following result  recalls the expression of the conditional expectation in this case.
	\begin{theo}[Theorem 5.2.3, ST94]
	\label{theo:cond_moments_stable_general_alp1_expect}
	Let $(X_1,X_2)$ be $\alpha$-stable, with $\alpha=1$ and spectral representation $(\Gamma,\boldsymbol{0})$, where $\Gamma$ satisfies \eqref{eq:nu_cond} with $\nu>0$. Then, for almost every $x$,
	$$
	\mathbb{E}[X_2|X_1=x] = -a\sigma_1 q_0 + \kappa_1(x-\mu_1) + \dfrac{\lambda_1-\beta_1\kappa_1}{\beta_1}\bigg[(x-\mu_1) - \sigma_1\dfrac{U(x)}{\pi f_{X_1}(x)}\bigg],
	$$
	if $\beta_1\ne0$, and 
	$$
	\mathbb{E}[X_2|X_1=x] = -a\sigma_1 q_0 + \kappa_1(x-\mu_1) - a\sigma_1\lambda_1\dfrac{V(x)}{\pi f_{X_1}(x)},
	$$
	if $\beta_1=0$.
	Here, $a=2/\pi$, $\sigma_1$, $\beta_1$, $\kappa_1$ and $\lambda_1$ are as in \eqref{def:sbkl},  $f_{X_1}$ is the marginal density of $X_1\sim\mathcal{S}(1,\beta_1,\sigma_1,\mu_1)$, $q _0  = \dfrac{1}{\sigma_1} \int_{S_2} s_2\ln |s_1|\Gamma(d\bs)$, $\mu_1 = - a \int_{S_2} s_1 \ln|s_1| \Gamma(d\bs)$, and U, V are given in \eqref{def:uvw1}-\eqref{def:uvw2} in the Supplementary File.
	\end{theo}
	We next provide our result for the second order conditional moment when $\alpha=1$. 
	As for the conditional expectation, two different expressions hold according to whether the marginal distribution of $X_1$ is skewed or symmetric. 
	\begin{theo}\label{theo:cond_moments_stable_general_alp1_scm}
		Let $(X_1,X_2)$ be $\alpha$-stable, with $\alpha=1$ and spectral representation $(\Gamma,\boldsymbol{0})$, where $\Gamma$ satisfies \eqref{eq:nu_cond} with $\nu>1$. Then, for almost every $x$,
		\begin{align*}
		\mathbb{E}\Big[X_2^2\Big|X_1=x\Big] & = \sigma_1^2(a^2q_0^2-\kappa_1^2) + \dfrac{2\sigma_1\lambda_1}{\beta_1}\Big(\sigma_1\kappa_1 - aq_0(x-\mu_1)\Big) + \dfrac{\lambda_2}{\beta_1}\Big((x-\mu_1)^2 - \sigma_1^2\Big)\\
		& \hspace{1cm} + \Big(a\sigma_1q_0(\lambda_1-\beta_1\kappa_1) + (\kappa_1\lambda_1-\lambda_2)(x-\mu_1)\Big)\dfrac{2\sigma_1U(x)}{\beta_1\pi f_{X_1}(x)}\\
		& \hspace{1cm} + \Big(\lambda_2+\beta_1\kappa_2-2\kappa_1\lambda_1 + a^2\sigma_1\beta_1(\lambda_1^2 - \beta_1\lambda_2) W(x) \Big) \dfrac{\sigma_1}{\beta_1 \pi f_{X_1}(x)},
		\end{align*}
		if $\beta_1\ne0$, and
		\begin{align*}
		&\mathbb{E}\Big[X_2^2\Big|X_1=x\Big]  = \sigma_1^2(\kappa_2 + a^2q_0^2 - \kappa_1^2 ) - 2a\sigma_1\kappa_1q_0(x-\mu_1) + \kappa_2(x-\mu_1)^2\\
		& \hspace{1cm} + a\sigma_1(\lambda_2-2\lambda_1\kappa_1)\dfrac{F_{X_1}(x) - 1/2}{f_{X_1}(x)}+ \dfrac{a\sigma_1\lambda_1}{\pi f_{X_1}(x)} \bigg[2\Big(a\sigma_1q_0-\kappa_1(x-\mu_1)\Big)V(x) + a\sigma_1\lambda_1W(x)\bigg],
		\end{align*}
		if $\beta_1=0$. Here, $a=2/\pi$, $\sigma_1$, $\beta_1$, the $\kappa_p$'s and the $\lambda_p$'s are as in \eqref{def:sbkl},  $f_{X_1}$ 
		and $F_{X_1}$ are respectively the marginal density and cumulative distribution function of $X_1\sim\mathcal{S}(1,\beta_1,\sigma_1,\mu_1)$, $q _0  = \dfrac{1}{\sigma_1} \int_{S_2} s_2\ln |s_1|\Gamma(d\bs)$, $\mu_1 = - a \int_{S_2} s_1 \ln|s_1| \Gamma(d\bs)$, and U, V and W are given in \eqref{def:uvw1}-\eqref{def:uvw3} in the Supplementary File. 
	\end{theo}
	\textit{Proof.} See Section \ref{sec:tha4} in the Supplementary File.\\
	The expressions of the conditional moments simplify when one considers the asymptotics with respect to the conditioning variable, as $X_1=x$ becomes large.
	\begin{prop}\label{prop:equiv}
		Let $p\in\{1,2,3,4\}$ and let $(X_1,X_2)$ be $\alpha$-stable with $\alpha\in(0,2)$, and spectral representation $(\Gamma,\boldsymbol{0})$ such that the conditional moment of order $p$ exists. If $|\beta_1|\ne1$, then
		\begin{align*}
		x^{-p}\,\mathbb{E}\Big[X_{2}^p\Big|X_1=x\Big]  \underset{x\rightarrow+\infty}{\longrightarrow} \dfrac{\kappa_p + \lambda_p}{1+\beta_1}, \hspace{2cm}	x^{-p}\,\mathbb{E}\Big[X_{2}^p\Big|X_1=x\Big]  \underset{x\rightarrow-\infty}{\longrightarrow} \dfrac{\kappa_p - \lambda_p}{1-\beta_1},
		\end{align*}
		and if $|\beta_1|=1$ and $\beta_1x\rightarrow+\infty$, then, $x^{-p}\,\mathbb{E}\Big[X_{2}^p\Big|X_1=x\Big]  {\longrightarrow} \kappa_p$.
	\end{prop}
	\textit{Proof.} See Section \ref{sec:prop_asymp} in the Supplementary File.

	\section{Conditional moments of noncausal $\alpha$-stable processes}\label{sec:processes}
	
	Operating the set of properties of bivariate $\alpha$-stable distributions provided in the previous section, we study the existence and expressions of the conditional moments of $\alpha$-stable infinite moving average processes. 
	Discussions on practical aspects as well as examples focusing on modelling practices of the empirical noncausal literature follow the main result.
	Let us consider $(X_t)$ a two-sided MA($\infty$) process as in \eqref{def:2sidedma} with $\alpha$-stable errors $\varepsilon_t\stackrel{i.i.d.}{\sim}\mathcal{S}(\alpha,\beta,\sigma,\mu)$ and coefficients $(a_k)$ satisfying
	\begin{gather}
	\sum_{k\in\mathbb{Z}}|a_{k}|^s < +\infty, \text{ for some } s\in(0,\alpha)\cap[0,1], \label{def:2sidedma_cond1}\\ 
	\text{ and in addition for } \alpha=1, \beta\ne0, \hspace{0.3cm} \sum_{k\in\mathbb{Z}}|a_{k}| \Big|\!\ln|a_k|\Big| < +\infty. \label{def:2sidedma_cond2}
	\end{gather}
	Conditions \eqref{def:2sidedma_cond1}-\eqref{def:2sidedma_cond2} ensure that $\sum_{k\in\mathbb{Z}} a_k \varepsilon_{t+k}$ converges absolutely almost surely so that $(X_t)$ is well defined and strictly stationary.
	A moving average process of the form \eqref{def:2sidedma} satisfying the above conditions is said to be purely causal if $a_k=0$ for $k>0$ and purely noncausal if $a_k=0$ for $k<0$.
	Noncausality is found to be crucial for the existence of conditional moments higher than order $\alpha$.
	An important class of models that we shall consider and which admits MA($\infty$) representations satisfying the above conditions is the class of ARMA processes.
	General ARMA processes --causal, noncausal, invertible or non-invertible-- are strictly stationary  solutions of stochastic recursive  equations of the form
	\begin{equation}\label{eq:marma_intro}
	\psi(F)\phi(B)X_t = \Theta(F)H(B)\varepsilon_t,
	\end{equation}
	where $F$ (resp. $B=F^{-1}$) denotes the forward (resp. backward) operator, $\psi(z):=1-\psi_1z-\ldots-\psi_pz^p$ and $\phi(z):=1-\phi_1z-\ldots-\phi_qz^q$ are polynomials of degrees $p$ and $q$, and $H$ and $\Theta$ are two polynomials of respective degrees $r$ and $s$ with roots on or outside the unit circle.
	Equation \eqref{eq:marma_intro} admits a unique strictly stationary solution provided that $\psi(z)\ne0$, $\phi(z)\ne0$ for $|z|\le1$, and that $\psi$ (resp. $\phi$) has no common root with $\Theta$ (resp. $H$). 
	The stationary solution is noncausal if $p\ge1$.
	
	\subsection{Spectral representation of $(X_t,X_{t+h})$}
	
	Because the error sequence $(\varepsilon_t)$ is $\alpha$-stable distributed, the bivariate vector $(X_t,X_{t+h})$, for $(X_t)$ satisfying \eqref{def:2sidedma}, \eqref{def:2sidedma_cond1} and \eqref{def:2sidedma_cond2}, is itself $\alpha$-stable for any horizon $h$ and the results from the previous section apply.
	This is a consequence of the following lemma, which provides the spectral representation of discrete time vectors of linear moving averages driven by $\alpha$-stable i.i.d. errors. 
	\begin{lem}
		\label{prop:multistable}
		Let $0<\alpha<2$.
		For $\varepsilon_t\stackrel{i.i.d.}{\sim}\mathcal{S}(\alpha,\beta,\sigma,\mu)$ and real deterministic sequences $(a_{k,i})_k$, $i=1,2$, both satisfying \eqref{def:2sidedma_cond1}-\eqref{def:2sidedma_cond2}, let $\boldsymbol{X_t}=(X_{1,t},X_{2,t})$, with $X_{i,t} = \sum_{k\in\mathbb{Z}}a_{k,i}\varepsilon_{t+k}$,
		and denote $\ba_{k}=(a_{k,1},a_{k,2})$ for $k\in\mathbb{Z}$.
		Then, $\boldsymbol{X_t}$ is an $\alpha$-stable random vector in $\mathbb{R}^{2}$, with spectral representation $(\Gamma,\bmu^0)$ given by 
		\begin{equation}\label{def:spectral}
		\begin{aligned}
		\Gamma(A) & = \sigma^\alpha\sum_{s=\pm1}\sum_{k\in\mathbb{Z}}\dfrac{1+s\beta}{2} \|\ba_{k}\|^\alpha \delta_{\left\{\dfrac{s\ba_{k}}{\|\ba_{k}\|}\right\}}(A), \hspace{0.7cm} & \bmu^0 & = \sum_{k\in\mathbb{Z}}\ba_{k}\mu- \mathds{1}_{\{\alpha=1\}}\frac{2}{\pi}\sigma\beta\sum_{k\in\mathbb{Z}}\ba_{k}\ln\|\ba_{k}\|,
		\end{aligned}
		\end{equation}
		for any Borel set $A\subset S_2$, where $\delta_{\{\bx\}}(A)=1$ if $\bx\in A$, else $\delta_{\{\bx\}}(A)=0$, is the Dirac measure at point $\bx\in\mathbb{R}^{2}$, $\|\cdot\|$ stands for the Euclidean norm, and by convention, if for some $k\in\mathbb{Z}$, $\ba_{k}=\boldsymbol{0}$, i.e., $\|\ba_k\|=0$, then the $k$\textsuperscript{th} term vanishes from the sums.
	\end{lem}
	\textit{Proof.} See Section \ref{sec:lem_spec_rep_MA} in the Supplementary File.\\
	
	\subsection{Conditional moments}
	
	The results on bivariate stable vectors immediately apply to $\bX_t=(X_t,X_{t+h})$ with $\ba_k = (a_k,a_{k-h})$.
	A sufficient condition for the existence of conditional moments is given in the following proposition as well as their expressions.
	Without loss of generality, we will assume in the rest of this section that the stable errors have zero location parameter, i.e., $\mu=0$, unless stated otherwise.
	\begin{prop}
		\label{prop:2sidedma_condmom}
		Let $(X_t)$ be an $\alpha$-stable two-sided MA($\infty$) process, $0 < \alpha < 2$, $\beta\in[-1,1]$, $\sigma>0$, satisfying \eqref{def:2sidedma}, \eqref{def:2sidedma_cond1}-\eqref{def:2sidedma_cond2} and let $h\ge1$.\\
		$\iota)$ Assume there is $\nu\ge0$ such that
		
		\vspace*{-0.8cm}
		
		\begin{equation}
		\label{eq:nu_cond_ma}
		\sum_{k\in\mathbb{Z}}  \big(a_{k}^2+a_{k-h}^2\big)^{\frac{\alpha+\nu}{2}} |a_k|^{-\nu} < +\infty.
		\end{equation}
		Then $\mathbb{E}[|X_{t+h}|^\gamma|X_t]< +\infty$ for $0 \le \gamma < \min(\alpha+\nu,2\alpha+1)$.\\
		$\iota\iota)$ For $\alpha\ne1$, the moments $\mathbb{E}[X_{t+h}^p|X_t]$, $p\le 4$, when they exist, are given by Theorems \ref{theo:cond_expect_stable_general_alpnot1}-\ref{theo:cond_moments_stable_general_alpnot1} with
		\vspace*{-0.15cm}
		\begin{align*}
		\sigma_1^\alpha & = \sigma^\alpha \sum\limits_{k\in\mathbb{Z}}|a_k|^\alpha, & \beta_1   & = \beta \hspace{0.1cm} \dfrac{\sum\limits_{k\in\mathbb{Z}}a_k^{<\alpha>}}{\sum\limits_{k\in\mathbb{Z}}|a_k|^\alpha}, & \kappa_p    & = \dfrac{\sum\limits_{k\in\mathbb{Z}}|a_k|^{\alpha}\left(\dfrac{a_{k-h}}{a_k}\right)^p}{\sum\limits_{k\in\mathbb{Z}}|a_k|^\alpha},           & \lambda_p & = \beta \hspace{0.1cm} \dfrac{\sum\limits_{k\in\mathbb{Z}}a_k^{<\alpha>}\left(\dfrac{a_{k-h}}{a_k}\right)^p}{\sum\limits_{k\in\mathbb{Z}}|a_k|^\alpha}.
		\end{align*}
		$\iota\iota\iota)$ For $\alpha=1$, let $(\tilde{X}_t,\tilde{X}_{t+h}):=(X_t,X_{t+h})-\bmu^0$ where $\bmu^0$ is the shift vector as in Lemma \ref{prop:multistable}.
		Then, the first- and second-order moments of $\tilde{X}_{t+h}|\tilde{X}_t$ are respectively given by Theorems \ref{theo:cond_moments_stable_general_alp1_expect}-\ref{theo:cond_moments_stable_general_alp1_scm}
		with the $\kappa_p$'s, $\lambda_p$'s, $\sigma_1$, $\beta_1$ as in $\iota\iota)$ and 
		\begin{align*}
		q_0 & = \beta \hspace{0.1cm} \sum\limits_{k\in\mathbb{Z}}a_{k-h} \ln\left(\dfrac{|a_k|}{a_k^2+a_{k-h}^2}\right)/\sum\limits_{k\in\mathbb{Z}}|a_k|, & \mu_1   & = -\dfrac{2 \sigma\beta}{\pi} \hspace{0.1cm} \sum\limits_{k\in\mathbb{Z}}a_k \ln\left(\dfrac{|a_k|}{a_k^2+a_{k-h}^2}\right).
		\end{align*}
		By convention, in all the points above, if $(a_k,a_{k-h})=(0,0)$, then the $k$\textsuperscript{th} term vanishes from the sums.
	\end{prop}
	\begin{rem}[Existence of moments]
		\rm Point $\iota)$ provides a sufficient condition for the existence of conditional moments. 
		Notice that the left-hand side of \eqref{eq:nu_cond_ma} is an increasing function of $\nu$.
		Thus, if \eqref{eq:nu_cond_ma} holds for some $\nu_0>0$, it then holds for any $0\le\nu \le\nu_0$, and if it fails for $\nu_0$, it then fails for all $\nu\ge\nu_0$. 
		Causal processes, say of the form $\sum_{k\le0}a_k\varepsilon_{t+k}$ with $a_0=1$, automatically fail condition \eqref{eq:nu_cond_ma} for all $\nu>0$, 
		as $(a_h,a_{0})=(0,1)$ and the $h$\textsuperscript{th} term of the sum is finite only if $\nu=0$.
		In the case of symmetric errors ($\beta=0$), Theorem 1.1 by \cite{ct95b} allows to conclude that \eqref{eq:nu_cond_ma} is also necessary and hence that causal processes do not have finite conditional moments for orders higher than $\alpha$.
		Conversely, \eqref{eq:nu_cond_ma} may hold for some $\nu>0$ for noncausal processes 
		provided the coefficients $(a_k)$ do not decay too fast as $k\rightarrow+\infty$.
		In fact, the slower the decay of $(a_k)$ as $k\rightarrow+\infty$, the higher the values of $\nu$ for which \eqref{eq:nu_cond_ma} will hold. 
		In other terms, the stronger the dependence on <<future>> errors, the higher the order at which conditional moments will exist: hence the intuition that higher-order conditional moments may exist provided that the process is \textit{anticipative} or \textit{noncausal enough}.
		It is easy to show that \eqref{eq:nu_cond_ma} holds for any $\nu\ge0$ as soon as $(a_k)$ decays geometrically or hyperbolically, guaranteeing the existence of conditional moments up to order $2\alpha+1$ at all prediction horizons for noncausal ARMA and fractionally integrated processes. 
		Consider for instance a noncausal process $(X_t)$ of the form \eqref{def:2sidedma} such that $a_k=0$ for $k<0$, $a_k\ne0$ for $k\ge0$ and $a_k\underset{k\rightarrow+\infty}{\sim} c\lambda^k$, for some non-zero constant $c$ and $\lambda\in(-1,1)$. 
		Letting $\nu\ge0$,
		\begin{align*}
		(a_k^2+a_{k-h}^2)^{\frac{\alpha+\nu}{2}}|a_k|^{-\nu} & = |a_k|^\alpha\bigg(1+\dfrac{a_{k-h}^2}{a_{k}^2}\bigg)^{\frac{\alpha+\nu}{2}}  \underset{k\rightarrow+\infty}{\sim} |c|^\alpha (1 + \lambda^{-2h})^{\frac{\alpha+\nu}{2}} |\lambda|^{\alpha k},
		\end{align*}
		and since $|\lambda|^\alpha<1$, the summability condition \eqref{eq:nu_cond_ma} holds for any $\nu\ge0$. In particular, it holds for $\nu=\alpha+1$ and therefore, Point $\iota$ of Proposition \ref{prop:2sidedma_condmom} ensures that $(X_t)$ admits finite conditional moments up to order $2\alpha+1$.
		It is possible to find noncausal processes for which \eqref{eq:nu_cond_ma} holds only up to some $\nu\in[0,\alpha+1)$, i.e., entailing that conditional moments are finite only up to order $\gamma$ strictly within $(\alpha,2\alpha+1)$, with $\gamma$ moreover depending on the prediction horizon. 
		Such processes are necessarily noncausal and typically feature extremely short range dependence on future errors.
		See Section \ref{sec:intermediatenu} in the Supplementary File for an example.
	\end{rem}
	\begin{rem}[Computational aspects]
		\rm
		From a computational perspective, the conditional moments of $X_{t+h}$ given $X_t=x$ given in Proposition \ref{prop:2sidedma_condmom} can be inexpensively calculated for various horizons $h$ and conditioning values $x$. 
		In the case $\alpha\ne1$, computing these moments requires evaluating the functions $\mathcal{H}\big(n,\btheta;x\big)$, $n=2,3,4$, appearing in Theorem \ref{theo:cond_moments_stable_general_alpnot1}, which depend both on $x$ and on $h$ through the $\kappa_p$'s and $\lambda_p$'s given in point $\iota\iota$ of Proposition \ref{prop:2sidedma_condmom}.
		These functions can be decomposed into $a_h u_n(x) + b_h v_n(x)$, where $a_h$ and $b_h$ are constants depending only on $h$ and fixed parameters of the process, while $u_n(x)=\mathcal{H}(n,(0,1);x)$ and $v_n=\mathcal{H}(n,(1,0);x)$ are integrals of a single variable which need only to be computed once for a given conditioning value $x$. 
		Computing these integrals requires paying attention to two main hurdles. First, these are improper integrals on $(0,\infty)$, which requires truncating the integral using a high enough cutoff value $\overline{U}>0$. This will typically yield a good approximation as the integrand vanishes at exponential speed. 
		Notice that the speed of the decay does not depend on $h$ nor $x$ and a single sufficiently high threshold will do for all horizons and conditioning values.
		Second, the integrand contains an oscillatory term, whose <<frequency>> increases with $|x|$. This requires choosing a sufficiently fine subdivision of the truncated integration interval $(0,\overline{U})$. 
		For lower magnitudes of $|x|$, coarser subdivisions will suffice.
		As $|x|$ grows larger, one might fear that the required fineness of the subdivision will lead to prohibitively expensive computational costs: in this large conditioning value regime, one can however avoid the computation of the integral altogether and favour the asymptotic approximations given by Proposition \ref{prop:equiv}.
		Similar considerations hold for the moments in the case $\alpha=1$.
		More details can be found in \cite{st94} Section 5.5 on numerical techniques for computing the moment of order 1, which recommendations are still relevant for higher orders. 
	\end{rem}

	\subsection{Examples}
	
	\subsubsection{Mixed ARMA processes}
	\label{ex:marma}
	
		Mixed causal-noncausal AR (MAR) processes are often invoked in the empirical noncausal literature for speculative bubble modelling. 
		Their conditional distribution and moments are known analytically only in special cases [see \cite{fri17} for details], and, beyond these special cases, practical forecasting relies on the simulation- and sample-based methods by \cite{lls12} and \cite{gj16}.
		Mixed causal-noncausal ARMA processes with in addition a possibly non-invertible MA components (MARMA) as in \eqref{eq:marma_intro} however, have not yet taken up as much as MAR processes for speculative bubble modelling. 
		This is probably due to the absence of analytical results regarding their conditional distribution.
		Estimation procedures for such ARMA processes focus on providing estimators of the coefficients of the AR and MA polynomials, whereas the results of Proposition \ref{prop:2sidedma_condmom} rely on the coefficients of the MA($\infty$) representation $X_t=\sum_{k\in\mathbb{Z}}a_k\varepsilon_{t+k}$.
		Fortunately,  the coefficients $(a_k)$ can be recovered exactly from the AR and MA polynomials. 
		For $(X_t)$ a MARMA process solution of Equation \eqref{eq:marma_intro}, we have from \cite{gj16} Section 2.3 the following decomposition:
		\begin{equation}\label{eq:arma_decomp}
		X_t = B^p b_1(B) v_t + b_2(B) u_t,
		\end{equation}
		where $b_1(B):=\sum_{i=0}^{q-1}b_{1,i}B^i$ and $b_2(B):=\sum_{j=0}^{p-1}\,b_{2,j}B^j$ are the two polynomials resulting from the partial fraction decomposition
		$$
		\dfrac{1}{\phi(B) \Big(B^p \psi(F)\Big)} = \dfrac{b_1(B)}{\phi(B)} + \dfrac{b_2(B)}{B^p\psi(F)},
		$$
		and where $(v_t)$ and $(u_t)$ are defined by $v_t:=\psi(F)X_t$ and $u_t:=\phi(B)X_t$. 
		Letting $Z_t := \Theta(F)H(B)\varepsilon_t$, the processes $(v_t)$ and $(u_t)$ furthermore satisfy the recursions $\phi(B)v_t=Z_t$ and $\psi(F)u_t=Z_t$.
		When $\Theta=H=1$, $(X_t)$ reduces to a MAR process and $(v_t)$ and $(u_t)$ are respectively called the causal and noncausal components of $(X_t)$.
		Identifying the MA($\infty$) representations in $(\varepsilon_t)$ of the left- and right-hand side of \eqref{eq:arma_decomp} yields a general expression of the coefficients $(a_k)$ as
		\begin{equation}\label{eq:a_k_marma}
		\forall k\in\mathbb{Z}, \quad a_k = \sum_{\ell=-r}^s\vartheta_\ell \left[\sum_{i=0}^{q-1} b_{1,i}\, c_{1,k+p+i-\ell} + \sum_{j=0}^{p-1}b_{2,j}c_{2,k+j-\ell}\right],
		\end{equation}
		where $\sum_{\ell=-r}^s\vartheta_\ell F^\ell := \Theta(F)H(B)$, $(c_{1,k})$ and $(c_{2,k})$ are the coefficients of the Laurent expansions  of $1/\phi(z)$ and $1/\psi(z)$ [\cite{con78} p.107], which are such that  $c_{1,k}=0$ for $k>0$ ;  $c_{2,k}=0$ for $k<0$ ; $c_{1,0}=c_{2,0}=1$ and otherwise recursively obtained from the AR polynomials as
		\begin{align*}
		\forall k < 0, \quad c_{1,k} & =  \sum_{i=1}^q \phi_i c_{1,k+i}, & \forall k > 0, \quad c_{2,k} & =  \sum_{j=1}^p \psi_j c_{2,k-j}.
		\end{align*}
		Proposition \ref{prop:2sidedma_condmom} then applies to the MARMA process $(X_t)$ with coefficients sequence $(a_k)$ as in \eqref{eq:a_k_marma}.
		For practical purposes, the infinite sums in Proposition \ref{prop:2sidedma_condmom} can be truncated. 
		For MARMA processes,  $(a_k)$ vanishes geometrically fast as $|k|\rightarrow\infty$ and truncation will typically yield a good approximation.\\
		
\begin{figure}[!h]
			\centering
			\includegraphics[scale=0.332]{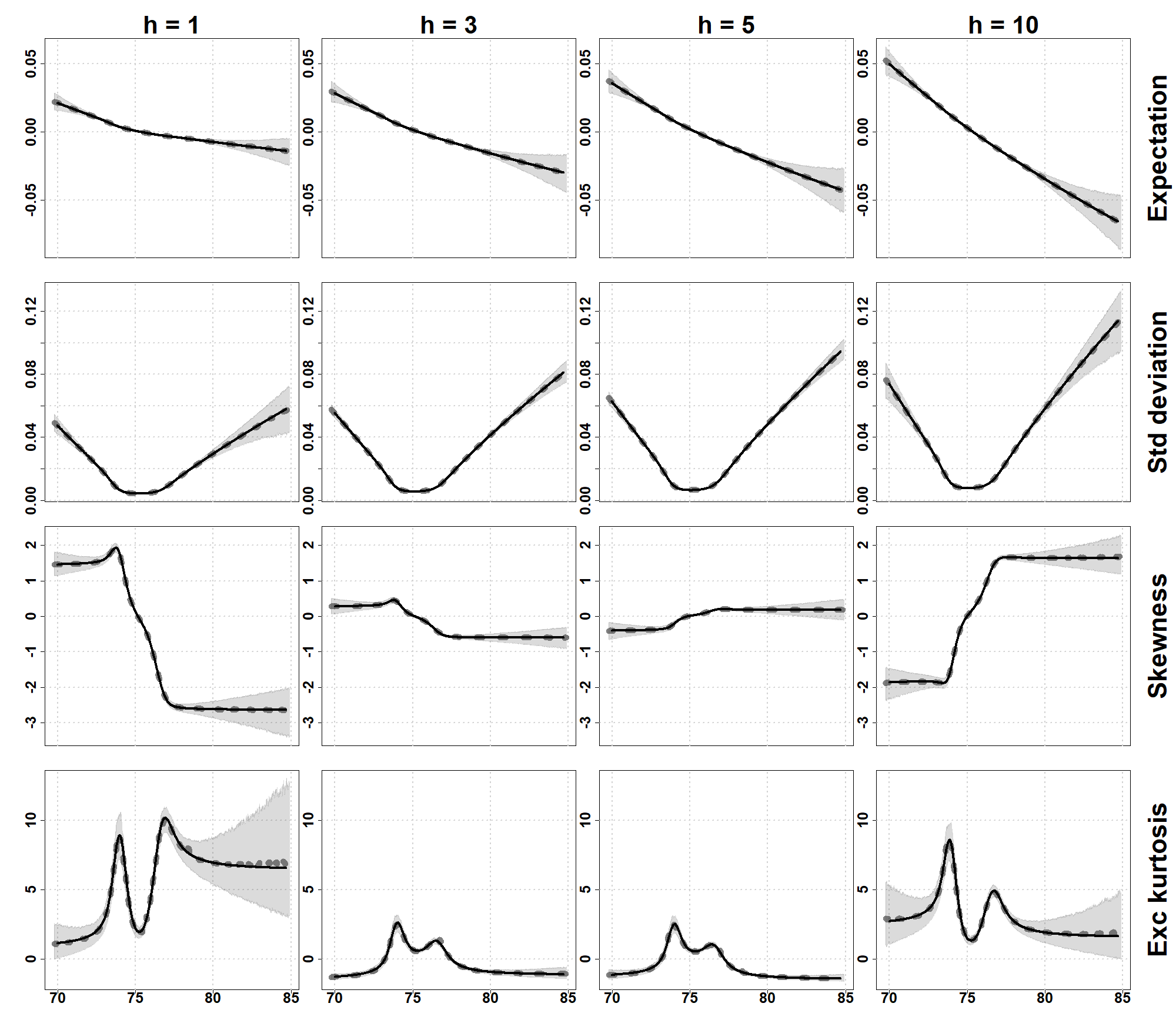}
			\caption{\footnotesize Conditional expectation, standard deviation, skewness and excess kurtosis (in rows) of the returns $R_{t+h}=(X_{t+h}-X_t)/X_t$ at horizons $h=1,3,5,10$ (in columns) of the ARMA process $(1-0.9F)(1+0.3B)X_t=(1+0.4F)(1-0.3B)\varepsilon_t$, $\varepsilon_t\stackrel{i.i.d.}{\sim}\mathcal{S}(1.8,0.5,0.2,10)$ for conditioning values $X_t=x\in(70,85)$ (x-axis of each plot, 99.9\% of the probability mass of the marginal distribution of $X_t$ is supported on (70,85)).
				Black solid lines: theoretical moments \eqref{eq:def_moments} given by Proposition \ref{prop:2sidedma_condmom} and \eqref{eq:a_k_marma}; Grey dotted lines: average of Nadaraya-Watson estimators (bandwidth=0.1) across 2000 simulated trajectories of $10^7$ observations each;
				Grey shaded areas: empirical 0.05-0.95 interquantile interval across simulations.
			}
			\label{fig:arma_moments}
		\end{figure}

		A simulation experiment was conducted to illustrate the results of Proposition \ref{prop:2sidedma_condmom} in the case of a MARMA process. 
		The theoretical conditional moments are compared to model-free non-parametrically estimated counterparts in order to assess the validity of the analytical formulae.
		Let us consider, for expository purposes, that the price series $(X_t)$ of an asset is modelled by the MARMA process defined as the strictly stationary solution of $(1-0.9F)(1+0.3B)X_t=(1+0.4F)(1-0.3B)\varepsilon_t$, $\varepsilon_t\stackrel{i.i.d.}{\sim}\mathcal{S}(1.8,0.5,0.2,10)$. 
		We will focus on the conditional moments of the returns at horizon $h$ of $(X_t)$, denoted $R_{t+h}=\dfrac{X_{t+h}-X_t}{X_t}$.
		On the one hand, we use the formulae of Proposition \ref{prop:2sidedma_condmom} to compute the theoretical expectation, standard deviation, skewness and excess kurtosis of the returns, conditional on the level $X_t=x$: 
		\begin{equation}\label{eq:def_moments}
		\begin{aligned}
		\mu(x,h) & := \mathbb{E}[R_{t+h}|X_t=x], & \sigma^2(x,h) & :=  \mathbb{E}\Big[\big(R_{t+h}-\mu(x,h)\big)^2\Big|X_t=x\Big],\\
		\gamma_1(x,h) & := \mathbb{E}\bigg[\bigg(\dfrac{R_{t+h}-\mu(x,h)}{\sigma(x,h)}\bigg)^3\bigg|X_t=x\bigg], & \gamma_2(x,h) & :=  \mathbb{E}\bigg[\bigg(\dfrac{R_{t+h}-\mu(x,h)}{\sigma(x,h)}\bigg)^4\bigg|X_t=x\bigg]-3.
		\end{aligned}
		\end{equation}
		It is just a matter of expanding the powers in the definitions above to express the conditional moments  of $R_{t+h}$ in terms of $\mathbb{E}[X_{t+h}^p|X_t]$, $p\in\{1,2,3,4\}$, where $X_{t}=\sum_{k\in\mathbb{Z}}a_k\varepsilon_{t+k}$ admits an $\alpha$-stable MA($\infty$) representation whose coefficients are given by \eqref{eq:a_k_marma}.
		On the other hand, we simulate $M=2000$ trajectories $x_1^{(m)},x_2^{(m)},\ldots,x_N^{(m)}$, $m=1,\ldots,M$, with $N=10^7$ observations of the aforementioned MARMA process and obtain model-free estimates of the conditional power moments $\mathbb{E}[X_{t+h}^p|X_t=x]$ using Nadaraya-Watson estimator
		$$
		\hat{E}^{(m)}(X_{t+h}^p|X_t=x) := \dfrac{\sum_{i=1}^{N-h} K_w\big(x-x_i^{(m)}\big)\big(x^{(m)}_{i+h}\big)^p}{\sum_{j=1}^{N-h} K_w\big(x-x_j^{(m)}\big)},
		$$
		where $K_w$ is the Gaussian kernel with bandwidth $w$. 
		Empirical counterparts $\hat{\mu}^{(m)}(x,h)$, $\hat{\sigma}^{(m)}(x,h)$, $\hat{\gamma}_{1}^{(m)}(x,h)$, $\hat{\gamma}_{2}^{(m)}(x,h)$, $m=1,\ldots,M$, of $\mu(x,h)$, $\sigma(x,h)$, $\gamma_1(x,h)$ and $\gamma_2(x,h)$ are obtained by substituting the non-parametric estimates $\hat{E}^{(m)}(X_{t+h}^p|X_t=x)$ in \eqref{eq:def_moments} in place of $\mathbb{E}[X_{t+h}^p|X_t=x]$.
		We considered prediction horizons $h=1,3,5,10$, conditioning values $x$ in the interval $(70,85)$ --corresponding to the 0.0005 and 0.9995 quantiles of the marginal distribution of $X_t$: 99.9\% of the probability mass of $X_t$ is supported on $(70,85)$-- and used a bandwidth of $w=0.1$.
		Letting $NW^{(m)}$ denote generically any of $\hat{\mu}^{(m)}$, $\hat{\sigma}^{(m)}$, $\hat{\gamma}_1^{(m)}$, $\gamma_2^{(m)}$, we compute for each quantity the point-wise average of Nadarya-Watson estimators as $\overline{NW}(x,h) := \frac{1}{M}\sum_{m=1}^M NW^{(m)}(x,h)$ as well as the point-wise 0.05 and 0.95 quantiles across simulations. 
		Figure \ref{fig:arma_moments} compares the theoretical conditional moments obtained using Proposition \ref{prop:2sidedma_condmom} and \eqref{eq:a_k_marma} with their empirical non-parametric counterparts. 
		We notice that the average $\overline{NW}$ is very closely matching the theoretical moments curves, and that the theoretical moments lie everywhere within the empirical 0.05-0.95 interquantile.
		This provides evidence for the sanity of Theorem \ref{theo:cond_moments_stable_general_alpnot1} and Proposition \ref{prop:2sidedma_condmom}. 
		In addition, we notice that the dispersion of the model-free non-parametric estimators is rather important for $X_t=x$ far from central values, despite the length of the simulated trajectories ($N=10^7$ observations).
		This suggests that the analytical formulae can hardly be traded for purely data-driven methods when it comes to estimating the dynamics during extreme events,  even with massive amounts of data.\\
		\indent To compute the conditional moments in practice, one can now overlook the model-free non-parametric approach and resort to a parametric plug-in strategy: e.g., estimate the MARMA and stable parameters by maximum likelihood and plug the parameter estimates in the formulae of Proposition \ref{prop:2sidedma_condmom}.
		An additional experiment reported in Section \ref{sec:marmaillus_supp} of the Supplementary File   illustrates the reliability of the latter parametric plug-in strategy, and its ability to accurately recover the conditional moments curves for practically-relevant sample sizes.
		Illustrations of the shape of the conditional moments for various parameterisations of MARMA processes are also provided. 
	
	\begin{rem}
	\rm The asymptotic properties of the Nadaraya-Watson for strongly mixing sequences with bounded second order marginal moment have been established by \cite{h08}.
	In our context, where higher-order conditional moments may be bounded in spite of infinite marginal variance, the validity of the Nadarya-Watson estimator is an open issue.
	The agreement between the theoretical moment curves and the empirical ones obtained with the Nadaraya-Watson estimator also suggests that the latter's validity may extend. 
	This is left for further research.
	\end{rem}
	
	\subsubsection{Cauchy MA($\infty$) processes}
	\label{ex:cauchy}
		MAR processes with Cauchy errors (stable with  $\alpha=1$ and  $\beta=0$) are a popular benchmark for speculative bubble modelling in the noncausal literature [e.g., \cite{hen15,hec16,fri17,ghj19,cav17,vh20}].
		An attractive feature of this class of models is that the Cauchy distribution is one of the special cases in the stable family for which a closed-form density is available.
		For Cauchy MAR processes with a single noncausal root, i.e., as in \eqref{eq:marma_intro} with $\Theta=H=1$, $p = 1$ and $q\ge0$, the decomposition into causal and noncausal components \eqref{eq:arma_decomp} allows to obtain the conditional moments and density in closed-form. 
		The techniques based on decomposition \eqref{eq:arma_decomp} do not extend however, and no result is available for more general Cauchy noncausal processes. \\
		\hspace*{0.5cm} Let us apply Proposition \ref{prop:2sidedma_condmom} to $X_t=\sum_{k\in\mathbb{Z}}a_k\varepsilon_{t+k}$ with $\varepsilon_t\stackrel{i.i.d.}{\sim}\mathcal{S}(1,0,\sigma,0)$ and $(a_k)$ such that Point $\iota)$ guarantees the existence of the first- and second-order moments (e.g., a Cauchy MARMA process).
		Then, invoking Point $\iota\iota\iota)$ with $\beta_1=\lambda_1=\mu_1=q_0=0$ since $\beta=0$, we have for any $x\in\mathbb{R}$ and $h\ge1$,
		
		\vspace*{-0.9cm}
		
		\begin{align*}
		\mathbb{E}[X_{t+h}|X_t=x] & = \kappa_1 x, & 
		\mathbb{V}(X_{t+h}|X_t=x) & = (\kappa_2 - \kappa_1^2)(x^2+\sigma_1^2).
		\end{align*}
		
		\vspace*{-0.2cm}
		
		\noindent In particular, if $a_k\ge0$ for all $k$, then $\kappa_1 = \sum_{k\in\mathbb{Z}} |a_k| (a_{k-h}/a_k) / \sum_{k\in\mathbb{Z}} |a_k| = 1$ and $\mathbb{E}[X_{t+h}|X_t] = X_t$.

	\begin{rem}[Conditional heteroscedasticity of noncausal processes]
	\rm 
	
	\hfill \\
	\cite{gz17} and \cite{fri17} highlighted that the Cauchy noncausal AR(1) and MAR($1,q$) processes exhibit GARCH effects in calendar time, although seemingly defined based on i.i.d. errors. 
	The above result shows that this property extends to Cauchy MA($\infty$) processes. 
	Figure \ref{fig:arma_moments} further illustrates that this is not a specific feature of the Cauchy distribution, and that modelling prices with noncausal $\alpha$-stable processes also induces conditional heteroscedasticity in the returns for other values of $\alpha$.
	In the Cauchy case, the conditional volatility is quadratic  in the past values and the authors underlined that $(X_t)$ admits a semi-strong double autoregressive representation \textit{à la} \cite{ling07}.
	The conditional first and second  moments in Proposition \ref{prop:2sidedma_condmom} suggests that a more complex representation may hold in general for $\alpha\ne1$. 
	Proposition \ref{prop:equiv} ensures nevertheless that the variance of $X_{t+h}|X_t$ is still asymptotically quadratic in the conditioning value.
	This can be noticed in the  example of the following section.
	\end{rem}

	\subsubsection{$\alpha$-stable noncausal  AR(1)}
	\label{ex:ar1}
		Let $(X_t)$ be the $\alpha$-stable noncausal  AR(1) solution of $X_t = \rho X_{t+1 } + \varepsilon_t$, $\varepsilon_t \stackrel{i.i.d.}{\sim}\mathcal{S}(\alpha,\beta,\sigma,0)$
		with $\alpha\ne1$ (for simplicity), $\beta\in[-1,1]$, $\sigma>0$ and $0<|\rho|<1$.
		Then $\mathbb{E}[|X_{t+h}|^\gamma|X_t]<+\infty$ for $0\le \gamma < 2\alpha+1$ and $h\ge1$, and the conditional moments, when they exist, are given by Proposition \ref{prop:2sidedma_condmom} with
		\begin{align*}
		\sigma_1^\alpha & = \dfrac{\sigma^\alpha}{1-|\rho|^\alpha}, & \beta_1   & = \beta \dfrac{1-|\rho|^\alpha}{1-\rho^{<\alpha>}}, & \kappa_p    & = |\rho|^{\alpha h} \rho^{-hp},           & \lambda_p & = \beta_1 \big(\rho^{<\alpha>}\big)^{ h} \rho^{-hp},
		\end{align*}
		for $p\in\{1,2,3,4\}$.
		For $\rho>0$, a clear interpretation of the distribution $X_{t+h}|X_t=x$ appears during bubble episodes, that is, as $x$ becomes large relative to the central values of process $(X_t)$. 
		Letting $\mu(x,h)$, $\sigma^2(x,h)$, $\gamma_1(x,h)$ and $\gamma_2(x,h)$ denote the conditional expectation, variance, skewness and excess kurtosis of $X_{t+h}$ given $X_t = x$ respectively (as in \eqref{eq:def_moments} with $R_{t+h}$ replaced by $X_{t+h}$), when they exist, we have
		\begin{align*}
		\mu(x,h)  & \sim (\rho^{-h}x)\rho^{\alpha h}, & \gamma_1(x,h) & \longrightarrow s\dfrac{1-2\rho^{\alpha h}}{\sqrt{\rho^{\alpha h}(1-\rho^{\alpha h})}},\\
		\sigma^2(x,h)  & \sim (\rho^{-h}x)^2\rho^{\alpha h}(1- \rho^{\alpha h}), &  \gamma_2(x,h) & \longrightarrow \dfrac{1}{\rho^{\alpha h}} + \dfrac{1}{1-\rho^{\alpha h}} - 6,
		\end{align*}
		as $\beta_1 x\rightarrow+\infty$ if $|\beta_1|=1$, $x\rightarrow\pm\infty$ if $|\beta_1|\ne1$, and $s=1$ ($s=-1$) if $x\rightarrow+\infty$ ($x\rightarrow-\infty$). See Section \ref{sec:ex_ar1} in the Supplementary File  for the proof.

	\section{Forecasting noncausal bubble crashes}
	\label{sec:crash_odds}
	
	
	\noindent For practical econometric purposes, financial bubbles in stock prices, market indexes and price-dividend ratios are typically characterised as short-lived explosive episodes followed by abrupt or gradual collapses, and are analysed using reduced form models [\cite{ps18}].  
	In this section, we focus on the dynamics of noncausal processes during such explosive episodes, that is, when the conditioning level of the trajectory takes on large positive or negative values.
	The strikingly simplistic forms of the conditional moments of the $\alpha$-stable noncausal AR(1) during such events, as given in Section \ref{ex:ar1}, are characteristic of a weighted Bernoulli distribution charging probability $\rho^{\alpha h}$ to the value $\rho^{-h}x$ and probability $1-\rho^{\alpha h}$ to $0$. 
	In the framework of this model, it is thus natural to interpret $\rho^{\alpha h}$ as the probability that the bubble survives at least $h$ more time steps, conditionally on having reached the level $X_t=x$.\\ 
	\indent Such simplification of the dynamics during extreme events is actually not limited to the $\alpha$-stable noncausal AR(1).
	 We derive here closed-form expressions of the \textit{ex ante} crash odds of bubbles generated by noncausal processes. 
	We first formally establish in the case of the noncausal AR(1) that the intuition described above holds. 
	We then show that this intuition non-trivially extends to processes featuring noncausal AR(1)-type bubbles followed by almost arbitrarily shaped collapses after the peak.
	We end this section by obtaining an expression of the crash odds in the case of noncausal MA($\infty$) processes.\\
	\indent As we  focus on the extreme events, we do not need to fully specify a parametric distribution for the errors $(\varepsilon_t)$ as in Section \ref{sec:processes}, but only require that their probability tails are similar to those of an $\alpha$-stable distribution in that they decay as power-laws.
	Formally, we assume that $(\varepsilon_t)$ is an i.i.d. error sequence with regularly varying tails:
    \begin{align}\label{eq:reg_var_errors}
    \mathbb{P}(|\varepsilon_0|>x) & = x^{-\alpha} L(x), & \dfrac{\mathbb{P}(\varepsilon_0>x)}{\mathbb{P}(|\varepsilon_0|>x) } \underset{x\rightarrow+\infty}{\longrightarrow} \dfrac{1+\beta}{2} \in[0,1],
    \end{align}
    with tail parameter $\alpha>0$, asymmetry $\beta\in[-1,1]$ and $L$ any slowly varying function at infinity, i.e., such that $L(tx)/L(x)\rightarrow 1$ as $x\rightarrow+\infty$ for all $t>0$. 
    The $\alpha$-stable distribution, with $\alpha\in(0,2)$ and asymmetry parameter $\beta$, is a typical example of distribution whose tails are power-law as in \eqref{eq:reg_var_errors}.
    However, the more general assumption above and the results in the rest of this section encompass not only noncausal processes with $\alpha$-stable errors, for which we derived the moments in the previous section, but noncausal processes with any power-law tailed errors, including (skewed) $t$-student errors often invoked in the empirical noncausal literature. Note furthermore that the tail exponent $\alpha$ (or degrees of freedom in the case of the $t$-student) is not restricted to be below 2 in this section but can take any positive value.
    
    
    
    \subsection{Crash odds of noncausal AR(1)-type bubbles}
    \label{sec:ar1_odds}
    
    \subsubsection{Purely noncausal AR(1) : exponential bubbles with instant collapses}
    
    The following proposition provides the conditional distribution of the noncausal AR(1) during explosive bubble episodes.
	
	
    
    \begin{prop}\label{prop:ar1_bernoulli}
    Let $(X_t)$ be the noncausal AR(1) process solution of
    $X_t = \rho X_{t+1} +\varepsilon_t$    with $0<\rho<1$, i.i.d. errors $(\varepsilon_t)$ satisfying \eqref{eq:reg_var_errors} for some tail exponent $\alpha>0$ and asymmetry $\beta\in[-1,1]$.
    Then, for any $h\ge1$, any $\delta\in(0,\rho^{-h})$, we have as $x\rightarrow+\infty$
    \begin{align*}
    \mathbb{P}\bigg(\dfrac{X_{t+h}}{X_t}\in[\rho^{-h}-\delta,\rho^{-h}+\delta]\bigg|sX_t>x\bigg) & \longrightarrow \rho^{\alpha h}, & 
    \mathbb{P}\bigg(\dfrac{X_{t+h}}{X_t}\in[-\delta,\delta]\bigg|sX_t >x\bigg) & \longrightarrow 1-\rho^{\alpha h},
    \end{align*}
    for any $s\in\{-1,+1\}$ if $\beta\in(-1,1)$, and $s=\beta$ if $|\beta|=1$. 
    \end{prop}
    \textit{Proof.} See Section \ref{sec:prop_ar1_bernoulli} in the Supplementary File.\\
    The proposition formalises the intuition that bubbles generated by a noncausal AR(1) with regularly varying errors feature a geometric survival distribution with probability parameter $\rho^{\alpha}$.
    This interpretation implies that the survival probability does not depend on the current scale of the bubble. 
	Surprisingly, given that the noncausal AR(1) is a Markov process, it further implies that the survival probability of bubbles does not depend at all on the past history: such bubbles display a \textit{memory-less} property.
    Several statistics of interest can be easily computed to describe their survival distribution, e.g., crash probability at horizon $h$, hazard rate, expected lifetime. 
	As the bubbles are memory-less, their survival distribution can be fully characterised by the so-called half-life, or median survival time: the duration $h_{1/2}$ such that the crash probability at horizon $h_{1/2}$ is $1/2$. 
	More generally, one can be interested in the $q$-survival quantile, $q\in[0,1]$, that is, the duration $h_q$ such that the survival probability at horizon $h_q$ is equal to $1-q$. 
	Table \ref{tab:survival_stat_ar1} summarises the expressions of these descriptive survival statistics for bubbles generated by a noncausal AR(1) model with regularly varing errors.
		\begin{table}[h]\footnotesize
		\centering
		\begin{tabular}{cccc}
			\midrule\midrule
			Crash probability at hor. $h$ & Hazard rate & Expected life & $q$-Survival quantile  \\\cmidrule(lr){1-1} \cmidrule(lr){2-2}\cmidrule(lr){3-3} \cmidrule(lr){4-4}
			$1-\rho^{\alpha h}$ & $1-\rho^{\alpha}$ & $\displaystyle\dfrac{1}{1-\rho^{\alpha}}$ & $\dfrac{\ln(1-q)}{\alpha \ln \rho}$ \\
			\midrule\midrule
		\end{tabular}
		\caption{\footnotesize Descriptive survival statistics of bubbles generated by a heavy-tailed noncausal AR(1) with  AR coefficient $\rho\in(0,1)$ and tail exponent $\alpha>0$.}
		\label{tab:survival_stat_ar1}
	\end{table}   
	Computing these statistics only requires the knowledge of the AR coefficient $\rho$ and of the tail exponent $\alpha$.
	Typically, bubbles with smaller growth rates ($\rho$ closer to unity) and driven by heavier-tailed shocks (smaller $\alpha$) are likely to last longer.
	
	On the one hand, the memory-less property of these bubbles could be appealing from a financial and economic perspective as it implies that the crash date cannot be known with certainty by traders, hence ensuring a form of no-arbitrage condition.
	Bubbles with crash dates arising according to a constant hazard rate --another feature of the geometric distribution-- appear moreover compatible with the implications of game theoretic settings where arbitrageurs attempt to time exponentially increasing bubbles and induce the crash at a random date when the selling pressure they exert is high enough [\cite{m13}].
	On the other hand, the memory-less property also implies that no sophisticated method could allow a forecaster to say anything more regarding the future of AR(1) bubbles than <<growth or crash>> with the probabilities above. 
	In the case of non-exponentially shaped bubbles or if the extreme errors driving bubbles are assumed to be endogenous rather than i.i.d. (as in \cite{bkn18}), past history could however play a more central role for prediction.
	
	\begin{rem}[Parallel with \cite{bw82}]
	\rm The dynamics of the noncausal AR(1) during bubble episodes is reminiscent of the classical model proposed by \cite{bw82}:
	\begin{equation}\label{eq:bw82}
	 	X_{t+1} = \rho^\ast c_t X_{t} + \epsilon_t, \hspace{1cm} X_0 = 0,   
	\end{equation}
	where $\rho^\ast>1$, $(\epsilon_t)$ is an i.i.d. zero-mean  and finite variance error sequence, and $(c_t)$ are i.i.d. Bernoulli distributed random variables such that $\mathbb{P}(c_t=1)=1-\mathbb{P}(c_t=0)=p\in(0,1)$.
	This model recurrently generates exponentially-shaped explosive bubbles: the trajectory follows an explosive path while $c_t=1$ and ends in a crash when $c_t=0$.
	In view of Proposition \ref{prop:ar1_bernoulli}, the bubble episodes generated by a noncausal AR(1) with regularly varying errors follow a dynamics \textit{à la} Blanchard and Watson with $\rho^{\ast}=\rho^{-1}$ and $p=\rho^{\alpha}$.
	Interestingly, while Blanchard and Watson's model is explicitly designed to feature successive bubble/bust cycles, where the bust probability is a free parameter, the noncausal AR(1) generates trajectories where bubble events intersperse calmer periods.
	The dynamics \eqref{prop:ar1_bernoulli} only emerges during bubble events and the crash probability is rather a function of the other model parameters.
    The structural constraint on the survival probability $p=\rho^\alpha$ --specific to the noncausal AR(1), a linear process shown to be suitable to describe bubble components of solutions to rational expectation price models [see \cite{gjm16}]--
    has important statistical implications.
    In the framework of Blanchard and Watson's model, statistical information about $p$ can only be gathered from the observed durations of past bubbles that have already collapsed.
    Assuming $m$ bubbles of durations $T_1,\ldots,T_m$ are observed on a given time series, say, generated by \eqref{eq:bw82}, one could propose $\hat{p}=m/\sum_{i=1}^m T_i$ as an estimator for the parameter $p$ of the Bernoulli variables $(c_t)$.
    In bubble modelling applications 
    it is however not uncommon to face very small $m$ situations, or even $m=0$ in cases where a single explosive and uncollapsed trend is observed. 
    This renders accurate estimation of $p$ difficult at best, and unfeasible at worst.
    \cite{w87} even considered $p$ not to be an identifiable parameter.
    In contrast, the estimation of $\rho^{\alpha}$ 
    can exploit more information present in the data: the sample autocorrelations of the time series and the bubble growth rates provide information about $\rho$, while the tail heaviness of the time series and of the residuals (obtained after estimation of $\rho$) provide information about $\alpha$.
    A maximum likelihood estimation of the noncausal AR(1) assuming a parametric distribution for the errors, such as $\alpha$-stable or $t$-student, would suffice to obtain an estimate of $\rho^{\alpha}$. Semi-parametric approaches could be operative as well, e.g., estimating $\rho$ by Least Squares and $\alpha$ using the Hill estimator.
	\end{rem}
	
	\subsubsection{Mixed causal-noncausal AR(1) : exponential bubbles with arbitrary collapses}
	
	To encompass explosive exponential bubble patterns followed by more complex post-peak dynamics, the noncausal literature considered adding a causal component to the noncausal AR(1), resulting in the much-invoked MAR($1,q$) processes [see for instance \cite{vh20,ghj19}].
	We show here that whatever the form of the causal component adjoined to the noncausal AR(1), i.e., whatever the shape of the collapse after the exponential growth episode, the crash probability --or more accurately, the probability of reaching the end of the exponential growth-- still follows from a geometric distribution with parameter $\rho^{\alpha}$.
	We do not restrict to the case of MAR processes but actually consider any process $(X_t)$ satisfying \eqref{def:2sidedma} with $a_k=\rho^k$ for all $k\ge0$.
	Such a process satisfies the autoregression $X_t = \rho X_{t+1} + Z_t$, where $Z_t := \sum_{\ell \le 0} b_\ell \varepsilon_{t+\ell}$, with  $b_\ell = a_\ell-\rho a_{\ell-1}$ for all $\ell \le 0$.
	Letting $m\ge1$, $h\ge1$, and $\brho:=(\rho^m,\rho^{m-1},\ldots,\rho,1)$, we will state our result in the context of a forecaster observing an ongoing explosive exponential  episode, that is, observing $(X_{t-m},\ldots,X_{t-1},X_t)$ being close to colinear with $\brho$, and wishing to forecast the future path $(X_{t+1},\ldots,X_{t+h})$.  
	The only restriction that we impose on $(a_k)_{k\le-1}$ is the one ruling out <<collapses>> that would be of similar shapes as the initial exponential growth. 
	This assumption is formalised below and the forecasting result follows.
	\begin{assumption}
	   \label{asu:not_exp_colin}
	   There is $\epsilon>0$ such that for all $k\le-1$  and $\lambda\in\mathbb{R}$,  $\big|\big|\lambda(a_{k+m},\ldots,a_{k+1},a_k) -  \brho\big|\big| > \epsilon.$
	\end{assumption}
    \begin{prop}\label{prop:mar1q_bernoulli}
    Let $m\ge1$, $h\ge1$, $\rho\in(0,1)$, $\alpha>0$ and $(X_t)$ a two-sided MA($\infty$) process with i.i.d. errors satisfying \eqref{def:2sidedma}, \eqref{def:2sidedma_cond1}-\eqref{def:2sidedma_cond2} and \eqref{eq:reg_var_errors} with $a_k = \rho^{k}$ for all $k\ge0$.
    Denote $\boldsymbol{\underline{X}}_t := (X_{t-m},\ldots,X_t)$, $\boldsymbol{\overline{X}}_{t+h} := (X_{t+1},\ldots,X_{t+h})$, $\bA_k:=(a_{k-1},\ldots,a_{k-h})/|a_k|=(\rho^{-1},\ldots,\rho^{-k},\rho^{-k}a_{-1},\ldots,\rho^{-k}a_{k-h})$ for $k\in\{0,\ldots,h\}$, $\|\cdot\|$ any norm and $d:=\min\limits_{\substack{k,\ell\in\{0,\ldots,h\}\\ k\ne\ell}} \|\bA_k-\bA_\ell\|$.
    If Assumption \ref{asu:not_exp_colin} holds for some $\epsilon>0$, we then have that $d>0$, and for any $\eta\in(0,\epsilon)$, $\delta\in(0,d)$,
    \begin{align*}
    \mathbb{P}\Bigg(\bigg|\bigg|\dfrac{\boldsymbol{\overline{X}}_{t+h}}{|X_t|}-s\bA_k\bigg|\bigg|<\delta \Bigg||X_t|>x,\bigg|\bigg|\dfrac{\boldsymbol{\underline{X}}_{t}}{|X_t|}-s\brho\bigg|\bigg|<\eta\Bigg) \underset{x\rightarrow+\infty}{\longrightarrow}  \left\{\begin{array}{cc}
         \rho^{\alpha k}(1-\rho^\alpha), & \text{ if } k\in\{0,\ldots,h-1\},  \\
         \rho^{\alpha h}, & \text{ if } k=h.
    \end{array}\right.
    \end{align*}
    for any $s\in\{-1,+1\}$ if $\beta\in(-1,1)$, and $s=\beta$ if $|\beta|=1$.
    \end{prop}
    \textit{Proof.} See Section \ref{sec:prop_mar1q_bernoulli} in the Supplementary File.\\
    The above result enjoys a very intuitive pattern interpretation. 
    We illustrate this on the example of the MAR(1,1) below.
    Let us already highlight that the odds of reaching the end of an observed exponential growth episode at some future horizon are of the same form as the crash odds of a purely noncausal AR(1), i.e., geometric governed by $\rho^{\alpha}$.
    This drastically simplifies the peak-date prediction exercise for a forecaster, who only has to estimate two parameters and can even afford to stay agnostic as to whatever form  the collapse following the peak will take.
    
\begin{figure}[h!]
        \hspace*{-0.25cm}\begin{tabular}{cc}
        \includegraphics[scale=0.8]{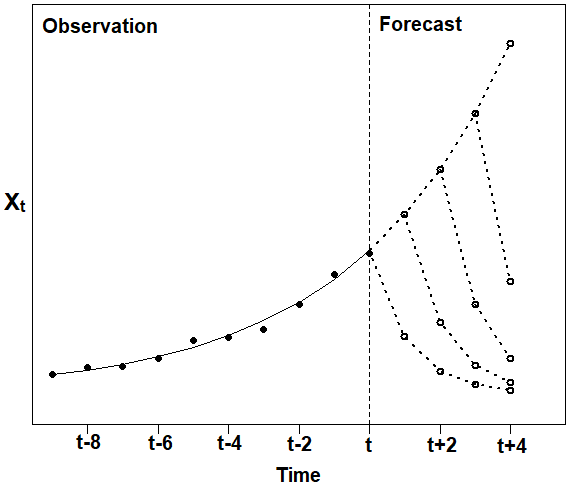} &  \begin{tikzpicture}[sloped]
       \node (root) at ( 0,0) [bbag] {$x$};
       \node (c) at ( 2.5,-1) [bag] {$\phi x$};
       \node (g) at ( 2.5,1) [bag] {$\rho^{-1}x$};
       \node (cc) at ( 5,-2) [bag] {$\phi^2x$};
       \node (gc) at ( 5,0) [bag] {$\rho^{-1}\phi x$};
       \node (gg) at ( 5,2) [bag] {$\rho^{-2}x$};
       \node (ccc) at ( 7.5,-3) [bag] {$\phi^3x$};
       \node (gcc) at ( 7.5,-1) [bag] {$\rho^{-1}\phi^2 x$};
       \node (ggc) at ( 7.5,1) [bag] {$\rho^{-2}\phi x$};
       \node (ggg) at ( 7.5,3) [bag] {$\rho^{-3}x$}; 
       \node (leg0) at ( 0,4) [bag] {$X_t$};
       \node (leg0) at ( 2.5,4) [bag] {$X_{t+1}$};
       \node (leg0) at ( 5,4) [bag] {$X_{t+2}$};
       \node (leg0) at ( 7.5,4) [bag] {$X_{t+3}$};
       \draw [->] (root) to node [below] {$1-\rho^{\alpha}$} (c);
       \draw [->] (root) to node [above] {$\rho^{\alpha}$} (g);
       \draw [->] (g) to node [above] {$\rho^{\alpha}$} (gg);
       \draw [->] (g) to node [below] {$1-\rho^{\alpha}$} (gc);
       \draw [->] (c) to node [below] {$1$} (cc);
       \draw [->] (cc) to node [below] {$1$} (ccc);
       \draw [->] (gg) to node [below] {$1-\rho^{\alpha}$} (ggc);
       \draw [->] (gg) to node [above] {$\rho^{\alpha}$} (ggg);
       \draw [->] (gc) to node [below] {$1$} (gcc);
    \end{tikzpicture}
    \end{tabular}
        \caption{\footnotesize Illustration of the likely future paths of a  bubble generated by a MAR(1,1) process with regularly varying errors as in \eqref{def:mar11}. Left panel: trajectorial interpretation with the past observed path in full points, the explosive exponential trend in solid line, and projected likely future paths in dotted lines and circles (graph drawn using $\rho=0.8$ and $\phi=0.4$). Right panel: probability tree interpretation of the projected future paths with outcomes at the origins and ends of arrows, and probabilities next to the arrows.}
        \label{fig:path_pred_illustration}
    \end{figure}

    \begin{exemple}[Forecasting  MAR(1,1) bubbles]
    \rm Consider the MAR(1,1) process defined as the strictly stationary solution of 
    \begin{equation}\label{def:mar11}
     (1-\rho F)(1-\phi B)X_t = \varepsilon_t,   
    \end{equation}
    where $|\phi|<1$ and $(\varepsilon_t)$ is an i.i.d. sequence of regularly varying errors as in \eqref{eq:reg_var_errors} with tail index $\alpha>0$.
    The process $(X_t)$ admits the MA($\infty$) representation $X_t=\sum_{k\in\mathbb{Z}} a_k \eta_{t+k}$, with $a_k=\rho^{k}$ for $k\ge0$ ; $a_k=\phi^{-k}$ for $k\le-1$ ; and $\eta_t := \varepsilon_{t}/(1-\rho\phi)$ for all $t\in\mathbb{Z}$. Note that $(\eta_t)$ is still an i.i.d. regularly varying sequence with index $\alpha$. Assumption \ref{asu:not_exp_colin} can be shown to hold and Proposition \ref{prop:mar1q_bernoulli} applies to $(X_t)$ with the sequence $(a_k)$ as described above and
    \begin{align*}
    \bA_k = \left\{\begin{array}{ll}
    (\phi,\phi^2,\ldots,\phi^{h}), & \text{ for } k=0, \\ 
    (\underbrace{\rho^{-1},\rho^{-2},\ldots,\rho^{-k}}_{k},\underbrace{\rho^{-k}\phi,\rho^{-k}\phi^2,\ldots,\rho^{-k}\phi^{h-k}}_{h-k}), & \text{ for } k\in\{1,\ldots,h-1\},\\
    (\rho^{-1},\rho^{-2},\ldots,\rho^{-h}), & \text{ for } k=h,
    \end{array}\right.
    \end{align*}
    If a forecaster observes during an extreme event of process $(X_t)$ that the recent past trajectory has approximately an exponential shape of growth rate $\rho^{-1}$, i.e., if one observes that $(X_{t-m},\ldots,X_t)$ is approximately colinear to $(\rho^m,\ldots,\rho,1)$, then the forecaster may assert that the exponential growth has probability $\rho^{\alpha h}$ to continue at least until horizon $h$, and probability $\rho^{\alpha k}(1-\rho^\alpha)$ to stop at an earlier date $k\in\{0,\ldots,h-1\}$.
    Whenever the exponential growth will reach a peak, 
    the trajectory will then enter a phase of exponential decay, with decay rate $\phi$.
    Figure \ref{fig:path_pred_illustration} illustrates the forecast interpretation from a trajectorial and probability tree perspectives.
    \end{exemple}
	\begin{rem}[Rational bubbles and fat tails]
	\rm \cite{ls02} showed that the marginal distribution of rational expectation bubble models \textit{à la} \cite{bw82} necessarily feature regularly varying tails.
	They further established that a necessary condition for any bubble process $(B_t)$ of the form \eqref{eq:bw82} to abide to the rational expectation condition $B_t = a \mathbb{E}_{t}[B_{t+1}]$, $0<a<1$, where $\mathbb{E}_{t}$ denotes the expectation conditional on all information available at date $t$, is that the tail index of the regular variations be strictly smaller than 1.
	Invoking evidence gathered by the empirical literature, which does not support such degrees of fat-tailedness, Lux and Sornette conclude that rational  bubble models  \textit{à la} Blanchard and Watson are incompatible with the observed statistical properties of financial data.\\
	\indent Interestingly, it appears that MAR processes could reconcile the rational expectations condition with tail indexes greater than 1.
	In the MAR(1,1) example above, we have that \textit{during the inflation phase of a bubble} generated by \eqref{def:mar11}, the one-step ahead conditional distribution is approximately behaved as
	$$
	X_{t+1} = \left\{\begin{array}{ll}
         \rho^{-1}X_t, & \text{ with probability } \rho^{\alpha},  \\
         \phi X_t, & \text{ with probability } 1-\rho^{\alpha}.
    \end{array}\right.
	$$
	Thus, during the inflation phase of a bubble, $\mathbb{E}_t[X_{t+1}] \approx \rho^{\alpha}(\rho^{-1}X_t) + (1-\rho^\alpha)\phi X_t$, and
	$$
	X_t \approx (\rho^{\alpha-1}+\phi - \rho^{\alpha}\phi )^{-1}\mathbb{E}_t[X_{t+1}].
	$$
	The rational expectations condition requires $(\rho^{\alpha-1}+\phi - \rho^{\alpha}\phi )^{-1}<1$, which can be rewritten as
	$$
	\alpha < \ln\bigg(\dfrac{1-\phi}{\rho^{-1}-\phi}\bigg)/\ln(\rho) := f_{\rho}(\phi).
	$$
	A straightforward analysis shows that, for any $\rho\in(0,1)$, the function $\phi\mapsto f_{\rho}(\phi)$ is strictly increasing on $[0,1)$ and that $f_\rho(0)=1$, $\lim\limits_{\phi\rightarrow1}f_\rho(\phi)=+\infty$. 
	For $\phi=0$, one retrieves Lux and Sornette's result. For $\phi>0$ however, values of $\alpha$ above 1 are admissible.
	This suggests that the MAR(1,1), as a process featuring Blanchard/Watson-like bubbles followed by gradual decays, can reconcile the rational expectations condition with regular variation tail indexes above 1.
	In fact, tail indexes arbitrarily large could be admissible provided the decay after the peak is slow enough ($\phi$ close enough to 1).
	\end{rem}
	
	\subsection{Crash odds of noncausal MA($\infty$) bubbles}
	
	Noncausal MA($\infty$) processes, which encompass general pre-peak bubble shapes, also feature a simplification of their dynamics during extreme events. 
	The following result generalises the second convergence in Proposition \ref{prop:ar1_bernoulli} to express the ex ante crash odds of bubbles generated by noncausal MA($\infty$) processes.
	\begin{prop}\label{prop:crashodds_noncausal_ma}
	Let $(X_t)$ be a MA($\infty$) process with i.i.d. errors as in \eqref{def:2sidedma}, \eqref{def:2sidedma_cond1}-\eqref{def:2sidedma_cond2} and \eqref{eq:reg_var_errors}, with $a_k = 0$ for all $k<0$ and $a_k>0$ for all $k\ge0$,  tail exponent $\alpha>0$ and asymmetry $\beta\in[-1,1]$. Assume also that there is some $\epsilon>0$ such $a_k/a_{k+1}>\epsilon$ for all $k\ge0$.
	Then, for any $h\ge1$, $\delta\in(0,\epsilon^h)$,
	\begin{equation}\label{eq:crash_odds_noncausal}
	    \mathbb{P}\bigg(\dfrac{X_{t+h}}{X_t}\in[-\delta,\delta]\bigg|sX_t>x\bigg) \underset{x\rightarrow+\infty}{\longrightarrow}  \sum\limits_{\ell=0}^{h-1}a_\ell^\alpha /\sum\limits_{k=0}^{+\infty}a_k^\alpha := p_{\infty,h},
	\end{equation}
	for any $s\in\{-1,+1\}$ if $\beta\in(-1,1)$, and $s=\beta$ if $|\beta|=1$.
	\end{prop}
	\textit{Proof.} See Section \ref{page:prop_noncausMA} in the Supplementary File.\\
	Similarly to the interpretation of the noncausal AR(1), one can notice that the crash probability of bubbles does not depend  on their current scale.
	Contrary to the noncausal AR(1) however, the survival probabilities could in general be different if the past history of the bubble was accounted for in the conditioning.
	To investigate this question, one has to characterise the conditional distribution of $X_{t+h}$ given more past information, e.g., $X_t,X_{t-1}$... This problem is out of the scope of the current paper and is addressed elsewhere [\cite{fries18}].
	To evaluate the asymptotic probability \eqref{eq:crash_odds_noncausal} in practice, only the knowledge of the coefficients $(a_k)$ and of $\alpha$ is needed, whereas asymmetry, scale or location have no role.\\ 
	\hspace*{0.5cm} We illustrate through simulations that the probability on the left-hand side of \eqref{eq:crash_odds_noncausal} indeed converges to the right-hand side limit as the conditioning value $x$ grows larger.
	We simulated $M=2000$ trajectories of $N=10^7$ observations of a noncausal AR(3) process.
	For each simulated trajectory $x_{1}^{(m)},\ldots,x_N^{(m)}$, $1\le m\le M$, we computed the following estimator of the probability \eqref{eq:crash_odds_noncausal}:
	\begin{align}\label{def:Xt+h_over_Xt_estimator}
	\hat{p}_{q,h}^{(m)} := 
	\left(\sum_{t=1}^{N-h} \mathds{1}_{\{|x_{t+h}^{(m)}/x_t^{(m)}|<\delta\} \cap \{x_t^{(m)}>q\}}\right)/\sum_{t=1}^{N-h} \mathds{1}_{\{x_t^{(m)}>q\}},
	\end{align}
	for several horizons $h$ and several quantiles $q$ of the marginal distribution of $X_t$.
	We perform this exercise twice, first assuming that the AR(3) process is driven by  $1.5$-stable errors, and then assuming $t$-student errors with 1.5 degrees of freedom.
	As our result holds for any heavy-tailed errors in the sense of \eqref{eq:reg_var_errors} and the tail exponents of the error sequences are equal, the estimated crash probabilities should tend to the same limit as $q$ increases.
	Table \ref{tab:Xt+h_overXt} gathers the average $\frac{1}{M}\sum_{m=1}^M \hat{p}_{q,h}^{(m)}$ of the empirical probabilities across the $M$ simulations along empirical 95\% confidence intervals.
	One notices that the empirical probabilities indeed come very close to the theoretical ones as $q$ increases, both for $\alpha$-stable and $t$-student errors.
	The dispersion of the non-parametric estimators across simulations again indicates that estimating the crash odds of bubble events by purely data-driven methods might be challenging, even with massive amount of data.
	The expressions given by Propositions \ref{prop:ar1_bernoulli}, \ref{prop:mar1q_bernoulli} and \ref{prop:crashodds_noncausal_ma}  thus offer the attractive alternative of computing plug-in estimators of crash odds after having estimated the model parameters. 
	\begin{table}\footnotesize
		\centering
		\begin{tabular}{clcc cc c cc c cc}
		\hline\hline
			& \hspace{0.1cm} &  & & \multicolumn{2}{c}{$h=1$} & & \multicolumn{2}{c}{$h=5$} & & \multicolumn{2}{c}{$h=10$}  \\  \cmidrule(lr){5-6} \cmidrule(lr){8-9} \cmidrule(lr){11-12} 
			&                    & & & Mean & $95\%$-CI  & &  Mean  & $95\%$-CI  & &   Mean   & $95\%$-CI  \\   \cmidrule(lr){5-5} \cmidrule(lr){6-6}  \cmidrule(lr){8-8} \cmidrule(lr){9-9} \cmidrule(lr){11-11} \cmidrule(lr){12-12}
		    \multirow{8}{*}{$\hat{p}_{q,h}$} & \multirow{2}{*}{$q_{0.9}$} & $t_{1.5}$ & \hspace*{0.25cm} & 8.37 & (8.33 , 8.41) & & 33.1 & (33.0 , 33.2) & & 42.5 & (42.3 , 42.7)   \\
		   & & \vspace*{0.2cm} $\mathcal{S}_{1.5}$ & & 9.12 & (9.08 , 9.17) & & 30.7 & (30.6 , 30.9) & & 37.6 & (37.4 , 37.7)   \\
			 & \multirow{2}{*}{$q_{0.99}$} & $t_{1.5}$ & & 17.7 & (17.5 , 17.8) & & 65.9 & (65.4 , 66.4) & & 83.4 & (83.0 , 83.8)  \\
			 & &  \vspace*{0.2cm} $\mathcal{S}_{1.5}$ & & 18.5 & (18.3 , 18.7) & & 69.1 & (68.7 , 69.6) & & 87.9 & (87.5 , 88.3)   \\
			& \multirow{2}{*}{$q_{0.999}$} & $t_{1.5}$ & & 20.5 & (19.9 , 21.1) & & 75.2 & (73.7 , 76.8) & & 94.4 & (93.3 , 95.3)   \\
			& & \vspace*{0.2cm} $\mathcal{S}_{1.5}$ & & 20.5 & (19.9 , 21.1) & & 75.6 & (74.1 , 77.1) & & 94.9 & (93.9 , 95.8)   \\
			& \multirow{2}{*}{$q_{0.9999}$} &  $t_{1.5}$  & & 20.7 & (18.9 , 22.9) & & 76.2 & (71.6 , 80.9) & & 95.4 & (92.3 , 98.0)   \\ 
			& & $\mathcal{S}_{1.5}$ & & 20.7 & (18.9 , 22.9) & & 76.2 & (71.4 , 81.1) & & 95.5 & (92.4 , 98.1)   \\ 
			\hline
			$p_{\infty,h}$ &  $\infty$  & & & 20.7 & -- & & 76.2 & -- & & 95.5 & --  \\
			\hline\hline
		\end{tabular}
		\caption{\footnotesize Comparison of  theoretical and empirical crash probabilities at horizons $h=1,5,10$ of bubbles generated by the noncausal AR(3) $X_t=0.9X_{t+1}+0.04X_{t+2}-0.096X_{t+3}+\varepsilon_t$ with 1.5-stable errors  $\varepsilon_t\stackrel{i.i.d.}{\sim}\mathcal{S}(1.5,1,0.25,0)$ ($\mathcal{S}_{1.5}$) and $t$-student errors with 1.5 degrees of freedom ($t_{1.5}$). 
		The theoretical crash probabilities are computed using \eqref{eq:crash_odds_noncausal}. Empirical average (Mean) and 95\% confidence intervals (95\%-CI) of the estimated probabilities are computed using \eqref{def:Xt+h_over_Xt_estimator} on $M=2000$ simulated trajectories of $N=10^7$ observations, with $\delta=0.2$ and for $q=q_a$ several $a$-quantiles of the marginal distribution of $X_t$. The quantiles of the marginal of $X_t$ in the case of $t$-student errors have been estimated by simulations.
		}
		\label{tab:Xt+h_overXt}
	\end{table}
	
	
	\begin{rem}[Tail dynamics and GARCH effects]
	\rm
	We here propose some intuition highlighting the connection between the tail dynamics derived in the previous propositions and the emerging GARCH effects of noncausal processes.  
	Consider for simplicity a purely noncausal process $X_t=\sum_{k\in\mathbb{Z}}a_k\varepsilon_{t+k}$ with $a_k>0$ for $k\ge0$ and $a_k=0$ for $k<0$.
The $\varepsilon_t$'s being heavy-tailed and i.i.d., if $X_t=\sum_{s\in\mathbb{Z}}a_{s-t}\varepsilon_{s}$ at some date $t$ is observed extreme, this likely results from one given $\varepsilon_{\tau}$ being extreme, for some random date $\tau$ in the neighbourhood of $t$ such that $a_{\tau-t}\ne0$, i.e., $\tau\ge t$.
Because of the i.i.d.-ness of the errors, it is likely that the extreme error $\varepsilon_\tau$ is isolated and outweights the other neighbouring $\varepsilon_s$'s contributing to $X_t$ in the sense that $\varepsilon_s/\varepsilon_\tau\approx 0$ for all $s\ne\tau$ such that $a_{s-t}\ne0$, i.e., $s\ge t$. 
Thus, we have the approximation
\begin{align*}
\dfrac{X_{t+1}}{X_t} & = \dfrac{a_{\tau-t-1} \varepsilon_\tau +\sum\limits_{s\ne \tau} a_{s-t-1}\varepsilon_{s}}{a_{\tau-t}\varepsilon_\tau +  \sum\limits_{s\ne \tau} a_{s-t}\varepsilon_{s}}  = \dfrac{a_{\tau-t-1} +\sum\limits_{s\ne \tau} a_{s-t-1}\varepsilon_{s}/\varepsilon_\tau}{a_{\tau-t} +  \sum\limits_{s\ne \tau} a_{s-t}\varepsilon_{s}/\varepsilon_\tau} \approx \dfrac{a_{\tau-t-1}}{a_{\tau-t}}.
\end{align*}
In the case of the noncausal AR(1), $a_k=\rho^{k}\mathds{1}_{k\ge0}$, and $a_{\tau-t-1}/a_{\tau-t} = \rho^{-1} \mathds{1}_{\tau \ge t+1}$ (recall that the random date $\tau$ satisfies $\mathds{1}_{\tau \ge t}=1$), which recovers the result of Proposition \ref{prop:ar1_bernoulli}: the conditional distribution of $X_{t+1}/X_t$ during extreme events concentrates on the points $0$ (crash) and $\rho^{-1}$ (growth), and the random date $\tau$ has to be interpreted as the peak date of the bubble.
Given the information at $t$, which is assumed to contain at least the value of $X_t$, the conditional variance of $X_{t+1}$ can now be approximated as 
\begin{align*}
\mathbb{V}(X_{t+1}|I_t) \approx \mathbb{V}\bigg(\dfrac{a_{\tau-t-1}}{a_{\tau-t}}X_{t}\bigg|I_t\bigg) = X_t^2 \hspace{0.2cm} \mathbb{V}\bigg(\dfrac{a_{\tau-t-1}}{a_{\tau-t}}\bigg|I_t\bigg).
\end{align*}
This analysis shows that provided the distribution of $\tau$ given $I_t$ is not degenerate (note that the existence of a non-zero constant such that $a_{k-1}/a_{k}=\text{const}$ for all $k\ge0$ is ruled out), then $\mathbb{V}\bigg(\dfrac{a_{\tau-t-1}}{a_{\tau-t}}\bigg|I_t\bigg) > 0$  and $(X_t)$ features GARCH effects during extreme events.
Continuing with the example of the noncausal AR(1), the above writes (we recognise the asymptotic variance in Section \ref{ex:ar1})
$$
\mathbb{V}(X_{t+1}|I_t) = X_t^2 \hspace{0.2cm} \mathbb{V}(\rho^{-1}\mathds{1}_{\tau \ge t+1}|I_t) = (\rho^{-1}X_t)^2 \hspace{0.2cm} \mathbb{P}(\tau \ge t+1 |I_t) \mathbb{P}(\tau < t+1 |I_t),
$$
highlighting that the conditional variance of $X_{t+1}$ stems from the uncertainty in the occurrence date of the peak given the available information.
Note that this heuristics does not necessarily presume that $I_t$ contains only information about the past values of $(X_t)$.
The set $I_t$ could contain information about other variables or noisy proxies of $\tau$ (insider information for instance).
Section \ref{ex:ar1} and Proposition \ref{prop:ar1_bernoulli} leads us to conclude that observing the infinite past of $(X_t)$ (recall that the noncausal AR(1) is Markov) does not induce $\tau|I_t$ to be degenerate and GARCH effects emerge.
Only in case of perfect foresight of $\tau$, i.e., $\tau|I_t = \text{const}$, does the GARCH phenomenon seem to vanish.
	\end{rem}

	

	\section{Evaluating the odds of crashes of real series}
	\label{sec:applications} 
	
\begin{figure}[H]
		\centering
		\includegraphics[scale=0.23]{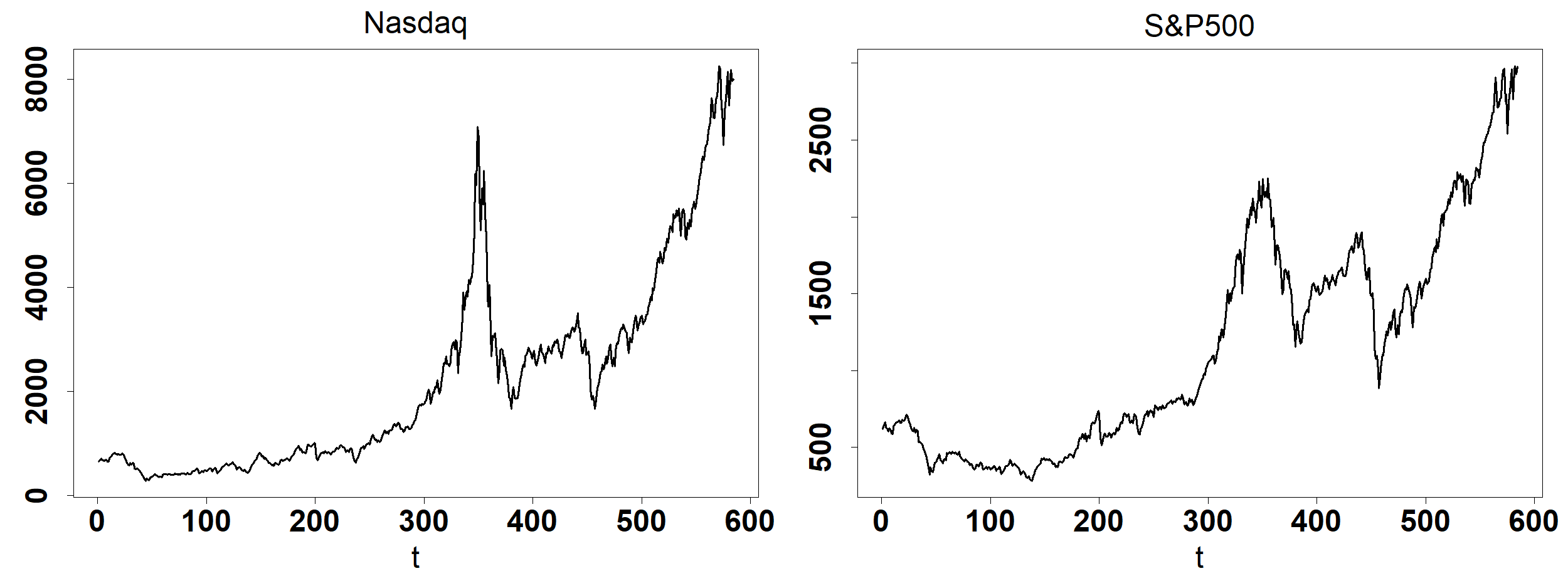}
		\caption{\footnotesize Monthly  Nasdaq  and S\&P500 real price indexes from 02/1971 to 09/2019 ($N=584$ observations each), in US dollars of 09/2019.
		}
		\label{fig:nasdaq_sp500}
	\end{figure}

	In this section, we consider two series commonly studied in the speculative bubble literature and in which evidence of explosive behaviour has been exhibited: the Nasdaq and S\&P500 indexes (see e.g. \cite{gz17,psy15,pwy11}).
	Figure \ref{fig:nasdaq_sp500} displays the monthly series of the Nasdaq and S\&P500 real prices from February 1971 to September 2019 ($N=584$ observations each), obtained by inflation-adjusting the nominal series using the Consumer Price Index 
	provided by the Federal Bank of Saint Louis (\url{fred.stlouisfed.org/series/CPIAUCSL}).
	\noindent The two series feature almost uninterrupted growth episodes since beginning of 2009 up to 2019.
	\cite{gz17} found evidence that a stable noncausal AR(1) bubble dynamics is compatible with the real Nasdaq trajectory.
	Unreported results using standard model selection methods from the noncausal literature (e.g., information criteria [\cite{lan11,hec17b}], coefficient testing [\cite{cav17}]) further confirm the reasonableness of this specification for both series of interest here.
	Starting on the premise that the explosive episodes in the data  can be modelled as ongoing realisations of noncausal AR(1) bubbles climbing towards exogenous power-law-scaled peaks, we will be in the position to provide estimates of ex ante crash odds of the recent growth trends based on the results of Section \ref{sec:crash_odds}.\\
	\indent In particular, we have shown in Section \ref{sec:ar1_odds} that bubbles generated by heavy-tailed noncausal AR(1) models feature geometric survival distributions with probability parameter $\rho^\alpha$.
	It thus suffices to provide values for the AR coefficient $\rho$ and the tail parameter $\alpha$.
	We obtain estimates of these parameters by $\alpha$-stable Maximum Likelihood, which has been shown to yield consistent estimators [\cite{and09}].
	We use the fast implementation of \cite{rsa17} to evaluate $\alpha$-stable densities, available in the R package \texttt{libstableR}.
	Contrary to $\hat{\alpha}$ which is asymptotically normal, the estimator $\hat{\rho}$ of the AR coefficient unfortunately features  an intractable asymptotic distribution.
	We resort to a parametric bootstrap procedure  to approximate the finite sample distribution of the estimators and compute confidence intervals.
	Table \ref{tab:spec_selec} reports the stable noncausal AR(1) fits and the value of the $-$log-likelihood at optima. 
	We note that the estimate of $\rho$ for the Nasdaq series is close to the one obtained by \cite{gz17}.
	The values of the $-$log-likelihood for fitted stable \textit{causal} AR(1) specifications, which estimate the AR coefficient to be 1 for both series, is provided for comparison purposes and confirm that the stationary noncausal options are to be preferred. \\
	\indent Under the stable noncausal AR(1) specification, the survival distributions of the recent explosive growth episodes can be completely estimated and characterised by plugging-in the obtained estimates $\hat{\rho}$ and $\hat{\alpha}$ in the statistics of Table \ref{tab:survival_stat_ar1}.
	For instance, an estimate of the crash probability at horizon $h$ can  be computed as $1-\hat{\rho}^{\hat{\alpha} h}$, while the $q$-life, that is, the duration $h_q$ such that the probability of an explosive episode lasting as long as $h_q$ is equal to  $1-q$, can be computed as $\ln(1-q)/\hat{\alpha}\ln \hat{\rho}$.
	Table \ref{tab:bubble_charac} displays a summary of bubble survival statistics for both series.
	In the case of the Nasdaq and S\&P500, our estimates indicate that bubbles generated by the corresponding stable noncausal  AR(1) processes should have $q=5$\% chance of lasting 8.3 and 10.6 years respectively, and $q=1$\% chance of lasting 12.7 and 16.3 years respectively.
	The observed durations of the growth episodes from 2009 to 2019 are therefore not abnormally long in that respect and appear very well compatible with the implied model properties.
	Irrespective of their past durations, as bubbles generated by such processes feature a memory-less property, this analysis suggests relatively important crash probabilities within one year ($h=12$) between 27.9 and 33.4\% for the Nasdaq series, and between 18.6 and 34.8\% for the S\&P500.
	Note that these estimates are robust to any behaviour the collapse after the peak may actually feature: as shown in Proposition \ref{prop:mar1q_bernoulli}, the shape of the collapse has no role in the ex ante probability of reaching the peak of a noncausal AR(1)-type growth episode.
	Last, similar values of the survival statistics are obtained if instead of the $\alpha$-stable assumption one opts for $t$-student or skewed-$t$ distributions, or if one proceeds to estimate $\rho$ by ordinary least squares and $\alpha$ by applying the Hill  estimator to the residuals.
	The results of the latter robustness checks are available in Section \ref{sec:robust} of the Supplementary File.
	
\begin{table}[h]\footnotesize
\centering
\begin{tabular}{cccccc}
\midrule\midrule
  &          \multicolumn{3}{c}{Stable Noncausal AR(1)}      & \hspace*{0.cm} &  Stable Causal AR(1)  \\ 
 & $\alpha$  & $\rho$ & $-$log-L & & $-$log-L \\\cmidrule(lr){2-4} \cmidrule(lr){6-6}
Nasdaq\hspace*{0.3cm} & 1.01 & 0.971 & {3608.640} & & 3642.348 \\
 & (0.925 , 1.11) & (0.969 , 0.972)  & & \\
 & & &  &  &  \\
S\&P500  & 1.36 & 0.983 & {3067.574} &  & 3086.641 \\
          & (1.25 , 1.48) & (0.975 , 0.987) & & &  \\
\midrule\midrule
\end{tabular}
\caption{\footnotesize Maximum likelihood estimates of the tail index $\alpha$ and the AR coefficient $\rho$ of stable noncausal AR(1) specifications fitted on the real Nasdaq and S\&P500 series.
The penultimate column reports the value of the $-$Log-likelihood ($-$log-L) at optima of the stable noncausal AR(1) fits (a lower value indicates a better fit).
The last column reports the value of the $-$Log-likelihood at optima of stable causal AR(1) fits as benchmark.
Finite sample 95\% confidence intervals are reported in parentheses below the point estimates and have been calculated by parametric bootstrap of the estimated stable noncausal AR(1) processes with 2000 simulated trajectories each of 584 observations.
}
\label{tab:spec_selec}
\end{table}  
	
	\begin{sidewaystable}
	\footnotesize
		\centering
		\begin{tabular}{cc}
			\begin{tabular}{ccccc}
				\midrule\midrule
				\multicolumn{5}{c}{Nasdaq}  \\
				&   &  &   &    \\
				\multicolumn{5}{c}{\begin{tabular}{cccc}
						Expected life (Y) & Half-life (Y) &  95\%--life (Y) &   99\%--life (Y) \\ \cmidrule(lr){1-1} \cmidrule(lr){2-2} \cmidrule(lr){3-3} \cmidrule(lr){4-4}
						2.8 & 1.9 & 8.3 & 12.7   \\
						(2.5 , 3.1) & (1.7 , 2.1) & (7.4 , 9.2) & (11.3 , 14.1)  \\
				\end{tabular}}  \\ 
				&   &  &   &    \\
				& \multicolumn{4}{c}{Probability of crash within $h$ months (\%)}\\
				Hazard rate & $h=1$ & 3 & 6 & 12\\ \cmidrule(lr){1-1} \cmidrule(lr){2-5}
				0.030 &  3.0 & 8.7 & 16.6 & 30.3 \\
				(0.027 , 0.033) & (2.7 , 3.3) & (7.9 , 9.7) & (15.1 , 18.4) & (27.9 , 33.4) \\
				\midrule \midrule
			\end{tabular}
			&
			\begin{tabular}{ccccc}
				\midrule \midrule
				\multicolumn{5}{c}{S\&P500}  \\
				&   &  &   &    \\
				\multicolumn{5}{c}{\begin{tabular}{cccc}
						Expected life (Y) & Half-life (Y) &  95\%--life (Y) &   99\%--life (Y)  \\ \cmidrule(lr){1-1} \cmidrule(lr){2-2} \cmidrule(lr){3-3} \cmidrule(lr){4-4}
						3.6 & 2.4 & 10.6 & 16.3   \\
						(2.4 , 4.9) & (1.6 , 3.4) & (7.0 , 15.5) & (10.7 , 22.4)  \\
				\end{tabular}}  \\ 
				&   &  &   &    \\
				& \multicolumn{4}{c}{Probability of crash within $h$ months (\%)}\\
				Hazard rate & $h=1$ & 3 & 6 & 12\\ \cmidrule(lr){1-1} \cmidrule(lr){2-5}
				0.023 & 2.3 & 6.8 & 13.2 & 24.7 \\
				(0.017 , 0.035) & (1.7 , 3.5) & (5.0 , 10.2) & (9.8 , 19.3) & (18.6 , 34.8) \\
				\midrule\midrule
			\end{tabular}
		\end{tabular}
		\caption{\footnotesize 
		Summary of bubble survival statistics using the formulae from Table \ref{tab:survival_stat_ar1}, based on the estimates of $\alpha$ and $\rho$ of stable noncausal AR(1) specifications fitted on the real Nasdaq (left panel) and S\&P500 (right panel) series (see Table \ref{tab:spec_selec} for the parameter estimates).
			Finite sample 95\% confidence intervals are reported in parentheses below the point estimates and have been calculated by parametric bootstrap of the estimated stable noncausal AR(1) processes with 2000 simulated trajectories each of 584 observations.
			Expected life, half-life, 95\%-life, 99\%-life are reported in years (Y) and crash probabilities are reported in percents (\%). }
		\label{tab:bubble_charac}
	\end{sidewaystable}
	
	\clearpage

	\section{Concluding remarks}

	By embedding $\alpha$-stable two-sided MA($\infty$) processes into the framework of bivariate $\alpha$-stable random vectors, we described in detail the conditional dependence of $X_{t+h}$ on $X_t$. 
	We have shown that noncausality plays a crucial role in the existence of conditional moments, and provided expressions for the latter up to the fourth order, when they exist, as well	as their asymptotic behaviours when the conditioning variable takes extreme values.
	We have detailed practical implementation aspects of the conditional moments as well as the contribution of the results to current methodological practices of the empirical noncausal literature.
	A future empirical investigation could determine whether corresponding patterns in the conditional moments of real data can be identified.
	These results could serve as a basis to formulate a higher-order moments portfolio allocation problem (e.g., following Jondeau and Rockinger (2006,2012)) where bubble-timing investors optimise over quantities of speculative and safer assets as well as over the holding time through a bubble. 
	Some limitations of the provided conditional moments formulae  could be addressed in further research.
	This includes expanding the conditioning to a set of past values or the entire past, as opposed to conditioning only by the present level of the trajectory.
	Also, even though Lemma \ref{prop:multistable} allows to extend the formulae of Proposition \ref{prop:2sidedma_condmom} to the conditional moments of, say, $X_{2,t+h}$ given the present level of another process $X_{1,t}$, obtaining a   characterisation of the moments in the general multivariate case remains an open issue. 
	Furthermore, despite noncausal processes admitting more conditional moments, higher-order conditional moments may nevertheless not exist for smaller values of $\alpha$ -for instance, the conditional skewness and kurtosis when $\alpha=1$.
	Alternative dependence measures capturing, say, conditional asymmetry and heavy-tailedness in such cases could be investigated.\\
	\hspace*{0.5cm} 
	Focusing on explosive bubble episodes generated by heavy-tailed noncausal MA($\infty$) processes, we provided closed-form asymptotic formulae for the predictive distribution, which enjoy very intuitive patterns and probability tree interpretations.
	This surprisingly revealed that the noncausal AR(1) bubbles are memory-less with a dynamics \textit{à la} \cite{bw82}.
	The survival distribution of such bubbles is geometric and can be fully characterised by the given of the AR coefficient $\rho$ and the tail exponent $\alpha$, both of which can be estimated by classical methods from the data. 
	Even more surprising is the fact that the augmentation of a noncausal AR(1) bubble by an arbitrarily-shaped collapse after the peak does not alter the survival distribution of the exponential growth phase of the bubble.
	From the point of view of a forecaster observing that the past trajectory is approximately exponentially-shaped, the likelihood of the peak being reached at some future horizon has the same simple expression in terms of $\rho$ and $\alpha$ whatever is bound to happen after the peak.
	Of course, the speed of the collapse still impacts how much is at risk in case of downturn. 
	Interestingly, bubbles generated by mixed causal-noncausal processes, and those of a MAR(1,1) in particular, feature an extended Blanchard and Watson dynamics with gradual collapse which appears able to reconcile rational expectation bubble models with tail exponents greater than 1, a well-documented statistical property of financial time series [\cite{ls02}].
	We further demonstrated how the closed-form formulae of the predictive distribution and of crash odds can be applied on growth episodes of real data.
	Statistical methods for agnostically estimating the coefficients $(a_k)$ of the MA representation, e.g., under low dimensional restrictions, and for robustly estimating the tail index $\alpha$ in locally explosive events could enable more refined evaluation of the crash odds.

	\section*{Acknowledgments}
	
	The author is extraordinarily indebted to Jean-Michel Zakoïan, and further thanks Denisa-Georgiana Banulescu, Jean-Marc Bardet, Frédérique Bec, Francisco Blasques, Ophélie Couperier, Gilles De Truchis, Elena Dumitrescu, Christian Francq, Christian Gouriéroux, Alain Hecq, Jérémy Leymarie, Yang Lu, Andre Lucas, Anders Rahbek,  Li Sun,   Sean Telg,  Arthur Thomas and  Elisa Voisin for insightful discussions.
	The author gratefully acknowledges the support of the Groupe des Écoles Nationales d'Économie et Statistique (GENES), the  Agence Nationale de la Recherche (via the Project	MultiRisk  ANR CE26 2016 - CR), and the support of the European Commission (via the Project NONCAUSALBubble H2020-MSCA-IF-2019 896504).

	\clearpage
	
	\addtocontents{toc}{\protect\setcounter{tocdepth}{3}}
	
	\begin{appendices}
		
		\setcounter{page}{1}
		
		\begin{center}
			{\sc \Large Supplementary Materials}\\
			{\sc \Large [For Online Publication only]}\\
		\end{center}
		
		\begin{center}
			Conditional Moments of Noncausal Alpha-Stable Processes and the Prediction of Bubble Crash Odds\\
			{\large S. Fries}
		\end{center}
		
		\tableofcontents


		
		\section{Complementary result}
		
		\subsection{Existence of moments and superexponential decay of $(a_k)$: a boundary case}
		\label{sec:intermediatenu}
		
		As pointed after Proposition \ref{prop:2sidedma_condmom}, noncausal ARMA and fractionally integrated processes whose MA coefficients decay at geometric and hyperbolic speed satisfy condition \eqref{eq:nu_cond_ma} for all $\nu>0$ (provided there are no index $k$ such that $a_{k-h}\ne0$ and $a_k=0$).
		Such processes hence admit finite conditional moments at least up to order $2\alpha+1$.
		Theorem 5.1.3 by Samorodnitsky and Taqqu, Theorems 1.1, 1.2 in \cite{ct95b} however point to the fact that intermediate cases may arise where moments are finite at most up to order $\alpha+\nu$ for some value of $\nu$ such that $\alpha < \alpha +\nu < 2\alpha + 1$.
		We propose here a noncausal MA($\infty$) process with super-exponentially decaying MA coefficients which can reach any intermediate value of the boundary.
		Consider the noncausal process defined for all $t\in\mathbb{Z}$ by $X_t = \sum_{k=0}^{+\infty} a_k \varepsilon_{t+k}$ with $a_k = \exp\{1-e^{a k}\}$, $a>0$, for all $k\ge0$, and let $(\varepsilon_{t})$ be an i.i.d. symmetrically distributed $\alpha$-stable error sequence.
		Letting $\nu\ge0$, the general term of the series in \eqref{eq:nu_cond_ma} reads for all $k\ge h$
		\begin{align*}
		(a_k^2 + a_{k-h}^2)^{\frac{\alpha+\nu}{2}}|a_k|^{-\nu} & = \big(1 + (a_{k-h}/a_k)^2\big)^{\frac{\alpha+\nu}{2}}|a_k|^{\alpha}\\
		& = \Big(1 + \exp\{2 e^{a k} (1-e^{-ah})\}\Big)^{\frac{\alpha+\nu}{2}} \exp\{-\alpha(1- e^{a k })\}\\
		& \underset{k \rightarrow +\infty}{\sim}  \exp\Big\{e^{a k} \big[ (1-e^{-ah}) (\alpha+\nu)- \alpha\big] + \alpha\Big\},
		\end{align*}
		which is the term of an absolutely convergent series if and only if $(1-e^{-ah}) (\alpha+\nu)- \alpha<0$, hence if and only if
		\begin{equation}\label{eq:superexpnucond}
		\nu < \alpha \left(\dfrac{1}{1-e^{-a h}} - 1\right).  
		\end{equation}
		Because we assume $(\varepsilon_t)$ to be symmetrically distributed, Theorems 1.1 and 1.2 in \cite{ct95b} allow to consider \eqref{eq:nu_cond_ma} and \eqref{eq:superexpnucond} as sufficient and necessary conditions for the finiteness of $\mathbb{E}[|X_{t+h}|^{\gamma}|X_t]$, $0 \le \gamma < \min(\alpha+\nu,2\alpha+1)$, in most configurations of $\alpha$ and $\nu$ (see within \cite{ct95b} for details).
		In particular, one can see that for a fixed prediction horizon $h\ge1$, the upper bound \eqref{eq:superexpnucond} on $\nu$ can lie anywhere between 0 and $+\infty$ according to the parameter $a$.
		The smaller $a>0$, i.e., the slower the decay, the higher the bound on $\nu$, and conversely, the greater $a$ (faster decay), the smaller the upper bound on $\nu$ for the existence of conditional moments.\\
		Furthermore, contrary to the case where $(a_k)$ decays at geometric or hyperbolic speeds, the finiteness of $\mathbb{E}[|X_{t+h}|^{\gamma}|X_t]$ also depends on the prediction horizon $h$.
		Most notably, for any fixed decay speed $a$, on can see that the bound \eqref{eq:superexpnucond} tends to 0 as $h\rightarrow+\infty$.
		For a decay parameter $a$ small enough, the moments $\mathbb{E}[|X_{t+h}|^{\gamma}|X_t]$ may thus be finite up to order $2\alpha+1$ for short-term prediction horizons while being finite only up to order $\alpha$ for longer-term prediction horizons.
		
		\newpage
		
		 \subsection{ Moments of MARMA processes : Complementary simulations and illustrations}
		 \label{sec:marmaillus_supp}
		 
		 \subsubsection{Plug-in estimation of the conditional moments}

\begin{figure}[H]
		     \centering
		     \includegraphics[scale=0.5]{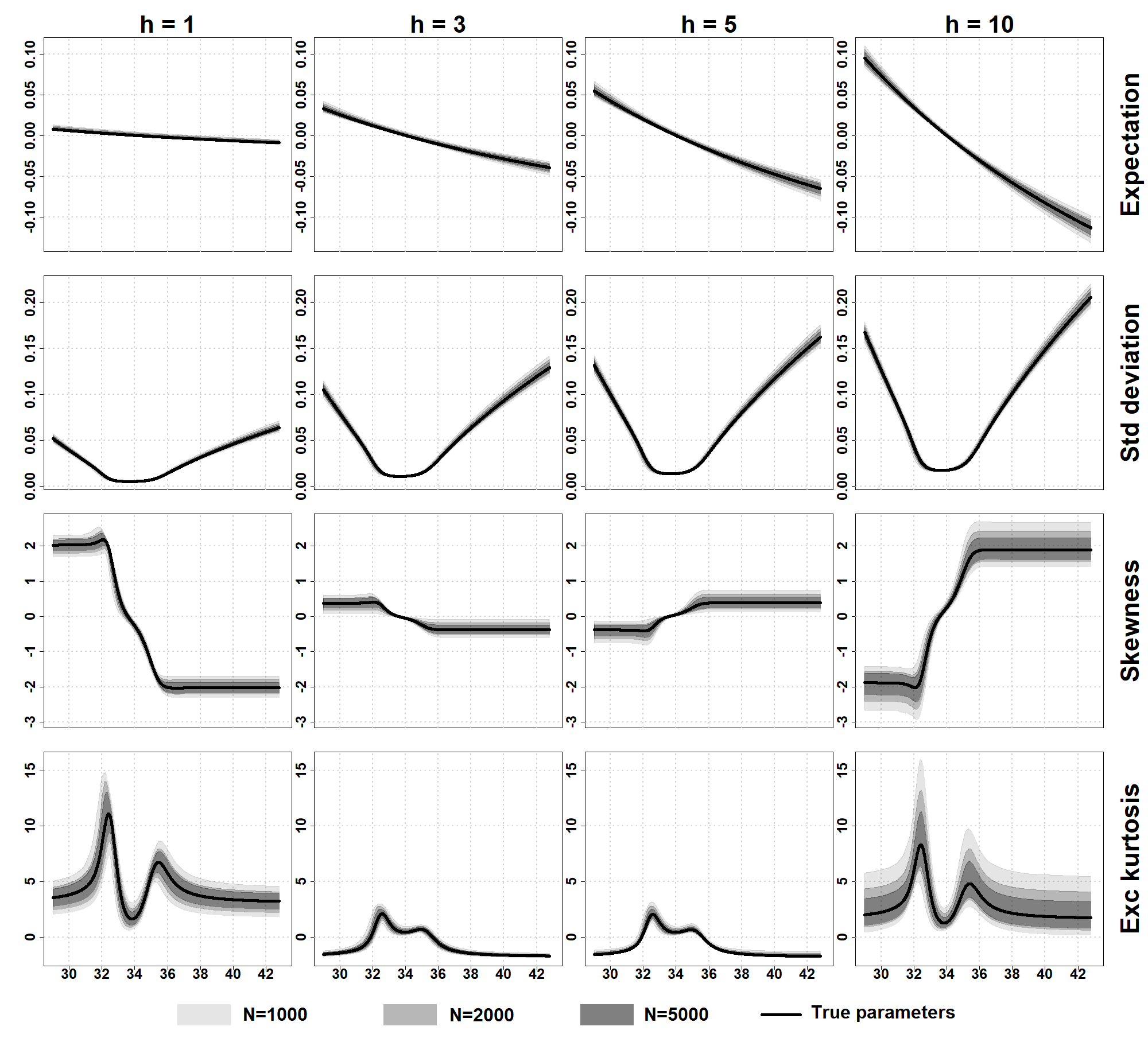}
		     \caption{Conditional expectation, standard deviation, skewness and excess kurtosis (in rows) of the returns $R_{t+h}=(X_{t+h}-X_t)/X_t$ at horizons $h=1,3,5,10$ (in columns) of the MARMA process $
		    X_t = 0.9X_{t+1} + \varepsilon_t+0.7\varepsilon_{t+1}, \hspace{0.2cm} \varepsilon_t\stackrel{i.i.d.}{\sim}\mathcal{S}(1.8,0.5,0.1,2)$, for conditioning values $X_t=x\in(29,43)$ (x-axis of each plot). Black solid lines: conditional moments \eqref{eq:def_moments} given by Proposition \ref{prop:2sidedma_condmom} computed using the true parameter values.
		    Grey shaded areas: 0.05-0.95 interquantile intervals across $M=2000$ simulations of the conditional moments \eqref{eq:def_moments} obtained by estimating the parameters of $(X_t)$ by maximum likehood and  plugging-in the estimates into the formulae of Proposition \ref{prop:2sidedma_condmom}. Estimation performed with sample sizes $N=1000$ (light grey), $N=2000$ (middle grey), $N=5000$ (dark grey).}
		     \label{fig:CondMom_ARMA_ML}
		 \end{figure}

\begin{figure}[H]
		     \centering
		     \includegraphics[scale=0.335]{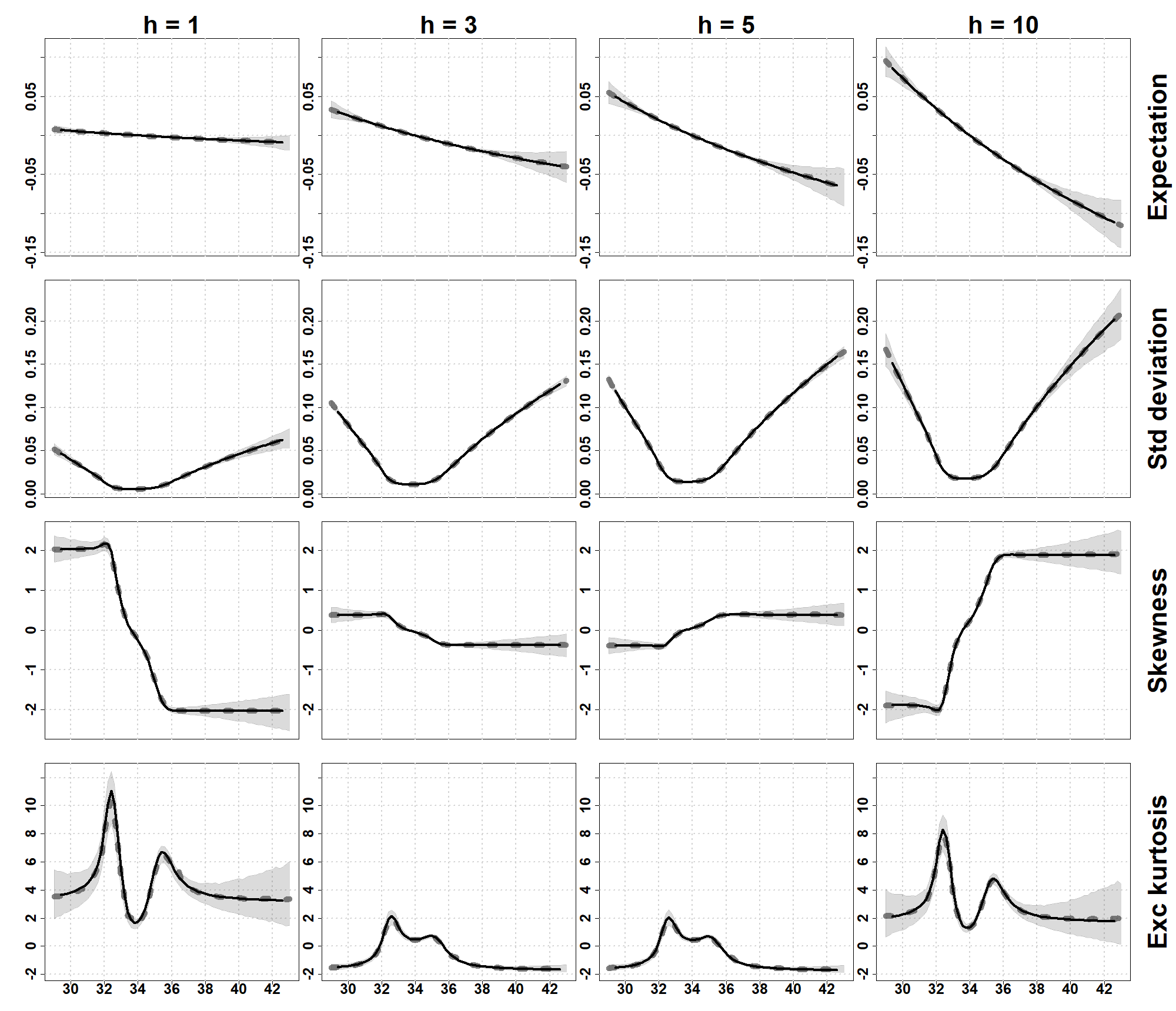}
		     \caption{Conditional expectation, standard deviation, skewness and excess kurtosis (in rows) of the returns $R_{t+h}=(X_{t+h}-X_t)/X_t$ at horizons $h=1,3,5,10$ (in columns) of the MARMA process $
		    X_t = 0.9X_{t+1} + \varepsilon_t+0.7\varepsilon_{t+1}, \hspace{0.2cm} \varepsilon_t\stackrel{i.i.d.}{\sim}\mathcal{S}(1.8,0.5,0.1,2)$, for conditioning values $X_t=x\in(29,43)$ (x-axis of each plot). 
		    Black solid lines: theoretical moments \eqref{eq:def_moments} given by Proposition \ref{prop:2sidedma_condmom}. Grey dotted lines: average of Nadaraya-Watson estimators (bandwidth=0.1) across 2000 simulated trajectories of $10^7$ observations each.
		    Grey shaded areas: 0.05-0.95 interquantile interval across simulations.}
		     \label{fig:CondMom_ARMA_nonpar}
		 \end{figure}

		 The simulation experiment of Section \ref{ex:marma} and Figure \ref{fig:arma_moments} illustrated the validity of the formulae of Proposition \ref{prop:2sidedma_condmom} by showing their match with model-free, data-driven non-parametric estimates of the conditional moments based on simulated trajectories.
		 To compute the conditional moments in practice, one can now overlook the model-free nonparametric approach and resort to a parametric plug-in estimation approach based on the formulae of Proposition \ref{prop:2sidedma_condmom} as follows:
		 \begin{enumerate}
		     \item Estimate the parameters of the stable MARMA process, for instance by $\alpha$-stable Maximum Likelihood [\cite{and09}] or M-estimation [\cite{wu13}].
		     \item Plug the parameter estimates in the formulae of Proposition \ref{prop:2sidedma_condmom} and compute the conditional moments.
		 \end{enumerate}
		 Provided the estimators of step 1 are consistent, as is the case for the two mentioned above, the plug-in estimators of the conditional moments will also be consistent.
		 We provide here the methodology and results of an additional simulation experiment designed to illustrate and gauge the parametric plug-in estimation of the conditional moments.
		 Consider again that the price $(X_t)$ of an asset is modelled by a MARMA process, say the noncausal-noninvertible solution of 
		 \begin{equation}\label{eq:condmom_arma_ml}
		    X_t = \psi_0X_{t+1} + \varepsilon_t+\theta_0\varepsilon_{t+1}, \hspace{1cm} \varepsilon_t\stackrel{i.i.d.}{\sim}\mathcal{S}(\alpha_0,\beta_0,\sigma_0,\mu_0), 
		 \end{equation}
		 where $\bvartheta_0:=(\psi_0,\theta_0,\alpha_0,\beta_0,\sigma_0,\mu_0)=(0.9,0.7, 1.8, 0.5, 0.1,2)$ is the vector of true parameter values.
		 We simulate $M=2000$ trajectories $x_1^{(m)},\ldots,x_N^{(m)}$, $m=1,\ldots,M$ from the above process for sample sizes $N=1000, 2000$ and $5000$ and estimate all the parameters by maximum likelihood as follows.\\
		 For any candidate vector of parameters $\bvartheta=(\psi,\theta,\alpha,\beta,\sigma,\mu)$, we follow \cite{wu13} to compute the residuals and evaluate the likelihood. 
		 To fix ideas, let us focus on simulation $m$. For any given candidate vector $\bvartheta$, we compute the residuals $z_1^{(m)}(\bvartheta),\ldots,z_N^{(m)}(\bvartheta)$ as: first compute $v_1^{(m)}(\bvartheta),\ldots,v_{N}^{(m)}(\bvartheta)$   according to
		 $$
		 v_t^{(m)}(\bvartheta) = x_t^{(m)} - \psi x_{t+1}^{(m)},
		 $$
		  for $t=1,\ldots,N$. Compute then the residuals backwards as
		 $$
		 z_{t-1}^{(m)}(\bvartheta) = \theta z_{t}^{(m)}(\bvartheta) - \psi v_t^{(m)}(\bvartheta),
		 $$
		 $t=N,N-1,\ldots,2$. In both steps, one can set initial (terminal) conditions, e.g., $z_{N}^{(m)}(\bvartheta)=x_{N+1}^{(m)}:=0$, and burn residuals close to the boundary.
		 Based on residuals $z_1^{(m)}(\bvartheta),\ldots,z_N^{(m)}(\bvartheta)$ (possibly accounting for some burn), we can then compute the $-$log-likelihood
		 $$
		 \mathcal{L}\Big(\bvartheta=(\psi,\theta,\alpha,\beta,\sigma,\mu);x_1^{(m)},\ldots,x_N^{(m)}\Big) := -\sum_{t=1}^N \ln f_{\alpha,\beta,\sigma,\mu}\Big(z_{t}^{(m)}(\bvartheta)\Big),
		 $$
		 where $f_{\alpha,\beta,\sigma,\mu}$ denotes the density of the stable distribution with corresponding parameters. 
		 For each simulated trajectory $m=1,\ldots,M$ and each sample size $N=1000, 2000, 5000$, we compute the residuals and the likelihood by following the steps above, and we numerically find the vector of parameters $\hat{\bvartheta}^{(m,N)}_{\text{ML}}$ minimising the $-$log-likelihood:
		 $$
		 \hat{\bvartheta}^{(m,N)}_{\text{ML}} := \argmin_{\bvartheta\in\mathbb{R}^6} \mathcal{L}\Big(\bvartheta;x_1^{(m)},\ldots,x_N^{(m)}\Big).
		 $$
		 The obtained estimators $\hat{\bvartheta}^{(1,N)}_{\text{ML}},\ldots,\hat{\bvartheta}^{(M,N)}_{\text{ML}}$ are then plugged in the formulae of Proposition \ref{prop:2sidedma_condmom} to compute the plug-in estimators $\hat{\mu}^{(m,N)}_{\text{ML}}(x,h)$, $\hat{\sigma}^{(m,N)}_{\text{ML}}(x,h)$, $\hat{\gamma}^{(m,N)}_{1,\text{ML}}(x,h)$, $\hat{\gamma}^{(m,N)}_{2,\text{ML}}(x,h)$, $m=1,\ldots,M$, $N=1000, 2000, 5000$, of the conditional expectation, standard deviation, skewness and excess kurtosis given in \eqref{eq:def_moments}.
		 Figure \ref{fig:CondMom_ARMA_ML} represents the pointwise 0.05-0.95 interquantile intervals of the plug-in conditional moments estimators across the $M=2000$ simulations, alongside the conditional moments computed using the true values of the parameters $\bvartheta_0$.
		 Three interquantile intervals appear on Figure \ref{fig:CondMom_ARMA_ML}, one for each sample size $N=1000, 2000, 5000$. 
		 It can be noticed that even for the smallest sample size, the interquantile interval is extremely narrow around most of the true conditional moments curves.
		 For higher-order moments and at furthest horizons, the interquantile intervals are slightly larger for $N=1000$ but narrow down fast as the sample size increases.
		 For comparison purposes, the model-free non-parametric estimators of the conditional moments for model \eqref{eq:condmom_arma_ml} have been computed following the same procedure as in Section \ref{ex:marma} and Figure \ref{fig:arma_moments}, using the Nadaraya-Watson estimator.
		 The non-parametric estimators have been computed based on $M=2000$ simulated trajectories of $N=10^7$ observations.
		 Figure \ref{fig:CondMom_ARMA_nonpar} represents the 0.05-0.95 interquantile intervals of the non-parametric estimator  alongside the true conditional moments.
		 When comparing Figures \ref{fig:CondMom_ARMA_ML} and \ref{fig:CondMom_ARMA_nonpar}, one can notice the dramatic efficiency gain from using the parametric plug-in estimators compared to the model-free non-parametric approach. 
		 With sample sizes four orders of magnitudes smaller, the parametric plug-in approach achieves a comparable or better accuracy.
		 Importantly, the plug-in approach is able to extrapolate well the conditional moments for conditioning values $X_t=x$ far away from the central values of the process $(X_t)$, while the error of the model-free non-parametric approach explodes in these regions because of the scarcity of extreme-valued data points.
		 
		 \clearpage
		 
		 \subsubsection{Illustrating the effects of parameters on the conditional moments}
		
\begin{figure}[H]
		    \centering
		    \includegraphics[scale=0.25]{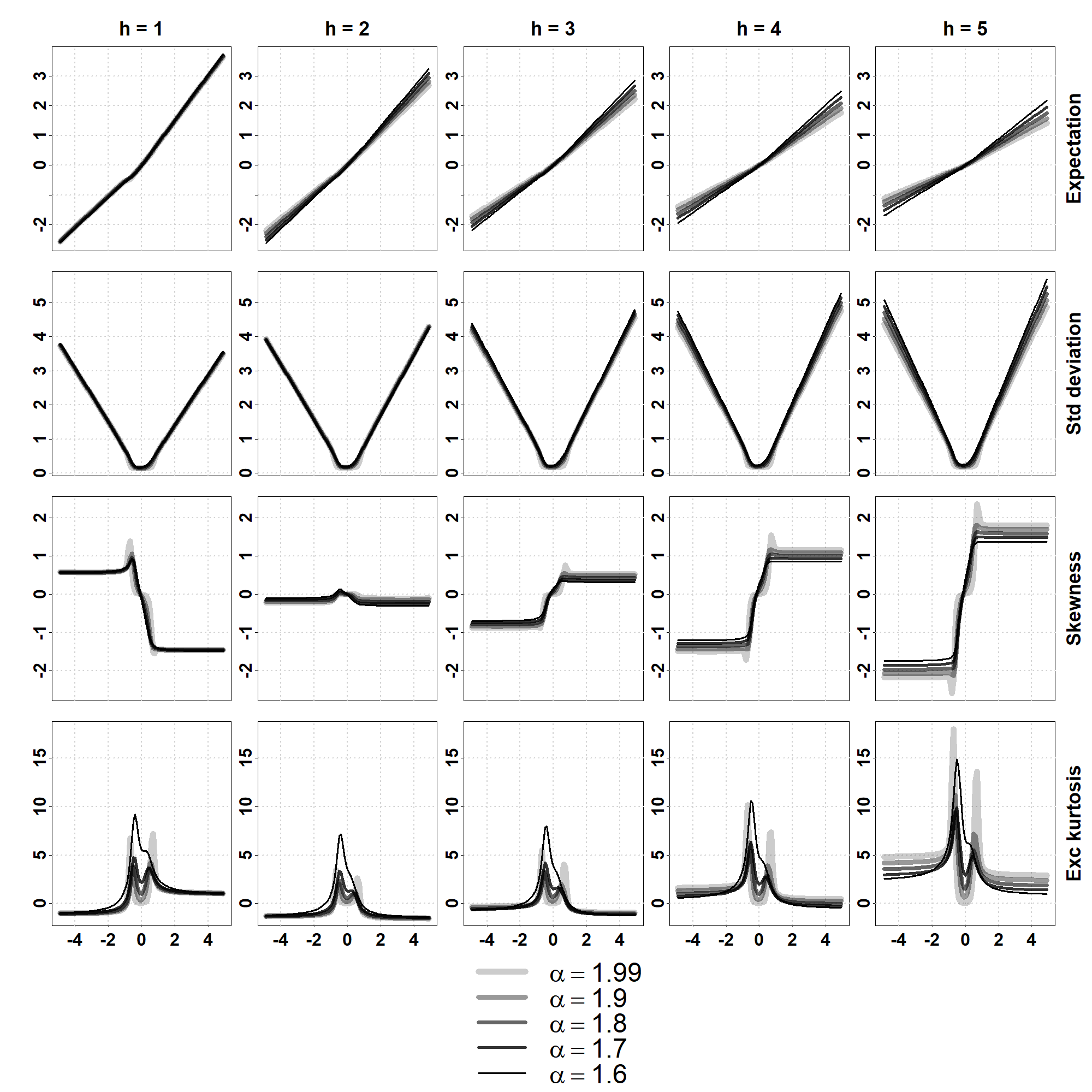}
		    \caption{\textbf{Conditional moments of a stable MARMA(1,1,1,1)  for different values of $\boldsymbol{\alpha}$.}\\ Conditional expectation, standard deviation, skewness and excess kurtosis (in rows) of $X_{t+h}$ given $X_t=x$, for horizons $h=1,2,3,4,5$ (in columns) and conditioning values $X_t=x\in(-5,5)$ (x-axis of each plot), computed using the formulae of Proposition \ref{prop:2sidedma_condmom}, where $(X_t)$ is the strictly stationary solution of  $(1-0.8F)(1+0.3B)X_t=(1+0.4F)(1-0.3B)\varepsilon_t$, $\varepsilon_t\stackrel{i.i.d.}{\sim}\mathcal{S}(\alpha,0.5,0.1,0)$, $\alpha\in\{1.99,1.9,1.8,1.7,1.6\}$.}
		    \label{fig:condmom_various_alp}
		\end{figure}

\begin{figure}[H]
		    \centering
		    \includegraphics[scale=0.25]{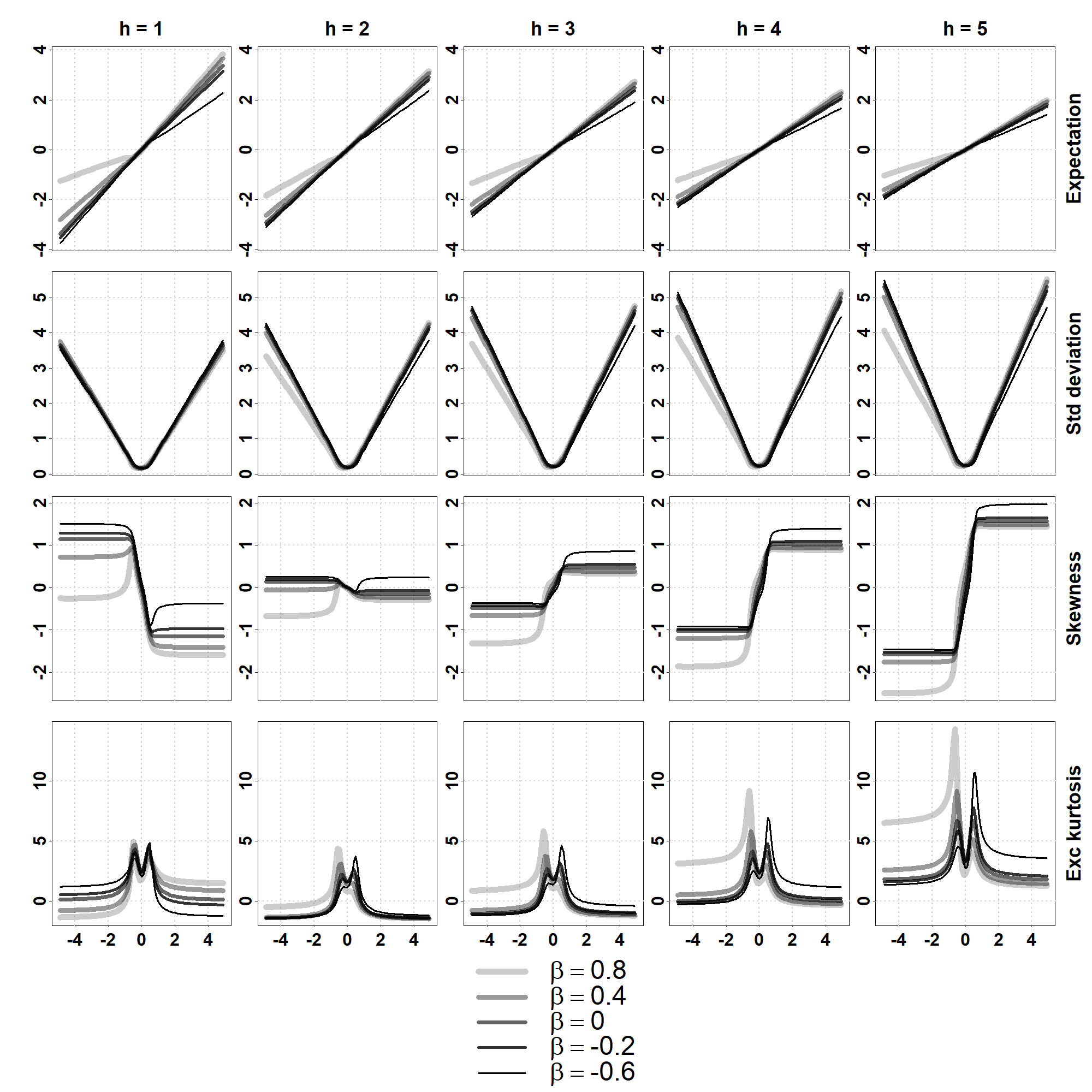}
		    \caption{\textbf{Conditional moments of a stable MARMA(1,1,1,1)  for different values of $\boldsymbol{\beta}$.}\\ Conditional expectation, standard deviation, skewness and excess kurtosis (in rows) of $X_{t+h}$ given $X_t=x$, for horizons $h=1,2,3,4,5$ (in columns) and conditioning values $X_t=x\in(-5,5)$ (x-axis of each plot), computed using the formulae of Proposition \ref{prop:2sidedma_condmom}, where $(X_t)$ is the strictly stationary solution of  $(1-0.8F)(1+0.3B)X_t=(1+0.4F)(1-0.3B)\varepsilon_t$, $\varepsilon_t\stackrel{i.i.d.}{\sim}\mathcal{S}(1.7,\beta,0.1,0)$, $\beta\in\{0.8,0.4,0,-0.2,-0.6\}$.}
		    \label{fig:condmom_various_bet}
		\end{figure}

\begin{figure}[H]
		    \centering
		    \includegraphics[scale=0.25]{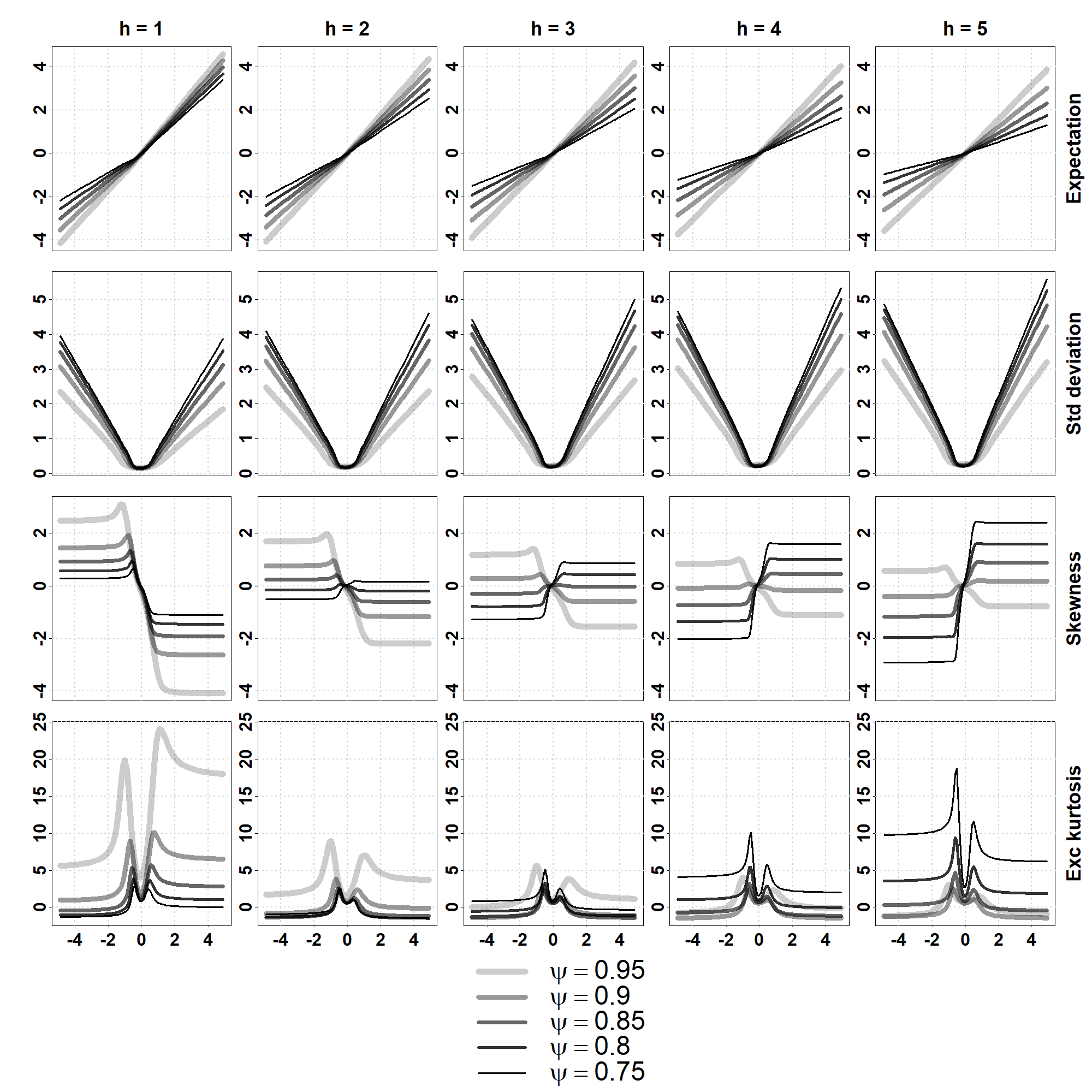}
		    \caption{\textbf{Conditional moments of a stable MARMA(1,1,1,1)  for different values of $\boldsymbol{\psi}$.}\\ Conditional expectation, standard deviation, skewness and excess kurtosis (in rows) of $X_{t+h}$ given $X_t=x$, for horizons $h=1,2,3,4,5$ (in columns) and conditioning values $X_t=x\in(-5,5)$ (x-axis of each plot), computed using the formulae of Proposition \ref{prop:2sidedma_condmom}, where $(X_t)$ is the strictly stationary solution of  $(1-\psi F)(1+0.3B)X_t=(1+0.4F)(1-0.3B)\varepsilon_t$, $\varepsilon_t\stackrel{i.i.d.}{\sim}\mathcal{S}(1.7,0.5,0.1,0)$, $\psi\in\{0.95,0.9,0.85,0.8,0.75\}$.}
		    \label{fig:condmom_various_psi}
		\end{figure}

\begin{figure}[H]
		    \centering
		    \includegraphics[scale=0.25]{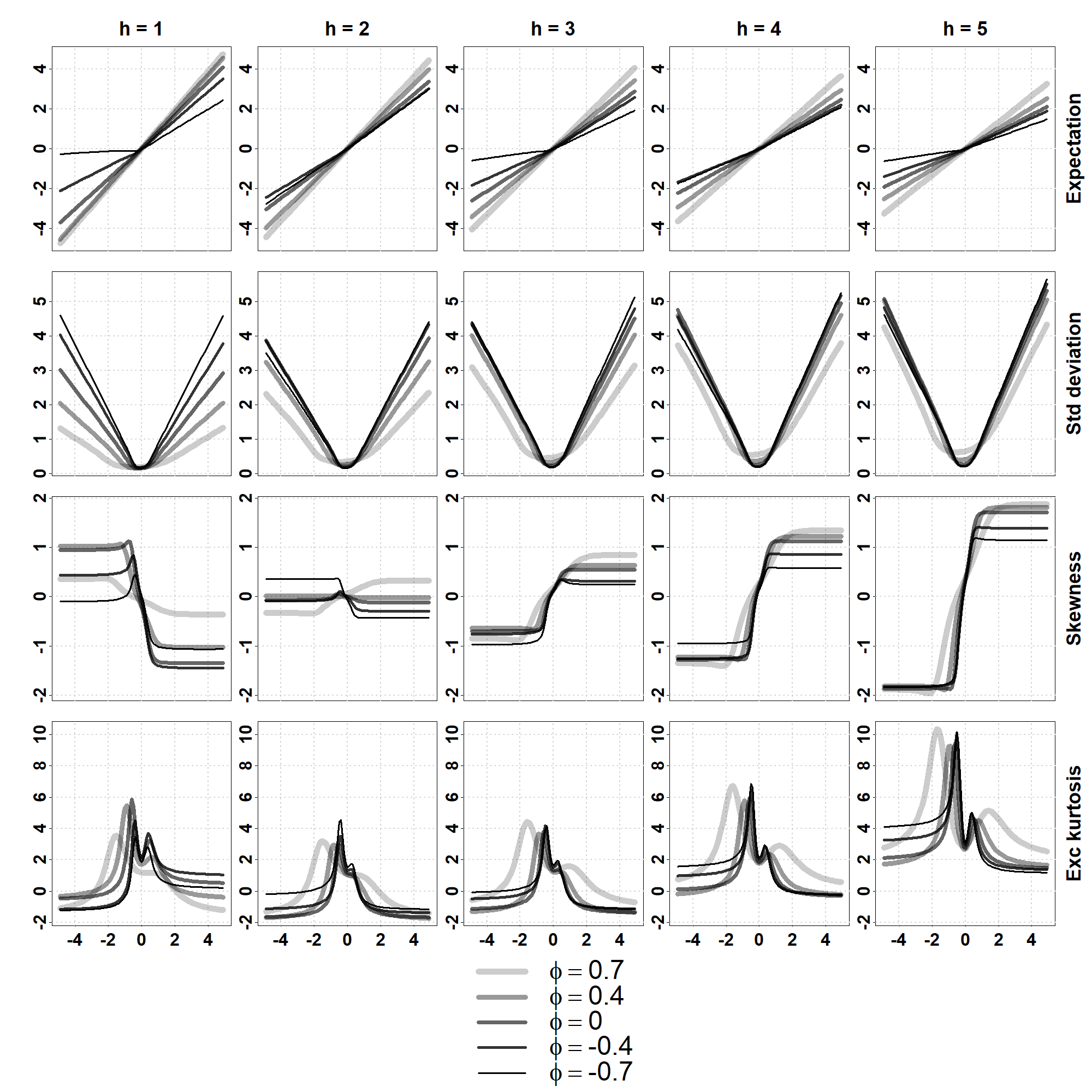}
		    \caption{\textbf{Conditional moments of a stable MARMA(1,1,1,1)  for different values of $\boldsymbol{\phi}$.}\\ Conditional expectation, standard deviation, skewness and excess kurtosis (in rows) of $X_{t+h}$ given $X_t=x$, for horizons $h=1,2,3,4,5$ (in columns) and conditioning values $X_t=x\in(-5,5)$ (x-axis of each plot), computed using the formulae of Proposition \ref{prop:2sidedma_condmom}, where $(X_t)$ is the strictly stationary solution of  $(1-0.8 F)(1-\phi B)X_t=(1+0.4F)(1-0.3B)\varepsilon_t$, $\varepsilon_t\stackrel{i.i.d.}{\sim}\mathcal{S}(1.7,0.5,0.1,0)$, $\phi\in\{0.7,0.4,0,-0.4,-0.7\}$.}
		    \label{fig:condmom_various_phi}
		\end{figure}

\begin{figure}[H]
		    \centering
		    \includegraphics[scale=0.5]{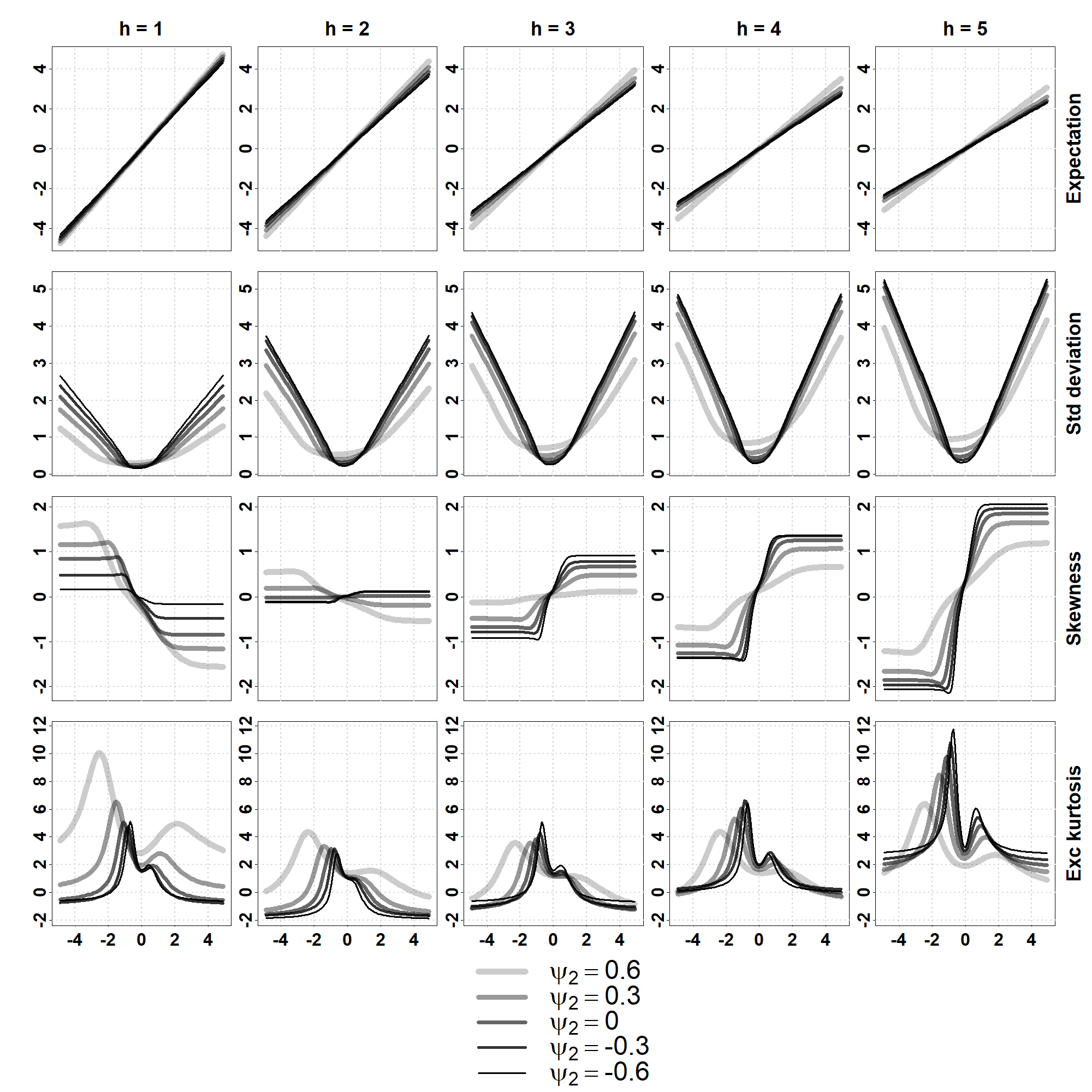}
		    \caption{\textbf{Conditional moments of a stable MAR(2,0,0,1) for different values of $\psi_2$}\\ Conditional expectation, standard deviation, skewness and excess kurtosis (in rows) of $X_{t+h}$ given $X_t=x$, for horizons $h=1,2,3,4,5$ (in columns) and conditioning values $X_t=x\in(-5,5)$ (x-axis of each plot), computed using the formulae of Proposition \ref{prop:2sidedma_condmom}, where $(X_t)$ is the strictly stationary solution of  $(1-\psi_1 F)(1-\psi_2 F)X_t=(1+0.9 B)\varepsilon_t$, $\varepsilon_t\stackrel{i.i.d.}{\sim}\mathcal{S}(1.7,0.5,0.1,0)$, $\psi_1=0.8$, $\psi_2\in\{0.6,0.3,0,-0.3,-0.6\}$.}
		    \label{fig:condmom_various_lambda1}
		\end{figure}

\begin{figure}[H]
		    \centering
		    \includegraphics[scale=0.5]{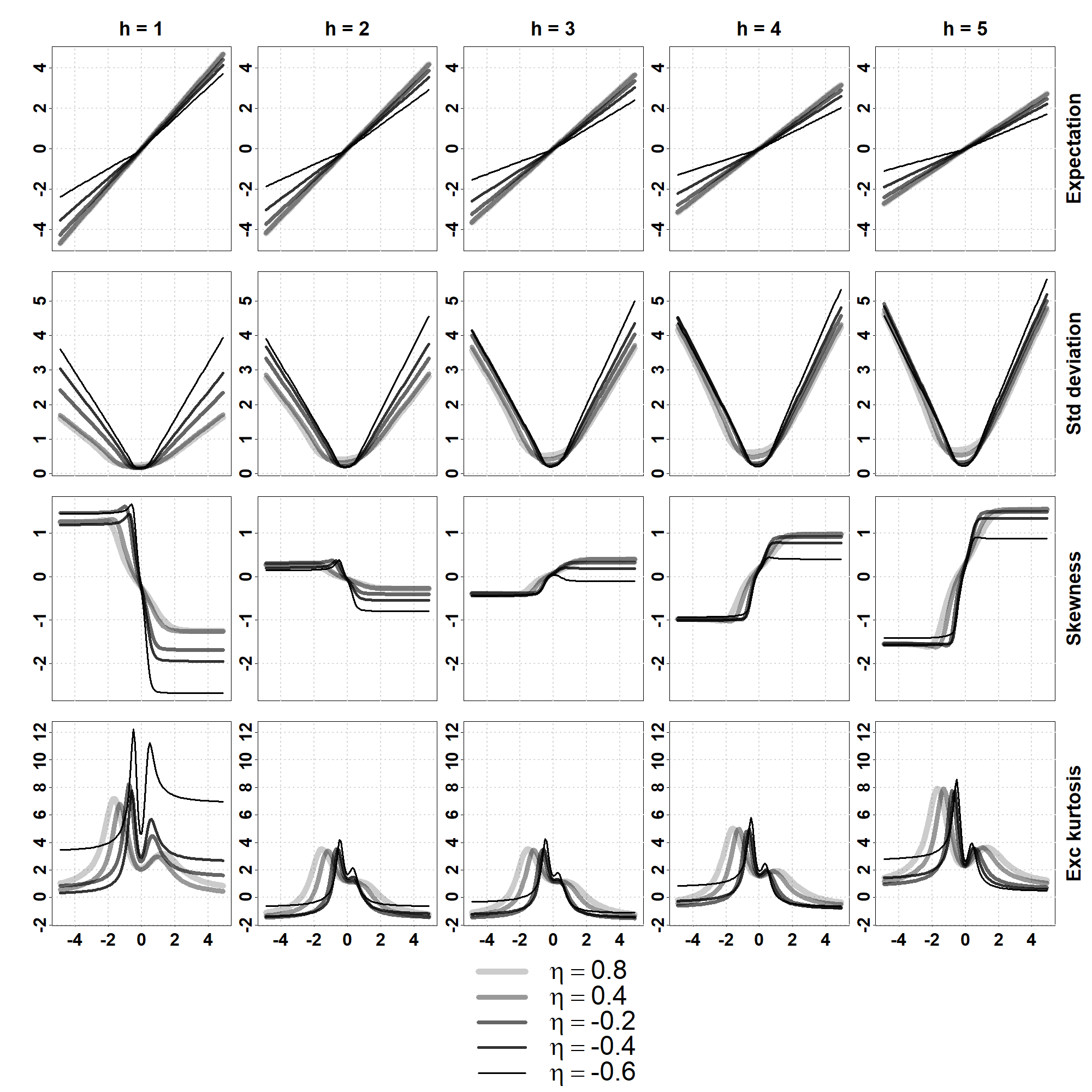}
		    \caption{\textbf{Conditional moments of a stable MAR(2,0,0,1)   for different values of $\boldsymbol{\eta}$.}\\ Conditional expectation, standard deviation, skewness and excess kurtosis (in rows) of $X_{t+h}$ given $X_t=x$, for horizons $h=1,2,3,4,5$ (in columns) and conditioning values $X_t=x\in(-5,5)$ (x-axis of each plot), computed using the formulae of Proposition \ref{prop:2sidedma_condmom}, where $(X_t)$ is the strictly stationary solution of  $(1-0.4 F)(1-0.8 F)X_t=(1+\eta B)\varepsilon_t$, $\varepsilon_t\stackrel{i.i.d.}{\sim}\mathcal{S}(1.7,0.5,0.1,0)$, $\eta\in\{0.8,0.4,-0.2,-0.4,-0.6\}$.}
		    \label{fig:condmom_various_eta}
		\end{figure}

		We provide here figures illustrating how the shape of the conditional moments of MARMA processes may be affected as we let parameter values vary.
		We introduce some notations and shorthands: in the rest of the section, we will denote processes solution of
		$$
		(1-\psi F)(1-\phi B)X_t = (1+\theta F)(1+\eta B)\varepsilon_t,
		$$
		as MARMA($1,1,1,1$), and processes solution of
		$$
		(1-\psi_1 F)(1-\psi_2 F)X_t = (1+\eta B)\varepsilon_t,
		$$
		as MARMA($2,0,0,1$). For the errors, we assume as in the previous sections that $\varepsilon_t\stackrel{i.i.d.}{\sim}\mathcal{S}(\alpha,\beta,\sigma,\mu)$.
		Figures \ref{fig:condmom_various_alp}-\ref{fig:condmom_various_eta} illustrate the effects of the different parameters on the shape of the conditional moments for the two types of MARMA processes above.

		\subsection{Real series application: robustness checks of the crash odds estimation}
		\label{sec:robust}
		
The crash odds estimates of Section \ref{sec:applications} are obtained by maximum likelihood assuming $\alpha$-stable distributed errors.
We here propose some additional empirical results assessing the robustness of these estimates to alternative assumptions on the errors and different fitting methodologies. 
For the Nasdaq and S\&P500 series as in Section \ref{sec:applications}, we fit  noncausal AR(1) models by:
\begin{enumerate}
    \item $t$-student maximum likelihood using the function \texttt{marx} implemented in the R package \texttt{MARX} [\cite{hec17b}]. This approach assumes $t$-distributed errors instead of $\alpha$-stable.
    \item skewed-$t$ regression using the function \texttt{selm} implemented in the R package \texttt{sn} [\cite{azz18}]. This approach assumes skewed-$t$-distributed errors instead of $\alpha$-stable.
    \item OLS estimation of the AR coefficient $\rho$, followed by Hill estimation of $\alpha$ on the residuals of the OLS step, using the function \texttt{hillplot} implemented in the R package \texttt{evmix} [\cite{hs18}]. This approach does not make any fully parametric assumption on the errors, but only assumes they are heavy-tailed in the sense of Equation \eqref{eq:reg_var_errors}.  
\end{enumerate}
Note that in the three approaches above, the errors are power-law tailed and the results of Section \ref{sec:crash_odds} apply.
Table \ref{tab:robust_estim} gathers the estimates of $\rho$ and $\alpha$, while Tables \ref{tab:nasdaq_robustsummary}-\ref{tab:sp500_robustsummary} display the survival statistics of bubbles generated by the corresponding heavy-tailed noncausal AR(1) for the Nasdaq and S\&P500.
One can notice that the survival statistics show similar values to those obtained in Table \ref{tab:bubble_charac}. 

\begin{table}[H]	
\centering
\begin{tabular}{cccc}
\midrule\midrule 
 & & \multicolumn{2}{c}{Noncausal AR(1)} \\
                                       &  & $\alpha$  & $\rho$   \\\cmidrule(lr){3-3} \cmidrule(lr){4-4} 
\multirow{3}{*}{Nasdaq}\hspace*{0.3cm} & $t$-student & 1.22 & 0.979  \\
                                       & \hspace{0.3cm} skewed-$t$ \hspace{0.3cm} & \hspace{0.3cm} 1.18 \hspace{0.3cm} & \hspace{0.3cm}  0.972 \hspace{0.3cm}   \\
                                       & OLS+Hill   & 1.80 & 0.988   \\
  & & &    \\
\multirow{3}{*}{S\&P500}              & $t$-student & 2.02 & 0.987   \\
                                      & skewed-$t$  & 2.12 & 0.983   \\
                                      & OLS+Hill    & 2.50   & 0.992    \\
\midrule\midrule
\end{tabular}
\caption{\footnotesize Estimates of the tail index $\alpha$ and the AR coefficient $\rho$ of heavy-tailed noncausal AR(1) specifications fitted on the real Nasdaq and S\&P500 series using three methodologies: 1) $t$-student maximum likelihood (function \texttt{marx} implemented in the R package \texttt{MARX} [\cite{hec17b}]), 2) skewed-$t$ regression (function \texttt{selm} implemented in the R package \texttt{sn} [\cite{azz18}]), 3) OLS estimation of $\rho$ followed by Hill estimator of $\alpha$ on the residuals of the OLS step (function \texttt{hillplot} implemented in the R package \texttt{evmix} [\cite{hs18}]).
}
\label{tab:robust_estim}
\end{table}

\begin{table}[H]
\centering
		\begin{tabular}{cccccc}
\midrule \midrule
\multicolumn{6}{c}{Nasdaq}  \\
& &   &  &   &    \\
 \multicolumn{6}{c}{\begin{tabular}{ccccc}
	&	Expected life (Y) & Half-life (Y) &  95\%--life (Y) &   99\%--life (Y)  \\  \cmidrule(lr){2-2} \cmidrule(lr){3-3} \cmidrule(lr){4-4} \cmidrule(lr){5-5}
$t$-student	&	3.3 & 2.3 & 9.7 & 14.9   \\
skewed-$t$	&	2.5 & 1.7 & 7.4 & 11.4  \\
OLS+Hill	&	4.0 & 2.7 & 11.7 & 18.0  \\
\end{tabular}}  \\ 
& &   &  &   &    \\
& & \multicolumn{4}{c}{Probability of crash within $h$ months (\%)}\\
& Hazard rate & \hspace*{0.3cm} $h=1$ \hspace*{0.3cm} & \hspace*{0.75cm} 3 \hspace*{0.75cm} & \hspace*{0.3cm} 6 \hspace*{0.3cm} & 12\\ \cmidrule(lr){2-2} \cmidrule(lr){3-6}
$t$-student & 0.025 & 2.5 & 7.4 & 14.3 & 26.5 \\
skewed-$t$ & 0.033 & 3.3 & 9.6 & 18.3 & 33.3 \\
OLS+Hill & 0.021 & 2.1 & 6.2 & 12.0 & 22.5\\
\midrule\midrule
\end{tabular}
\caption{Nasdaq: Summary of bubble survival statistics using the formulae from Table \ref{tab:survival_stat_ar1}, based on estimates of $\alpha$ and $\rho$ of noncausal AR(1) specifications fitted using three methodologies: 1) $t$-student maximum likelihood (function \texttt{marx} implemented in the R package \texttt{MARX} [\cite{hec17b}]), 2) skewed-$t$ regression (function \texttt{selm} implemented in the R package \texttt{sn} [\cite{azz18}]), 3) OLS estimation of $\rho$ followed by Hill estimator of $\alpha$ on the residuals of the OLS step (function \texttt{hillplot} implemented in the R package \texttt{evmix} [\cite{hs18}]). See Table \ref{tab:robust_estim} for the parameter estimates. Expected life, half-life, 95\%-life, 99\%-life are reported in years (Y) and crash probabilities are reported in percents (\%).}
\label{tab:nasdaq_robustsummary}
\end{table}

\begin{table}[H]
\centering
		\begin{tabular}{cccccc}
\midrule \midrule
\multicolumn{6}{c}{S\&P500}  \\
& &   &  &   &    \\
 \multicolumn{6}{c}{\begin{tabular}{ccccc}
	&	Expected life (Y) & Half-life (Y) &  95\%--life (Y) &   99\%--life (Y)  \\  \cmidrule(lr){2-2} \cmidrule(lr){3-3} \cmidrule(lr){4-4} \cmidrule(lr){5-5}
$t$-student	&	3.1 & 2.1 & 9.2 & 14.2   \\
skewed-$t$	&	2.3 & 1.6 & 6.7 & 10.3  \\
OLS+Hill	&	5.8 & 4.0 & 17.2 & 26.4  \\
\end{tabular}}  \\ 
& &   &  &   &    \\
& & \multicolumn{4}{c}{Probability of crash within $h$ months (\%)}\\
& Hazard rate & \hspace*{0.3cm} $h=1$ \hspace*{0.3cm} & \hspace*{0.75cm} 3 \hspace*{0.75cm} & \hspace*{0.3cm} 6 \hspace*{0.3cm} & 12\\ \cmidrule(lr){2-2} \cmidrule(lr){3-6}
$t$-student & 0.027 & 2.7 & 7.8 & 15.0 & 27.8 \\
skewed-$t$ & 0.037 & 3.7 & 10.6 & 20.0 & 36.0 \\
OLS+Hill & 0.014 & 1.4 & 4.3 & 8.4 & 16.0\\
\midrule\midrule
\end{tabular}
\caption{S\&P500: Summary of bubble survival statistics using the formulae from Table \ref{tab:survival_stat_ar1}, based on estimates of $\alpha$ and $\rho$ of noncausal AR(1) specifications fitted using three methodologies: 1) $t$-student maximum likelihood (function \texttt{marx} implemented in the R package \texttt{MARX} [\cite{hec17b}]), 2) skewed-$t$ regression (function \texttt{selm} implemented in the R package \texttt{sn} [\cite{azz18}]), 3) OLS estimation of $\rho$ followed by Hill estimator of $\alpha$ on the residuals of the OLS step (function \texttt{hillplot} implemented in the R package \texttt{evmix} [\cite{hs18}]). See Table \ref{tab:robust_estim} for the parameter estimates. Expected life, half-life, 95\%-life, 99\%-life are reported in years (Y) and crash probabilities are reported in percents (\%).}
\label{tab:sp500_robustsummary}
\end{table}
		
		\clearpage
		
		\label{page:theo22_pini}
		
		\section{Preliminary elements for the proof of the main results}
		
		\label{sec:prelim_elem}
		
		
		\subsection{Notations for the proofs of Theorem \ref{theo:cond_moments_stable_general_alpnot1} and Proposition \ref{prop:equiv}}
		
		The proof of Theorem \ref{theo:cond_moments_stable_general_alpnot1} is quite involved and relies on techniques used in [\cite{ct94,ct98}].
		It consists in differentiating the conditional characteristic function of $X_2|X_1$ up to the fourth derivation order and evaluating the derivatives at 0 to obtain the conditional moments.
		Formal computation of the derivatives yields divergent terms for the third and fourth order derivatives, as well as for the second order derivative when $1/2<\alpha<1$ and special manipulations are needed (in particular the <<appropriate integration by parts>> in \cite{ct94} (p.106) as well as an additional manipulation to obtain the fourth derivative). 
		We first introduce some notations to make the presentation of the proof as compact as possible, then provide the derivatives in Lemma \ref{lem:derivatives} and finally show Theorem \ref{theo:cond_moments_stable_general_alpnot1} by obtaining the functional forms of the conditional moments.

		
		Let $\bX=(X_1,X_2)$ be an $\alpha$-stable vector, with $0<\alpha<2$, $\alpha\ne1$, and spectral representation $(\Gamma,\boldsymbol{0})$. Its characteristic function will be denoted $\varphi_{\bX} (t,r)$ for any $(t,r)\in\mathbb{R}^2$, and reads
		\begin{equation}\label{def:char_proof}
		\varphi_{\bX} (t,r) = \exp\left\{-\int_{S_2} g_1(ts_1+rs_2)\Gamma(ds)\right\},  
		\end{equation}
		where $g_1(z)=|z|^\alpha-iaz^{<\alpha>}$ for $z\in\mathbb{R}$, and $a=\tg(\pi\alpha/2)$. As we assume $\sigma_1>0$ so that $X_1$ is not degenerate, the conditional characteristic function of $X_2$ given $X_1=x$, denoted $\phi_{X_2|x}(r)$ for $r\in\mathbb{R}$, equals
		\begin{equation}\label{def:cond_char_proof}
		\phi_{X_2|x} (r) := 1 + \dfrac{1}{2\pi f_{X_1}(x)} \int_{\mathbb{R}} e^{-itx}\Big(\varphi_{\bX} (t,r) - \varphi_{\bX} (t,0)\Big)dt.
		\end{equation}
		where $f_{X_1}$ denotes the density of $X_1\sim\mathcal{S}(\alpha,\beta_1,\sigma_1,0)$. The following notation of the $\mathcal{H}$ family function will be more handy than that in \eqref{def:hcal_res}: for any $y>-1$ and $\btheta=(\theta_1,\theta_2)\in\mathbb{R}^2$, define the function $\mathcal{H}(y,\btheta; \,\cdot\,)$ 
		for $x\in\mathbb{R}$ as
		\begin{align}\label{def:hcal}
		\mathcal{H}(y,\btheta;x) & =  \int_0^{+\infty} e^{-\sigma_1^\alpha u^\alpha} u^{y} \Big(\theta_1\cos(ux-a\beta_1\sigma_1^\alpha u^\alpha) + \theta_2\sin(ux-a\beta_1\sigma_1^\alpha u^\alpha)\Big)du,
		\end{align}
		For $z\in\mathbb{R}$, denote also,
		
		\vspace*{-1.4cm}
		
		\begin{align}
		g_2(z) & = z^{<\alpha-1>} -ia |z|^{\alpha-1},\label{def:g2}\\
		g_3(z) & = |z|^{\alpha-2} -ia z^{<\alpha-2>}.\label{def:g3}
		\end{align}
		
		\vspace*{-0.2cm}
		
		Often, we shall invoke functions of the form
		\begin{align}\label{def:notshorthand}
		r\longmapsto\int_{\mathbb{R}}e^{-itx}\varphi_{\bX}(t,r)f_1^{p_1}(t,r)\ldots f_m^{p_m}(t,r)dt,
		\end{align}
		where $m\le3$ and the $f_i$'s will be functions of the type $f_i(t,r)=\int_{S_2}g_{j_i}(ts_1+rs_2)s_1^{k_i}s_2^{\ell_i}\Gamma(d\bs)$, for $j_i=2,3$, $k_i,\ell_i\in\mathbb{Z}$ for which $f_i$ is well defined and positive integer exponents $p_i$'s. As a shorthand when no ambiguity is possible, we shall denote functions like \eqref{def:notshorthand} by
		\begin{align*}
		\Lambda\bigg(\int_{S_2}g_{j_1}s_1^{k_1}s_2^{\ell_1}\bigg)^{p_1}\bigg(\int_{S_2}g_{j_2}s_1^{k_2}s_2^{\ell_2}\bigg)^{p_2}\ldots
		\end{align*}
		
		\vspace*{-0.2cm}
		
		\noindent up to the $m$\textsuperscript{th} term. 
		
		\subsection{Lemma \ref{lem:derivatives} for the proof of Theorem \ref{theo:cond_moments_stable_general_alpnot1}}
		
		\begin{lem}
			\label{lem:derivatives}
			Let $(X_1,X_2)$ be an $\alpha$-stable vector, $0<\alpha<2$,$\alpha\ne1$, with conditional characteristic function $\phi_{X_2|x}$ as given in \eqref{def:cond_char_proof}. Let $r\in\mathbb{R}$. If $1<\alpha<2$, or if $0<\alpha<1$ and \eqref{eq:nu_cond} holds with $\nu>1-\alpha$, the first derivative of $\phi_{X_2|x}$ is given by
			\begin{align}
			\phi^{(1)}_{X_2|x}(r) & = \dfrac{-\alpha}{2\pi f_{X_1}(x)}\Lambda \bigg(\int_{S_2}g_2s_2\bigg).\label{eq:d1}
			\end{align}
			If $1/2<\alpha<2$ and \eqref{eq:nu_cond} holds with $\nu>2-\alpha$, the second derivative is given by
			\begin{align}
			\phi^{(2)}_{X_2|x}(r) & = \dfrac{-\alpha}{2\pi f_{X_1}(x)} \Bigg[ix\Lambda \bigg(\int_{S_2}g_2s_2^2s_1^{-1}\bigg)+\alpha\Bigg\{\Lambda \bigg(\int_{S_2}g_2s_2^2s_1^{-1}\bigg)\bigg(\int_{S_2}g_2s_1\bigg)-\Lambda \bigg(\int_{S_2}g_2s_2^2\bigg)^2\Bigg\}\Bigg],\label{eq:d2}
			\end{align}
			If $1<\alpha<2$ and \eqref{eq:nu_cond} holds with $\nu>3-\alpha$, the third derivative is given by
			\begin{align}
			\phi^{(3)}_{X_2|x}(r) & = \dfrac{-\alpha}{2\pi f_{X_1}(x)} \bigg(ix\Big((\alpha-1)I_1-\alpha I_2\Big)+\alpha^2(I_3-I_4)+\alpha(\alpha-1)(I_5+I_6-2I_7)\bigg),\label{eq:d3}
			\end{align}
			with 
			\begin{align*}
			I_1 & = \Lambda \bigg(\int_{S_2} g_3 s_2^3s_1^{-1}\bigg), & I_5 & = \Lambda \bigg(\int_{S_2} g_2 s_2^2s_1^{-1}\bigg)\bigg(\int_{S_2} g_3 s_2s_1\bigg),\\
			I_2 & = \Lambda \bigg(\int_{S_2} g_2 s_2\bigg)\bigg(\int_{S_2} g_2 s_2^2s_1^{-1}\bigg), & I_6 & = \Lambda \bigg(\int_{S_2} g_2 s_1\bigg)\bigg(\int_{S_2} g_3 s_2^3s_1^{-1}\bigg),\\
			I_3 & = \Lambda \bigg(\int_{S_2} g_2 s_2\bigg)^3, & I_7 & = \Lambda \bigg(\int_{S_2} g_2 s_2\bigg)\bigg(\int_{S_2} g_3 s_2^2\bigg),\\
			I_4 & = \Lambda \bigg(\int_{S_2} g_2 s_1\bigg)\bigg(\int_{S_2} g_2 s_2\bigg)\bigg(\int_{S_2} g_2 s_2^2s_1^{-1}\bigg).
			\end{align*}
			If $3/2<\alpha<2$ and \eqref{eq:nu_cond} holds with $\nu>4-\alpha$, the fourth derivative is given by
			\begin{align}
			\phi^{(4)}_{X_2|x}(r)  =   \dfrac{-\alpha}{2\pi f_{X_1}(x)}  \Bigg[ \hspace{0.2cm } & i\alpha x \bigg(\alpha\Big(3J_1-2J_2\Big) + (\alpha-1)\Big(2J_3-3J_4+J_5\Big)\bigg) + \alpha x^2 J_6 - (\alpha-1)x^2J_7 \nonumber\\
			& + \alpha^2(\alpha-1) \bigg(J_8+J_9+J_{10}-3\Big(2J_{11}+J_{12}-J_{13}\Big)\bigg) \label{eq:d4} \\
			& + \alpha(\alpha-1)^2 \bigg(4J_{14} - 3J_{15}-J_{16}\bigg)  + \alpha^3\bigg(3J_{17}-J_{18}-J_{19}\bigg)\Bigg],\nonumber
			\end{align}
			with
			\begin{align*}
			J_1 & = \Lambda \Big(\int_{S_2}g_2s_2^2s_1^{-1}\Big)\Big(\int_{S_2}g_2s_2\Big)^2, & J_{11} & = \Lambda \Big(\int_{S_2}g_2s_2^2s_1^{-1}\Big)\Big(\int_{S_2}g_3s_2s_1\Big)\Big(\int_{S_2}g_2s_2\Big),\\
			J_2 & = \Lambda \Big(\int_{S_2}g_2s_2^3s_1^{-2}\Big)\Big(\int_{S_2}g_2s_1\Big)\Big(\int_{S_2}g_2s_2\Big), & J_{12} & = \Lambda \Big(\int_{S_2}g_3s_2^3s_1^{-1}\Big)\Big(\int_{S_2}g_2s_1\Big)\Big(\int_{S_2}g_2s_2\Big),\\
			J_3 & = \Lambda \Big(\int_{S_2}g_3s_2^4s_1^{-2}\Big)\Big(\int_{S_2}g_2s_1\Big), & J_{13} & = \Lambda \Big(\int_{S_2}g_3s_2^2\Big)\Big(\int_{S_2}g_2s_2\Big)^2,\\
			J_4 & = \Lambda \Big(\int_{S_2}g_3s_2^3s_1^{-1}\Big)\Big(\int_{S_2}g_2s_2\Big), & J_{14} & = \Lambda \Big(\int_{S_2}g_3s_2^3s_1^{-1}\Big)\Big(\int_{S_2}g_3s_2s_1\Big),\\
			J_5 & = \Lambda \Big(\int_{S_2}g_2s_2^3s_1^{-2}\Big)\Big(\int_{S_2}g_3s_2s_1\Big), & J_{15} & =
			\Lambda \Big(\int_{S_2}g_3s_2^2\Big)^2,\\
			J_6 & = \Lambda \Big(\int_{S_2}g_2s_2^3s_1^{-2}\Big)\Big(\int_{S_2}g_2s_2\Big), & J_{16} & =
			\Lambda \Big(\int_{S_2}g_3s_2^4s_1^{-2}\Big)\Big(\int_{S_2}g_3s_1^{2}\Big),\\
			J_7 & = \Lambda \Big(\int_{S_2}g_3s_2^4s_1^{-2}\Big), & J_{17} & =
			\Lambda \Big(\int_{S_2}g_2s_2^2s_1^{-1}\Big)\Big(\int_{S_2}g_2s_1\Big)\Big(\int_{S_2}g_2s_2\Big)^2,\\
			J_8 & = \Lambda \Big(\int_{S_2}g_2s_2^3s_1^{-2}\Big)\Big(\int_{S_2}g_3s_1^2\Big)\Big(\int_{S_2}g_2s_2\Big), & J_{18} & = \Lambda \Big(\int_{S_2}g_2s_2\Big)^4, \\
			J_9 & = \Lambda \Big(\int_{S_2}g_2s_2^3s_1^{-2}\Big)\Big(\int_{S_2}g_3s_2s_1\Big)\Big(\int_{S_2}g_2s_1\Big), & J_{19} & = \Lambda \Big(\int_{S_2}g_2s_2^3s_1^{-2}\Big)\Big(\int_{S_2}g_2s_1\Big)^2\Big(\int_{S_2}g_2s_2\Big),\\
			J_{10} & = \Lambda \Big(\int_{S_2}g_3s_2^4s_1^{-2}\Big)\Big(\int_{S_2}g_2s_1\Big)^2. & &
			\end{align*}
		\end{lem}

		\section{Proof of Lemma \ref{lem:derivatives}}
		\label{sec:lemc2}
		
		
		For each of the derivatives, the proof involves two main steps: 1) justifying inversion of integral and derivation signs 2) computation of the derivative. 
		
		\subsection{Justifying inversion of integral and derivation signs}
		
		\subsubsection{Justifying inversion: First derivative}
		
		\noindent \textbf{Case $\boldsymbol{\alpha\in(0,1)}$}
		
		Assume $\alpha\in(0,1)$. We begin with the first derivative of the imaginary part of $\phi_{X_2|x}$. 
		\begin{align}
		&\dfrac{d}{dr} \Big(\text{Im}\phi_{X_2|x}(r)\Big) \nonumber \\
		& \hspace{1cm} =  \dfrac{-1}{2\pi f_{X_1}(x)}\lim_{h\rightarrow0} \dfrac{1}{h} \int_{\mathbb{R}}\Bigg[ e^{-\int_{S_2}|ts_1+(r+h)s_2|^\alpha\Gamma(ds)}\sin\left(tx-a\int_{S_2}(ts_1+(r+h)s_2)^{<\alpha>}\Gamma(ds)\right) \nonumber\\
		& \hspace{6cm} - e^{-\int_{S_2}|ts_1+rs_2|^\alpha\Gamma(ds)}\sin\left(tx-a\int_{S_2}(ts_1+rs_2)^{<\alpha>}\Gamma(ds)\right) \Bigg] dt\nonumber\\
		& \hspace{1cm} = \dfrac{-1}{2\pi f_{X_1}(x)} \lim_{h\rightarrow0} \dfrac{1}{h} \int_{\mathbb{R}}\Bigg[ \sin\left(tx-a\int_{S_2}(ts_1+(r+h)s_2)^{<\alpha>}\Gamma(ds)\right)\nonumber\\
		& \hspace{7cm}-\sin\left(tx-a\int_{S_2}(ts_1+rs_2)^{<\alpha>}\Gamma(ds)\right) \Bigg]\nonumber\\
		& \hspace{10cm}\times \exp\Big\{-\int_{S_2}|ts_1+rs_2|^\alpha\Gamma(ds)\Big\} dt \nonumber\\
		&  \hspace{1.5cm} - \dfrac{1}{2\pi f_{X_1}(x)} \lim_{h\rightarrow0} \dfrac{1}{h} \int_{\mathbb{R}}\Bigg[\exp\Big\{-\int_{S_2}|ts_1+(r+h)s_2|^\alpha\Gamma(ds)\Big\}\nonumber\\
		& \hspace{6.5cm}  -\exp\Big\{-\int_{S_2}|ts_1+rs_2|^\alpha\Gamma(ds)\Big\}\Bigg] \nonumber\\
		& \hspace{8cm}\times \sin\left(tx-a\int_{S_2}(ts_1+(r+h)s_2)^{<\alpha>}\Gamma(ds)\right) dt \nonumber\\
		& \hspace{1cm} := I_1 + I_2. \label{eq:I1_I2}
		\end{align}
		The integrand of $I_1$ converges to 
		$$
		-\alpha a\cos\Big(tx-a\int_{S_2}(ts_1+rs_2)^{<\alpha>}\Gamma(ds)\Big) \times \int_{S_2}|ts_1+rs_2|^{\alpha-1}s_2\Gamma(ds) \times\exp\Big\{-\int_{S_2}|ts_1+rs_2|^\alpha\Gamma(ds)\Big\}
		$$
		Using the mean value theorem, the triangle inequality and the inequality $-|x+y|^\alpha\le-|x|^\alpha+|y|^\alpha$ when $0<\alpha<1$, the integrand of $I_1$ can be bounded for any $h$, $|h|<|r|$, by
		\begin{align}
		& \Big|cos (y)\Big|\left(\Big|\dfrac{a}{h}\Big| \int_{S_2}\Big|(ts_1+(r+h)s_2)^{<\alpha>}-(ts_1+rs_2)^{<\alpha>}\Big|\Gamma(ds) \right)\exp\Big\{\int_{S_2} -|ts_1|^\alpha + |rs_2|^\alpha\Gamma(ds)\Big\}\nonumber \\
		& \le 2|a|e^{|r|^\alpha\sigma_2^\alpha}e^{-\sigma_1^\alpha |t|^\alpha} \int_{S_2}|ts_1+rs_2|^{\alpha-1} \Gamma(ds),\label{dem:bd1}
		\end{align}
		where $\sigma_2=\Big(\int_{S_2}|s_2|^\alpha\Gamma(d\bs)\Big)^{1/\alpha}$, $y\in\mathbb{R}$, and  we used the bound
		\begin{equation}\label{dem:bound215}
		\left|\dfrac{(ts_1+(r+h)s_2)^{<\alpha>}-(ts_1+rs_2)^{<\alpha>}}{h}\right| \le 2|ts_1+rs_2|^{\alpha-1}|s_2|,   
		\end{equation}
		for $ts_1+rs_2\ne0$, which is a consequence of 
		$||1+z|^{<\alpha>}-1| \le 2|z|$, for $z\in\mathbb{R}$ (see Lemma \ref{le:34} $(\iota\iota)$ below). Bound \eqref{dem:bd1} does not depend on $h$ and is integrable with respect to $t$. Indeed, invoking Lemma \ref{le:310} with $\eta=\alpha-1$, $b=p=0$,
		and \eqref{eq:nu_cond} with $\nu>2-\alpha>1-\alpha$
		\begin{align}
		& \Bigg|\int_{\mathbb{R}}e^{-\sigma_1^\alpha |t|^\alpha} \int_{S_2}\Big|t+r\frac{s_2}{s_1}\Big|^{\alpha-1}|s_1|^{\alpha-1}\Gamma(ds)dt-\int_{\mathbb{R}}\int_{S_2}e^{-\sigma_1^\alpha |t|^\alpha}|t|^{\alpha-1}|s_1|^{\alpha-1}\Gamma(d\bs)dt\Bigg|\nonumber\\
		& \hspace{2cm}\le \int_{S_2} |s_1|^{\alpha-1} \int_{\mathbb{R}}e^{-\sigma_1^\alpha |t|^\alpha}\Bigg|\Big|t+r\frac{s_2}{s_1}\Big|^{\alpha-1}-|t|^{\alpha-1}\Bigg|dt \Gamma(ds)\nonumber\\
		& \hspace{2cm}\le \text{const} \int_{S_2} |s_1|^{\alpha-1+\nu}|s_1|^{-\nu}\Gamma(ds)\nonumber\\
		& \hspace{2cm}  \le \text{const} \int_{S_2} |s_1|^{-\nu}\Gamma(ds)\nonumber\\
		& \hspace{2cm} < +\infty, \label{dem:bound_lem22}
		\end{align}
		and the integrability with respect to $t$ follows from the fact that $\int_{\mathbb{R}}e^{-\sigma_1^\alpha |t|^\alpha}|t|^{\alpha-1}dt<+\infty$. Hence the Lebesgue dominated convergence theorem applies to $I_1$ and we can invert integration and derivation.
		Focusing on $I_2$, its integrand tends to
		$$
		-\alpha \int_{S_2}(ts_1+rs_2)^{<\alpha-1>}s_2 \Gamma(ds) \exp\left\{-\int_{S_2}|ts_1+rs_2|^\alpha\Gamma(ds)\right\}\sin\left(tx-a\int_{S_2}|ts_1+rs_2|^{<\alpha>}\Gamma(ds)\right).
		$$
		Using the inequality 
		$$
		\left|\dfrac{(ts_1+(r+h)s_2)^{\alpha}-(ts_1+rs_2)^{\alpha}}{h}\right| \le |ts_1+rs_2|^{\alpha-1}|s_2|,
		$$
		for $ts_1+rs_2\ne0$, which is a consequence of $||1+z|^{\alpha}-1| \le |z|$, for $z\in\mathbb{R}$ (Lemma \ref{le:34} $(\iota)$ below) 
		and the inequality $|e^{-x}-e^{-y}|\le e^{-y}e^{|x-y|}|x-y|$, for $x,y\in\mathbb{R}$, we can bound the integrand of $I_2$ for any $|h|<|r|$ by
		\begin{align*}
		& \exp\left\{-\int_{S_2}|ts_1+rs_2|^\alpha\Gamma(ds)\right\} \exp\left\{\left|\int_{S_2}|ts_1+(r+h)s_2|^\alpha-|ts_1+rs_2|^\alpha\Gamma(ds)\right|\right\}\\ & \hspace{8cm} \times\left|\dfrac{1}{h}\int_{S_2}|ts_1+(r+h)s_2|^\alpha-|ts_1+rs_2|^\alpha\Gamma(ds)\right|\\
		& \le e^{2|r|^\alpha\sigma_2^\alpha}e^{-\sigma_1^\alpha |t|^\alpha} \int_{S_2}\Big|t+r\frac{s_2}{s_1}\Big|^{\alpha-1}|s_1|^{\alpha-1}\Gamma(ds).
		\end{align*}
		The integrability with respect to $t$ is deduced as for \eqref{dem:bound_lem22} using Lemma \ref{le:310} with $\eta=\alpha-1$, $b=p=0$. Thus, the Lebesgue-dominated convergence theorem applies to $I_2$ and we can invert integration and derivation. The real part of $\phi_{X_2|x}(r)$ can be treated in a similar way, allowing us to derivate under the integral. \\
		
		\noindent \textbf{Case $\boldsymbol{\alpha\in(1,2)}$}\\
		Assume $\alpha\in(1,2)$. Just as for the case $\alpha\in(0,1)$, the imaginary part of $\phi_{X_2|x}$ is given by \eqref{eq:I1_I2} 
		\begin{align*}
		&\dfrac{d}{dr} \Big(\text{Im}\phi_{X_2|x}(r)\Big) = I_1 + I_2.
		\end{align*}
		The integrands of $I_1$ and $I_2$ still converges to the same limits, however a different argument is needed to bound them. For $|h|<|r|$, the mean value theorem, the triangle inequality and the inequality of Lemma \ref{le:33}, yield the following bound for the integrand of $I_1$
		\begin{align}
		& \left(\Big|\dfrac{a}{h}\Big| \int_{S_2}\Big|(ts_1+(r+h)s_2)^{<\alpha>}-(ts_1+rs_2)^{<\alpha>}\Big|\Gamma(ds) \right)e^{|r|^\alpha\sigma_2^\alpha}e^{-2^{1-\alpha}\sigma_1^\alpha |t|^\alpha},\label{dem:bd1_alpsup1}
		\end{align}
		where $y\in\mathbb{R}$. By the triangle inequality and the mean value theorem, we have for some
		$u\in\bigg(\min\Big(ts_1+(r+h)s_2,ts_1+rs_2\Big),\max\Big(ts_1+(r+h)s_2,ts_1+rs_2\Big)\bigg)$
		\begin{align}
		\bigg|\int_{S_2}(ts_1+(r+h)s_2)^{<\alpha>}-(ts_1+rs_2)^{<\alpha>}\Gamma(d\bs)\bigg| & = \bigg|\int_{S_2} \alpha h s_2 |u|^{\alpha-1} \Gamma(d\bs)\bigg|\nonumber\\
		& \le \alpha |h| \bigg|\int_{S_2} |t|^{\alpha-1} + 2|r|^{\alpha-1} \Gamma(d\bs)\nonumber\\
		& \le \alpha |h| \Gamma(S_2) (|t|^{\alpha-1} + 2|r|^{\alpha-1})\label{eq:mvtcrochets}
		\end{align}
		Thus, \eqref{dem:bd1_alpsup1} can be bounded by
		\begin{align*}
		& \alpha|a| \Gamma(S_2)e^{|r|^\alpha\sigma_2^\alpha}e^{-2^{1-\alpha}\sigma_1^\alpha |t|^\alpha} (|t|^{\alpha-1} + 2|r|^{\alpha-1}),
		\end{align*}
		which is certainly integrable with respect to $t$ on $\mathbb{R}$ for $\alpha>1$. Let us now turn to $I_2$. We have again by the mean value theorem,
		\begin{align*}
		\bigg|\dfrac{|ts_1+(r+h)s_2|^\alpha-|ts_1+rs_2|^\alpha}{h}\bigg| \le \alpha(|t|^{\alpha-1}+2|r|^{\alpha-1}),
		\end{align*}
		if $|h|<|r|$, and thus
		\begin{align}
		& \bigg|\dfrac{e^{-\int_{S_2}|ts_1+(r+h)s_2|^\alpha\Gamma(d\bs)}-e^{-\int_{S_2}|ts_1+rs_2|^\alpha\Gamma(d\bs)}}{h}\bigg| \nonumber\\
		& \hspace{5cm} \le \max\bigg(e^{-\int_{S_2}|ts_1+(r+h)s_2|^\alpha\Gamma(d\bs)},e^{-\int_{S_2}|ts_1+rs_2|^\alpha\Gamma(d\bs)}\bigg)\nonumber\\
		& \hspace{6cm} \times \int_{S_2}\bigg|\dfrac{|ts_1+(r+h)s_2|^\alpha-|ts_1+rs_2|^\alpha}{h}\bigg|\Gamma(d\bs)\nonumber\\
		& \hspace{5cm} \le \Gamma(S_2)e^{|2r|^\alpha\sigma_2^\alpha}e^{-2^{1-\alpha}\sigma_1^\alpha|t|^\alpha}\alpha(|t|^{\alpha-1}+2|r|^{\alpha-1}),\label{eq:e-min}
		\end{align}
		by Lemma \ref{le:exp} \eqref{eq:exp95} and Lemma \ref{le:33}.
		The latter bound is again integrable with respect to $t$ on $\mathbb{R}$. Hence the dominated convergence theorem applies to $I_1$, $I_2$ and therefore to $\dfrac{d}{dr} \Big(\text{Im}\phi_{X_2|x}(r)\Big)$ and we can invert the integration and derivation signs. Similar arguments show the dominated convergence theorem applies to the real part of the conditional characteristic function as well. 
		
		\subsubsection{Justifying inversion: Second derivative}
		
		\noindent \textbf{Case $\boldsymbol{\alpha\in(1/2,1)}$}\\
		
		In an expanded fashion, $\phi^{(1)}_{X_2|x}(r)$ can be written,
		\begin{equation}\label{eq:phi'}
		\phi^{(1)}_{X_2|x}(r) = \dfrac{-\alpha}{2\pi f_{X_1}(x)} \Big[J_{1} - a J_2 -i(J_3+aJ_4)\Big],
		\end{equation}
		with, 
		\begin{align*}
		J_1(r) & = \int_{\mathbb{R}}e^{-\int_{S_2}|ts_1+rs_2|^{\alpha}\Gamma(ds)}\cos\Bigg(tx-a\int_{S_2}(ts_1+rs_2)^{<\alpha>}\Gamma(ds)\Bigg)\int_{S_2}(ts_1+rs_2)^{<\alpha-1>}s_2\Gamma(ds)dt,\\
		J_2(r) & = \int_{\mathbb{R}}e^{-\int_{S_2}|ts_1+rs_2|^{\alpha}\Gamma(ds)}\sin\Bigg(tx-a\int_{S_2}(ts_1+rs_2)^{<\alpha>}\Gamma(ds)\Bigg)\int_{S_2}|ts_1+rs_2|^{\alpha-1}s_2\Gamma(ds)dt,\\
		J_3(r) & = \int_{\mathbb{R}}e^{-\int_{S_2}|ts_1+rs_2|^{\alpha}\Gamma(ds)}\sin\Bigg(tx-a\int_{S_2}(ts_1+rs_2)^{<\alpha>}\Gamma(ds)\Bigg)\int_{S_2}(ts_1+rs_2)^{<\alpha-1>}s_2\Gamma(ds)dt,\\
		J_4(r) & = \int_{\mathbb{R}}e^{-\int_{S_2}|ts_1+rs_2|^{\alpha}\Gamma(ds)}\cos\Bigg(tx-a\int_{S_2}(ts_1+rs_2)^{<\alpha>}\Gamma(ds)\Bigg)\int_{S_2}|ts_1+rs_2|^{\alpha-1}s_2\Gamma(ds)dt.\\
		\end{align*}
		To obtain $\phi^{(2)}_{X_2|x}(r)$, we will show that the dominated convergence theorem applies to $J_1'$. Let us consider,
		\begin{align}
		J_1'(r) & = \lim_{h\rightarrow0} \frac{1}{h} \int_{\mathbb{R}} \Bigg[ \exp\Big\{-\int_{S_2}|ts_1+(r+h)s_2|^{\alpha}\Gamma(d\bs)\Big\}\cos\Big(tx-a\int_{S_2}(ts_1+(r+h)s_2)^{<\alpha>}\Gamma(ds)\Big)\nonumber\\
		& \hspace{7cm} \times  \int_{S_2} (ts_1+(r+h)s_2)^{<\alpha-1>}s_2\Gamma(ds)\nonumber\\
		& \hspace{3cm} - \exp\Big\{-\int_{S_2}|ts_1+rs_2|^{\alpha}\Gamma(d\bs)\Big\}\cos\Big(tx-a\int_{S_2}(ts_1+rs_2)^{<\alpha>}\Gamma(ds)\Big)\nonumber\\
		& \hspace{7cm} \times  \int_{S_2} (ts_1+rs_2)^{<\alpha-1>}s_2\Gamma(ds)
		\Bigg]dt\nonumber\\
		& = \lim_{h\rightarrow0} \frac{1}{h} \int_{\mathbb{R}} \Bigg[\exp\Big\{-\int_{S_2}|ts_1+(r+h)s_2|^{\alpha}\Gamma(d\bs)\Big\} - \exp\Big\{-\int_{S_2}|ts_1+rs_2|^{\alpha}\Gamma(d\bs)\Big\}\Bigg]\nonumber\\
		& \hspace{2cm} \times \cos\Big(tx-a\int_{S_2}(ts_1+rs_2)^{<\alpha>}\Gamma(ds)\Big) \int_{S_2} (ts_1+rs_2)^{<\alpha-1>}s_2\Gamma(ds) dt\nonumber\\
		&  + \lim_{h\rightarrow0} \frac{1}{h} \int_{\mathbb{R}} \exp\Big\{-\int_{S_2}|ts_1+(r+h)s_2|^{\alpha}\Gamma(d\bs)\Big\}\\
		&  \hspace{3cm}\times \Bigg[\cos\Big(tx-a\int_{S_2}(ts_1+(r+h)s_2)^{<\alpha>}\Gamma(ds)\Big)\nonumber\\
		& \hspace{6cm} - \cos\Big(tx-a\int_{S_2}(ts_1+rs_2)^{<\alpha>}\Gamma(ds)\Big)\Bigg]\nonumber\\
		& \hspace{7cm} \times  \int_{S_2} (ts_1+rs_2)^{<\alpha-1>}s_2\Gamma(d\bs) dt\nonumber\\
		& + \lim_{h\rightarrow0} \frac{1}{h} \int_{\mathbb{R}} \exp\Big\{-\int_{S_2}|ts_1+(r+h)s_2|^{\alpha}\Gamma(d\bs)\Big\}\cos\Big(tx-a\int_{S_2}(ts_1+(r+h)s_2)^{<\alpha>}\Gamma(ds)\Big)\nonumber\\
		& \hspace{2cm} \times  \Bigg[\int_{S_2} (ts_1+(r+h)s_2)^{<\alpha-1>}s_2\Gamma(ds)-\int_{S_2} (ts_1+rs_2)^{<\alpha-1>}s_2\Gamma(ds)\Bigg]dt\nonumber\\
		& := K_1 + K_2 + K_3. \label{eq:j1}
		\end{align}
		It can be shown that the dominated convergence theorem applies to $K_1$ following the proof in \cite{ct94} (p.105) for $I_1$. Consider $K_2$. The integrand converges to 
		\begin{align*}
		\alpha a \Bigg(\int_{S_2} & |ts_1+rs_2|^{\alpha-1}s_2\Gamma(ds)\Bigg)\Bigg(\int_{S_2} (ts_1+rs_2)^{<\alpha-1>}s_2\Gamma(ds)\Bigg)\\
		& \hspace{3cm} \times \sin \Big(tx-a\int_{S_2} (ts_1+rs_2)^{<\alpha>}\Gamma(ds)\Big)\exp\Big\{-\int_{S_2}|ts_1+rs_2|^{\alpha}\Gamma(ds)\Big\}.
		\end{align*}
		Using the mean value theorem, \eqref{dem:bound215} and the triangle inequality, we can bound the integrand for any $|h|<|r|$ by
		\begin{align}
		& \Bigg|\frac{1}{h}\int_{S_2} (ts_1+(r+h)s_2)^{<\alpha>}-(ts_1+rs_2)^{<\alpha>}\Gamma(ds)\Bigg|\nonumber\\
		& \hspace{2cm}\times |\sin(y)|e^{2|r|^\alpha\sigma_2^\alpha}e^{-|t|^\alpha\sigma_1^\alpha} \int_{S_2} \Big|t+r\frac{s_2}{s_1}\Big|^{\alpha-1}|s_2| |s_1|^{\alpha-1}\Gamma(ds)\nonumber\\
		& \le 2 e^{2|r|^\alpha\sigma_2^\alpha}\Bigg(\int_{S_2} \Big|t+r\frac{s_2}{s_1}\Big|^{\alpha-1} |s_1|^{\alpha-1}\Gamma(ds)\Bigg)^2 e^{-|t|^\alpha\sigma_1^\alpha} \label{dem:bdk2}
		\end{align}
		where $y\in\mathbb{R}$. The bound \eqref{dem:bdk2} does not depend on $h$ and is integrable with respect to $t$: invoking (2.9) Lemma 2.2 in \cite{ct94}, 
		\begin{align}
		& \Bigg|\int_{\mathbb{R}}\int_{S_2}\int_{S_2}e^{-\sigma_1^\alpha |t|^\alpha} \Big|t+r\frac{s_2}{s_1}\Big|^{\alpha-1}\Big|t+r\frac{s'_2}{s'_1}\Big|^{\alpha-1}|s'_1|^{\alpha-1}|s_1|^{\alpha-1}\Gamma(ds)\Gamma(ds')dt\\
		& \hspace{7cm}-\int_{\mathbb{R}}\int_{S_2}\int_{S_2}e^{-\sigma_1^\alpha |t|^\alpha}|t|^{2\alpha-2}dt\Gamma(d\bs)\Gamma(d\bs')\Bigg|\nonumber\\
		& \hspace{0.5cm} = \Bigg|\int_{S_2}\int_{S_2}|s'_1|^{\alpha-1}|s_1|^{\alpha-1}\int_{\mathbb{R}} e^{-\sigma_1^\alpha |t|^\alpha}\Bigg[\Big|t+r\frac{s_2}{s_1}\Big|^{\alpha-1}\Big|t+r\frac{s'_2}{s'_1}\Big|^{\alpha-1} - \Big|t+r\frac{s_2}{s_1}\Big|^{\alpha-1}\Big|t\Big|^{\alpha-1}
		\nonumber\\ 
		& \hspace{7.5cm} + \Big|t+r\frac{s_2}{s_1}\Big|^{\alpha-1}\Big|t\Big|^{\alpha-1} -|t|^{2\alpha-2}\Bigg]dt\Gamma(ds)\Gamma(ds')\Bigg|\nonumber\\
		& \hspace{0.5cm} \le \int_{S_2}\int_{S_2}|s'_1|^{\alpha-1}|s_1|^{\alpha-1}\int_{\mathbb{R}} e^{-\sigma_1^\alpha |t|^\alpha} \Bigg[ \bigg|\Big|t+r\frac{s'_2}{s'_1}\Big|^{\alpha-1}-|t|^{\alpha-1}\bigg| \Big|t+r\frac{s_2}{s_1}\Big|^{\alpha-1} \nonumber\\
		& \hspace{7.5cm} +\bigg|\Big|t+r\frac{s_2}{s_1}\Big|^{\alpha-1}-|t|^{\alpha-1}\bigg||t|^{\alpha-1}\Bigg]dt\Gamma(ds)\Gamma(ds')\nonumber\\
		& \hspace{0.5cm}\le \text{const} \Bigg(\int_{S_2}|s_1|^{\alpha-1}\Gamma(ds)\Bigg)^2\nonumber\\
		& \hspace{0.5cm} < +\infty,\label{dem:bdk2_29}
		\end{align}
		where const is a constant depending only on $\alpha$ and $\sigma_1^{\alpha}$. The integrability of \eqref{dem:bdk2} follows from \eqref{dem:bdk2_29}, the fact that $\int_{\mathbb{R}} e^{-\sigma_1^\alpha |t|^\alpha}|t|^{2\alpha-2}dt < +\infty$ and \eqref{eq:nu_cond} with $\nu>2-\alpha>1-\alpha$. Hence the dominated convergence theorem applies to $K_2$. Let us now turn to $K_3$: <<this is [a] case when appropriate "integration by parts" is needed>> (\cite{ct94}). With the change of variable $t'=t+\dfrac{hs_2'}{s_1'}$,
		\begin{align*}
		K_3 & = \lim_{h\rightarrow0} \frac{1}{h} \Bigg[ \int_{\mathbb{R}} \exp\Big\{-\int_{S_2}|ts_1+(r+h)s_2|^{\alpha}\Gamma(d\bs)\Big\}\cos\Big(tx-a\int_{S_2}(ts_1+(r+h)s_2)^{<\alpha>}\Gamma(ds)\Big)\\
		& \hspace{6cm} \times \int_{S_2} (t+\frac{hs'_2}{s'_1}+\frac{rs'_2}{s'_1})^{<\alpha-1>}s'_2{s'_1}^{<\alpha-1>}\Gamma(ds')dt \\
		& \hspace{2cm} - \int_{\mathbb{R}} \exp\Big\{-\int_{S_2}|ts_1+(r+h)s_2|^{\alpha}\Gamma(d\bs)\Big\}\cos\Big(tx-a\int_{S_2}(ts_1+(r+h)s_2)^{<\alpha>}\Gamma(ds)\Big)\\
		& \hspace{6cm} \times \int_{S_2} (t+\frac{rs'_2}{s'_1})^{<\alpha-1>}s'_2{s'_1}^{<\alpha-1>}\Gamma(ds')dt \Bigg]\\
		& = \lim_{h\rightarrow0} \frac{1}{h} \int_{\mathbb{R}} \int_{S_2} \Bigg[ \exp\bigg\{-\int_{S_2}\bigg|\Big(t-\frac{hs'_2}{s'_1}\Big)s_1+(r+h)s_2\bigg|^{\alpha}\Gamma(d\bs)\bigg\}\\
		& \hspace{4cm} \times \cos\bigg(\Big(t-\frac{hs'_2}{s'_1}\Big)x-a\int_{S_2}\bigg(\Big(t-\frac{hs'_2}{s'_1}\Big)s_1+(r+h)s_2\bigg)^{<\alpha>}\Gamma(ds)\bigg)\\
		& \hspace{3cm} - \exp\Big\{-\int_{S_2}|ts_1+(r+h)s_2|^{\alpha}\Gamma(d\bs)\Big\}\cos\Big(tx-a\int_{S_2}(ts_1+(r+h)s_2)^{<\alpha>}\Gamma(ds)\Big)\Bigg]\\
		& \hspace{7cm} \times   \bigg(t+r\frac{s'_2}{s'_1}\bigg)^{<\alpha-1>}s'_2 {s'_1}^{<\alpha-1>}\Gamma(ds')dt\\
		& = \lim_{h\rightarrow0} \frac{1}{h} \int_{\mathbb{R}} \int_{S_2} \dfrac{1}{\frac{hs'_2}{s'_1}} \Bigg[ \cos\bigg(\Big(t-\frac{hs'_2}{s'_1}\Big)x-a\int_{S_2}\bigg(\Big(t-\frac{hs'_2}{s'_1}\Big)s_1+(r+h)s_2\bigg)^{<\alpha>}\Gamma(ds)\bigg)\\
		& \hspace{4cm} - \cos\Big(tx-a\int_{S_2}(ts_1+(r+h)s_2)^{<\alpha>}\Gamma(ds)\Big)\Bigg]\\
		& \hspace{3cm} \times \exp\Big\{-\int_{S_2}|ts_1+(r+h)s_2|^{\alpha}\Gamma(d\bs)\Big\}\bigg(t+r\frac{s'_2}{s'_1}\bigg)^{<\alpha-1>}{s'_2}^2 |s'_1|^{\alpha-2}\Gamma(ds')dt\\
		& + \lim_{h\rightarrow0} \frac{1}{h} \int_{\mathbb{R}} \int_{S_2} \dfrac{1}{\frac{hs'_2}{s'_1}} \Bigg[ \exp\bigg\{-\int_{S_2}\bigg|\Big(t-\frac{hs'_2}{s'_1}\Big)s_1+(r+h)s_2\bigg|^{\alpha}\Gamma(d\bs)\bigg\}- \exp\Big\{-\int_{S_2}|ts_1+(r+h)s_2|^{\alpha}\Gamma(d\bs)\Big\}\Bigg]\\
		& \hspace{4cm} \times \cos\bigg(\Big(t-\frac{hs'_2}{s'_1}\Big)x-a\int_{S_2}\bigg(\Big(t-\frac{hs'_2}{s'_1}\Big)s_1+(r+h)s_2\bigg)^{<\alpha>}\Gamma(ds)\bigg)\\
		& \hspace{5cm}\times \bigg(t+r\frac{s'_2}{s'_1}\bigg)^{<\alpha-1>}{s'_2}^2 |s'_1|^{\alpha-2}\Gamma(ds')dt\\
		& = K_{31} + K_{32}.
		\end{align*}
		The case of $K_{32}$ is similar to that of $I_{22}$ in \cite{ct94} (p.106-108), the dominated convergence theorem applies. We focus on $K_{31}$. Its integrand converges to
		\begin{align*}
		& \sin\Bigg(tx - a\int_{S_2}(ts_1+rs_2)^{<\alpha>}\Gamma(ds) \Bigg)\exp\Bigg\{-\int_{S_2}|ts_1+rs_2|^{\alpha}\Gamma(ds)\Bigg\}\\
		& \hspace{2cm} \times \Bigg(x-\alpha a\int_{S_2}|ts_1+rs_2|^{\alpha-1}s_1\Gamma(ds)\Bigg)\Bigg(\int_{S_2}(ts'_1+rs'_2)^{<\alpha-1>}{s'_2}^2 {s'_1}^{-1}\Gamma(ds')\Bigg).
		\end{align*}
		Using the mean value theorem and Lemma \ref{le:34} $(\iota\iota)$, we can bound the integrand of $K_{31}$ for any $|h|<|r|$ by
		\begin{align*}
		& |\sin(y)|e^{2|r|^\alpha\sigma_2^\alpha} e^{-|t|^\alpha\sigma_1^\alpha} \int_{S_2} \bigg|t+r\frac{s'_2}{s'_1}\bigg|^{\alpha-1}{s'_2}^2 |s'_1|^{\alpha-2}\\
		& \hspace{0.5cm}\times \Bigg|\dfrac{1}{\frac{hs'_2}{s'_1}}\Bigg|\Bigg|-\frac{hs'_2}{s'_1}x - a\int_{S_2} \bigg(\Big(t-\frac{hs'_2}{s'_1}\Big)s_1+(r+h)s_2\bigg)^{<\alpha>} - (ts_1+(r+h)s_2)^{<\alpha>}\Gamma(ds)\bigg)\Bigg| \Gamma(ds')\\
		& \hspace{0.5cm} \le e^{2|r|^\alpha\sigma_2^\alpha} e^{-|t|^\alpha\sigma_1^\alpha} \int_{S_2} \bigg|t+r\frac{s'_2}{s'_1}\bigg|^{\alpha-1}{s'_2}^2 |s'_1|^{\alpha-2} \Bigg(|x| + 2a\int_{S_2}\Big|t + (r+h)\frac{s_2}{s_1}\Big|^{\alpha-1}|s_1|\Gamma(ds) \Bigg)\Gamma(ds')\\
		& \hspace{0.5cm} \le |x|e^{2|r|^\alpha\sigma_2^\alpha} e^{-|t|^\alpha\sigma_1^\alpha} \int_{S_2} \bigg|t+r\frac{s'_2}{s'_1}\bigg|^{\alpha-1}{s'_2}^2 |s'_1|^{\alpha-2}\Gamma(ds')\\
		& \hspace{1cm} + 2ae^{2|r|^\alpha\sigma_2^\alpha} e^{-|t|^\alpha\sigma_1^\alpha} \int_{S_2}\int_{S_2}\bigg|t+r\frac{s'_2}{s'_1}\bigg|^{\alpha-1}\Big|t + (r+h)\frac{s_2}{s_1}\Big|^{\alpha-1}|s_1|{s'_2}^2 |s'_1|^{\alpha-2}\Gamma(ds)\Gamma(ds').
		\end{align*}
		The integrability with respect to $t$ of the first (resp. second) term is obtained in the same way as for \eqref{dem:bound_lem22} (resp. \eqref{dem:bdk2_29}) and concluding using \eqref{eq:nu_cond} with $\nu>2-\alpha$.
		Thus, the dominated convergence theorem applies to $K_{31}$, which finally shows that the dominated convergence theorem applies to $J_1'$. The other $J$'s can be treated in a similar fashion. \\

		\noindent \textbf{Case $\boldsymbol{\alpha\in(1,2)}$}\\
		
		\noindent After derivation, $\phi^{(1)}_{X_2|x}(r)$ is given by \eqref{eq:phi'} with functions $J$'s of the form
		\begin{align*}
		\int_{\mathbb{R}}e^{-\int_{S_2}|ts_1+rs_2|^{\alpha}\Gamma(ds)}\text{trig}\Bigg(tx-a\int_{S_2}|ts_1+rs_2|^{<\alpha>}\Gamma(ds)\Bigg)\int_{S_2}(ts_1+rs_2)^{<\alpha-1> \text{ or } \alpha-1}s_2\Gamma(ds)dt,
		\end{align*}
		which are similar to deal with. Consider for instance $J_1(r)$. It's derivative can be written as in \eqref{eq:j1}
		\begin{align*}
		J_1'(r) & = K_1 + K_2 + K_3.
		\end{align*}
		For the integrand of $K_1$, we can use \eqref{eq:e-min} and the triangle inequality to bound it by
		\begin{align*}
		\Gamma(S_2)e^{|2r|^\alpha\sigma_2^\alpha}e^{-2^{1-\alpha}\sigma_1^\alpha|t|^\alpha}\alpha(|t|^{\alpha-1}+2|r|^{\alpha-1}) \int_{S_2} |ts_1+rs_2|^{\alpha-1}|s_2|\Gamma(ds).
		\end{align*}
		Since $0<\alpha-1<1$, we can further bound it by
		\begin{align*}
		\Gamma(S_2)e^{|2r|^\alpha\sigma_2^\alpha}e^{-2^{1-\alpha}\sigma_1^\alpha|t|^\alpha}\alpha(|t|^{\alpha-1}+2|r|^{\alpha-1})^2,
		\end{align*}
		which is integrable with respect to $t$. The same bound can be obtained for the integrand of $K_2$ using the mean value theorem, \eqref{eq:mvtcrochets} and Lemma \ref{le:33}. As for $K_3$, there is no need to perform "appropriate integration by parts" since $0<\alpha-1<1$. Its integrand converges to 
		\begin{align*}
		(\alpha-1)\exp\Big\{-\int_{S_2}|ts_1+rs_2|^{\alpha}\Gamma(d\bs)\Big\}\cos\Big(tx-a\int_{S_2}(ts_1+rs_2)^{<\alpha>}\Gamma(ds)\Big) \int_{S_2} |ts_1+rs_2|^{\alpha-2}s_2^2\Gamma(ds).
		\end{align*}
		Using Lemmas \ref{le:33} and \ref{le:34} $(\iota\iota)$, it can be bounded for any $|h|<|r|$ by
		\begin{align*}
		\dfrac{2}{|h|}\Gamma(S_2)e^{|2r|^\alpha\sigma_2^\alpha}e^{-2^{1-\alpha}\sigma_1^\alpha|t|^\alpha}\int_{S_2} |ts_1+rs_2|^{\alpha-2}|hs_2|\Gamma(ds),\\
		\le \Gamma(S_2)e^{|2r|^\alpha\sigma_2^\alpha}e^{-2^{1-\alpha}\sigma_1^\alpha|t|^\alpha} \int_{S_2} \Big|t+\dfrac{rs_2}{s_1}\Big|^{\alpha-2} |s_1|^{\alpha-2}\Gamma(ds).
		\end{align*}
		We can show that this bound is integrable with respect to $t$ using Lemma \ref{le:310} with $\eta=\alpha-2$, $b=0$ and $p=0$, the fact that $\int_{\mathbb{R}}e^{-2^{1-\alpha}\sigma_1^\alpha |t|^\alpha} |t|^{\alpha-2}dt<+\infty$ for $\alpha\in(1,2)$ and \eqref{eq:nu_cond} with $\nu>2-\alpha$. The dominated convergence theorem thus applies and we get
		\begin{align}
		& \phi^{(2)}_{X_2|x}(r)  = \dfrac{-\alpha}{2\pi f_{X_1}(x)}\Bigg[ -\alpha\int_{\mathbb{R}} e^{-itx}\varphi_{\bX}(t,r) \Big(\int_{S_2}g_2(ts_1+rs_2)s_2\Gamma(ds)\Big)^2dt \nonumber\\
		& \hspace{5cm} + (\alpha-1)\int_{\mathbb{R}}e^{-itx}\varphi_{\bX}(t,r)\Big(\int_{S_2}g_3(ts_1+rs_2)s_2^2\Gamma(ds)\Big)dt \Bigg],\label{eq:phi''_alpsup1}
		\end{align}
		with $g_3(z) = |z|^{\alpha-2}-iaz^{<\alpha-2>}$ for $z\in\mathbb{R}$. Integrating by parts the terms $|ts_1+rs_2|^{<\alpha-2> \text{ or } \alpha-2}$ involved in the expression $\int_{\mathbb{R}}e^{-itx}\varphi_{\bX}(t,r)\Big(\int_{S_2}g_3(ts_1+rs_2)s_2^2\Gamma(ds)\Big)dt$ yields the expression \eqref{eq:d2} obtained in the case $\alpha\in(1/2,1)$. Hence, 
		the same functional form for the second order conditional moment \eqref{eq:scm_gen} in Theorem \ref{theo:cond_moments_stable_general_alpnot1} holds when $\alpha>1$.
		
		\subsubsection{Justifying inversion: Third derivative}

		Let $\alpha\in(1,2)$ and let \eqref{eq:nu_cond} hold with $\nu>3-\alpha$. Starting from the second derivative of $\phi^{(2)}_{X_2|x}(r)$ given at \eqref{eq:d2}, with obvious notations
		\begin{align*}
		\phi^{(2)}_{X_2|x}(r) = \dfrac{-\alpha}{2\pi f_{X_1}(x)} \Big[ixI_1(r)+\alpha(I_3(r)-I_2(r))\Big]
		\end{align*}
		
		On the one hand, it can be shown that the dominated convergence theorem applies to $I_1'$ using the usual arguments the fact that \eqref{eq:nu_cond} holds with $\nu>3-\alpha$. On the other hand, after some elementary manipulations, we get that
		\begin{align*}
		I_3-I_2 & = \int_{\mathbb{R}} e^{-itx +i a \int_{S_2}(ts_1+rs_2)^{<\alpha>}\Gamma(d\bs)} e^{-\int_{S_2}|ts_1+rs_2|^\alpha\Gamma(d\bs)}\\
		& \hspace{1cm} \times \int_{S_2}\int_{S_2} \Bigg\{(ts_1+rs_2)^{<\alpha-1>}(ts_1'+rs_2')^{<\alpha-1>} - a^2 |ts_1+rs_2|^{\alpha-1}|ts_1'+rs_2'|^{\alpha-1}\\
		& \hspace{2cm} -ia \bigg(|ts_1+rs_2|^{\alpha-1}(ts_1'+rs_2')^{<\alpha-1>}+(ts_1+rs_2)^{<\alpha-1>}|ts_1'+rs_2'|^{\alpha-1}\bigg)\Bigg\}\\ &\hspace{8cm}\times \Big[s_2^2s_1^{-1}s_1'-s_2s_2'\Big]\Gamma(d\bs)\Gamma(d\bs')dt
		\end{align*}
		The previous expression can be decomposed into terms of the form 
		\begin{align*}
		& \int_{\mathbb{R}}\int_{S_2}\int_{S_2}\text{trig}\bigg(-tx + a \int_{S_2}(ts_1+rs_2)^{<\alpha>}\Gamma(d\bs)\bigg)\\
		& \hspace{2cm}\times e^{-\int_{S_2}|ts_1+rs_2|^\alpha\Gamma(d\bs)}\\
		& \hspace{3cm}\times |ts_1+rs_2|^{<\alpha-1>\text{ or } \alpha-1}\hspace{0.4cm}\times\hspace{0.4cm}|ts_1'+rs_2'|^{<\alpha-1>\text{ or } \alpha-1}\\
		& \hspace{4cm}\times\Big[s_2^2s_1^{-1}s_1'-s_2s_2'\Big]  \Gamma(d\bs)\Gamma(d\bs')dt,
		\end{align*}
		where <<trig>> is to be replaced by a sine or cosine function. Each of these terms can be treated in a similar way to show that the dominated convergence theorem applies. We will consider
		\begin{align*}
		J(r) & = \int_{\mathbb{R}}\int_{S_2}\int_{S_2}\cos\bigg(tx - a \int_{S_2}(ts_1+rs_2)^{<\alpha>}\Gamma(d\bs)\bigg)e^{-\int_{S_2}|ts_1+rs_2|^\alpha\Gamma(d\bs)}\\
		& \hspace{2cm}\times |ts_1+rs_2|^{\alpha-1}(ts_1'+rs_2')^{<\alpha-1>}\Big[s_2^2s_1^{-1}s_1'-s_2s_2'\Big]  \Gamma(d\bs)\Gamma(d\bs')dt.
		\end{align*}
		We have
		\begin{align*}
		J'(r) & = \lim_{h\rightarrow0}\dfrac{1}{h}\int_{\mathbb{R}}\int_{S_2}\int_{S_2} \Bigg[\cos\bigg(tx-a \int_{S_2}(ts_1+(r+h)s_2)^{<\alpha>}\Gamma(d\bs)\bigg)\\
		& \hspace{6cm}-\cos\bigg(tx-a \int_{S_2}(ts_1+rs_2)^{<\alpha>}\Gamma(d\bs)\bigg)\Bigg]\\
		& \hspace{2cm}\times e^{-\int_{S_2}|ts_1+(r+h)s_2|^\alpha\Gamma(d\bs)}|ts_1+(r+h)s_2|^{\alpha-1}(ts_1'+(r+h)s_2')^{<\alpha-1>}\\
		& \hspace{3cm}\times\Big[s_2^2s_1^{-1}s_1'-s_2s_2'\Big]  \Gamma(d\bs)\Gamma(d\bs')dt\\
		& + \lim_{h\rightarrow0}\dfrac{1}{h}\int_{\mathbb{R}}\int_{S_2}\int_{S_2}\cos\bigg(tx-a \int_{S_2}(ts_1+rs_2)^{<\alpha>}\Gamma(d\bs)\bigg)\\
		& \hspace{3cm} \times\Bigg[e^{-\int_{S_2}|ts_1+(r+h)s_2|^\alpha\Gamma(d\bs)}-e^{-\int_{S_2}|ts_1+rs_2|^\alpha\Gamma(d\bs)}\Bigg]\\
		& \hspace{2cm}\times |ts_1+(r+h)s_2|^{\alpha-1}(ts_1'+(r+h)s_2')^{<\alpha-1>}\Big[s_2^2s_1^{-1}s_1'-s_2s_2'\Big]  \Gamma(d\bs)\Gamma(d\bs')dt\\
		& + \lim_{h\rightarrow0}\dfrac{1}{h}\int_{\mathbb{R}}\int_{S_2}\int_{S_2}\cos\bigg(tx-a \int_{S_2}(ts_1+rs_2)^{<\alpha>}\Gamma(d\bs)\bigg)e^{-\int_{S_2}|ts_1+rs_2|^\alpha\Gamma(d\bs)}\\
		& \hspace{3cm}\times\Bigg[ |ts_1+(r+h)s_2|^{\alpha-1}-|ts_1+rs_2|^{\alpha-1}\Bigg]\\
		& \hspace{4cm}\times(ts_1'+(r+h)s_2')^{<\alpha-1>}\Big[s_2^2s_1^{-1}s_1'-s_2s_2'\Big]  \Gamma(d\bs)\Gamma(d\bs')dt\\
		& + \lim_{h\rightarrow0}\dfrac{1}{h}\int_{\mathbb{R}}\int_{S_2}\int_{S_2}\cos\bigg(tx-a \int_{S_2}(ts_1+rs_2)^{<\alpha>}\Gamma(d\bs)\bigg)e^{-\int_{S_2}|ts_1+rs_2|^\alpha\Gamma(d\bs)}\\
		& \hspace{3cm}\times\Bigg[ (ts_1'+(r+h)s_2')^{<\alpha-1>}-(ts_1'+rs_2')^{<\alpha-1>}\Bigg]\\
		& \hspace{4cm}\times|ts_1+rs_2|^{\alpha-1}\Big[s_2^2s_1^{-1}s_1'-s_2s_2'\Big]  \Gamma(d\bs)\Gamma(d\bs')dt\\
		& := K_1 + K_2 + K_3 + K_4.
		\end{align*}
		We will show that we can apply the dominated convergence theorem to the $K_i$'s. Let us begin with $K_1$. Its integrand converges to
		\begin{align*}
		& \alpha a \int_{S_2\times S_2\times S_2} \sin\bigg(tx-a \int_{S_2}(ts_1+rs_2)^{<\alpha>}\Gamma(d\bs)\bigg)e^{-\int_{S_2}|ts_1+rs_2|^\alpha\Gamma(d\bs)}\\
		& \hspace{2cm}\times |ts_1+rs_2|^{\alpha-1}(ts_1'+rs_2')^{<\alpha-1>}|ts_1''+rs_2''|^{\alpha-1}s_2''\Big[s_2^2s_1^{-1}s_1'-s_2s_2'\Big]\Gamma(d\bs)\Gamma(d\bs')\Gamma(d\bs'').
		\end{align*}
		For any $h$, $|h|<|r|$, the integrand of $K_1$ can be bounded using the mean value theorem on the cosine and Lemma \ref{le:33} by
		\begin{align}\label{eq:K11}
		& \dfrac{|a|}{|h|}  \bigg|\int_{S_2}(ts_1+(r+h)s_2)^{<\alpha>}-(ts_1+rs_2)^{<\alpha>}\Gamma(d\bs)\bigg|e^{2^\alpha|r|^\alpha\sigma_2^\alpha}e^{-2^{1-\alpha}\sigma_1^\alpha|t|^{\alpha}}\nonumber\\
		& \hspace{1cm} \times\bigg|\int_{S_2}\int_{S_2}|ts_1+(r+h)s_2|^{\alpha-1}(ts_1'+(r+h)s_2')^{<\alpha-1>}\Big[s_2^2s_1^{-1}s_1'-s_2s_2'\Big]  \Gamma(d\bs)\Gamma(d\bs')\bigg|.
		\end{align}
		
		Hence, by inequality \eqref{eq:mvtcrochets} and given that $0<\alpha-1<1$, the quantity \eqref{eq:K11} can be bounded by
		\begin{align*}
		& \alpha|a| \Gamma(S_2)e^{2^\alpha|r|^\alpha\sigma_2^\alpha}e^{-2^{1-\alpha}\sigma_1^\alpha|t|^{\alpha}} (|t|^{\alpha-1} + 2|r|^{\alpha-1}) \\
		& \hspace{1cm} \times\bigg|\int_{S_2}\int_{S_2}|ts_1+(r+h)s_2|^{\alpha-1}(ts_1'+(r+h)s_2')^{<\alpha-1>}\Big[s_2^2s_1^{-1}s_1'-s_2s_2'\Big]  \Gamma(d\bs)\Gamma(d\bs')\bigg|\\
		& \le \alpha|a| \Gamma(S_2)e^{2^\alpha|r|^\alpha\sigma_2^\alpha}e^{-2^{1-\alpha}\sigma_1^\alpha|t|^{\alpha}} (|t|^{\alpha-1} + 2|r|^{\alpha-1})^3 \bigg(\Gamma(S_2)+\int_{S_2}|s_1|^{-1}\Gamma(d\bs)\bigg)\\
		& \le \text{const } e^{-2^{1-\alpha}\sigma_1^\alpha|t|^{\alpha}} (|t|^{\alpha-1} + 2|r|^{\alpha-1})^3,
		\end{align*}
		where const is a finite nonnegative constant because of \eqref{eq:nu_cond} with $\nu>3-\alpha>1$ and the fact that $\Gamma$ is a finite measure. This last bound, independent of $h$, is integrable with respect to $t$ on $\mathbb{R}$. The dominated convergence theorem applies to $K_1$. Consider now $K_2$. Its integrand converges to
		\begin{align}
		& \alpha \int_{S_2\times S_2\times S_2}\cos\bigg(tx-a \int_{S_2}(ts_1+rs_2)^{<\alpha>}\Gamma(d\bs)\bigg)e^{-\int_{S_2}|ts_1+rs_2|^\alpha\Gamma(d\bs)}\label{eq:mvt}\\
		& \hspace{1cm}\times |ts_1+rs_2|^{\alpha-1}(ts_1'+rs_2')^{<\alpha-1>}(ts_1''+rs_2'')^{<\alpha-1>}s_2''\Big[s_2^2s_1^{-1}s_1'-s_2s_2'\Big]  \Gamma(d\bs)\Gamma(d\bs')\Gamma(d\bs'')\nonumber
		\end{align}
		
		By \eqref{eq:e-min}, the integrand of $K_2$ can be bounded by
		\begin{align*}
		& \Gamma(S_2)e^{|2r|^\alpha\sigma_2^\alpha}e^{-2^{1-\alpha}\sigma_1^\alpha|t|^\alpha}\alpha(|t|^{\alpha-1}+2|r|^{\alpha-1})\\
		& \hspace{2cm} \bigg|\int_{S_2}\int_{S_2} |ts_1+(r+h)s_2|^{\alpha-1}(ts_1'+(r+h)s_2')^{<\alpha-1>}\Big[s_2^2s_1^{-1}s_1'-s_2s_2'\Big]  \Gamma(d\bs)\Gamma(d\bs')\bigg|
		\end{align*}
		Which can be further bounded by an integrable function of $t$ in a similar way as for the integrand of $K_1$. The dominated convergence theorem applies to $K_2$. Consider now $K_3$. Its integrand converges to
		\begin{align*}
		& (\alpha-1)\int_{S_2}\int_{S_2}\cos\bigg(tx-a \int_{S_2}(ts_1+rs_2)^{<\alpha>}\Gamma(d\bs)\bigg)e^{-\int_{S_2}|ts_1+rs_2|^\alpha\Gamma(d\bs)}\\
		& \hspace{3cm}\times (ts_1+rs_2)^{<\alpha-2>}(ts_1'+(r+h)s_2')^{<\alpha-1>}s_2\Big[s_2^2s_1^{-1}s_1'-s_2s_2'\Big]  \Gamma(d\bs)\Gamma(d\bs')
		\end{align*}
		Using Lemmas \ref{le:33}, \ref{le:34} $(\iota)$ and the triangle inequality, the integrand of $K_3$ can be bounded by
		\begin{align*}
		& \dfrac{1}{|h|}e^{|r|^\alpha\sigma_2^\alpha} e^{-2^{1-\alpha}\sigma_1^\alpha|t|^\alpha}\int_{S_2}\int_{S_2} |h s_2| |ts_1+rs_2|^{\alpha-2}|ts_1'+(r+h)s_2'|^{\alpha-1}\Big|s_2^2s_1^{-1}s_1'-s_2s_2'\Big|  \Gamma(d\bs)\Gamma(d\bs')\\
		& \le e^{|r|^\alpha\sigma_2^\alpha} \Gamma(S_2)\int_{S_2}  e^{-2^{1-\alpha}\sigma_1^\alpha|t|^\alpha} |ts_1+rs_2|^{\alpha-2} (|t|^{\alpha-1} + 2|r|^{\alpha-1}) \Big|1+|s_1|^{-1}\Big|\Gamma(d\bs)
		\end{align*}
		To show the integrability with respect to $t$ of the last bound we make use of Lemma \ref{le:310} with $\eta=\alpha-2$, $b=0, \alpha-1$ and $p=0$ and the fact that with $1<\alpha<2$, $\int_{\mathbb{R}} e^{-2^{1-\alpha}\sigma_1^\alpha|t|^\alpha}|t|^{\alpha-2}dt<+\infty$ and $\int_{\mathbb{R}} e^{-2^{1-\alpha}\sigma_1^\alpha|t|^\alpha}|t|^{2\alpha-3}dt<+\infty$
		\begin{align*}
		& e^{|r|^\alpha\sigma_2^\alpha}\Gamma(S_2)\int_{S_2}\Big|1+|s_1|^{-1}\Big| \int_{\mathbb{R}} e^{-2^{1-\alpha}\sigma_1^\alpha|t|^\alpha}|s_1|^{\alpha-2}\Big|t+r\dfrac{s_2}{s_1}\Big|^{\alpha-2} (|t|^{\alpha-1} + 2|r|^{\alpha-1})dt \Gamma(d\bs)\\
		& \hspace{0.75cm} \le e^{|r|^\alpha\sigma_2^\alpha}\Gamma(S_2)\int_{S_2} \Big|1+|s_1|^{-1}\Big||s_1|^{\alpha-2} \Bigg[\int_{\mathbb{R}}e^{-2^{1-\alpha}\sigma_1^\alpha|t|^\alpha} \bigg|\Big|t+r\dfrac{s_2}{s_1}\Big|^{\alpha-2} - |t|^{\alpha-2} + |t|^{\alpha-2}\bigg||t|^{\alpha-1}dt\\
		& \hspace{4.5cm} + 2|r|^{\alpha-1}\int_{\mathbb{R}} e^{-2^{1-\alpha}\sigma_1^\alpha|t|^\alpha}\bigg|\Big|t+r\dfrac{s_2}{s_1}\Big|^{\alpha-2} - |t|^{\alpha-2} + |t|^{\alpha-2}\bigg|dt\Bigg]\Gamma(d\bs)\\
		& \hspace{0.75cm} \le e^{|r|^\alpha\sigma_2^\alpha}\Gamma(S_2)\int_{S_2} \Big|1+|s_1|^{-1}\Big||s_1|^{\alpha-2} \Bigg[\int_{\mathbb{R}}e^{-2^{1-\alpha}\sigma_1^\alpha|t|^\alpha} \bigg|\Big|t+r\dfrac{s_2}{s_1}\Big|^{\alpha-2} - |t|^{\alpha-2}\bigg||t|^{\alpha-1}dt\\
		& \hspace{7cm} + 2|r|^{\alpha-1}\int_{\mathbb{R}}e^{-2^{1-\alpha}\sigma_1^\alpha|t|^\alpha} \bigg|\Big|t+r\dfrac{s_2}{s_1}\Big|^{\alpha-2} - |t|^{\alpha-2} \bigg|dt\\
		& \hspace{7cm} + \int_{\mathbb{R}} e^{-2^{1-\alpha}\sigma_1^\alpha|t|^\alpha}|t|^{2\alpha-3}dt\\
		& \hspace{7cm} + 2|r|^{\alpha-1}\int_{\mathbb{R}} e^{-2^{1-\alpha}\sigma_1^\alpha|t|^\alpha}|t|^{\alpha-2}dt
		\Bigg]\Gamma(d\bs)\\
		& \hspace{0.75cm} \le \text{const } \int_{S_2} \Big|1+|s_1|^{-1}\Big||s_1|^{\alpha-2}\Gamma(d\bs)\\
		& \hspace{0.75cm} \le \text{const }  \Big(\int_{S_2}|s_1|^{\alpha-2}\Gamma(d\bs)+\int_{S_2}|s_1|^{\alpha-3}\Gamma(d\bs)\Big),
		\end{align*}
		which is finite because of \eqref{eq:nu_cond} with $\nu>3-\alpha$. Hence, the dominated convergence theorem applies to $K_3$. The case of $K_4$ is similar, using Lemma \ref{le:34} $(\iota\iota)$ instead of $(\iota)$ to bound the term $\bigg|(ts_1'+(r+h)s_2')^{<\alpha-2>} - (ts_1'+rs_2')^{<\alpha-2>} \bigg|$. The dominated convergence theorem applies to all the $K_i$'s and we can invert the integration and derivation signs in $J'$.
		
		\subsubsection{A special manipulation to obtain the fourth derivative}
		
		Before derivating $\phi^{(3)}_{X_2|x}$, we follow the advice stated in \cite{ct98} (p.48) 
		and integrate by parts the terms containing $\int_{S_2}g_3(ts_1+rs_2)s_2^3s_1^{-1}\Gamma(d\bs)$ and $\int_{S_2}g_3(ts_1+rs_2)s_2^2\Gamma(d\bs)$, namely $I_1$, $I_6$ and $I_7$. This is done in order to guarantee the validity of the representation of the fourth derivative when \eqref{eq:nu_cond} holds for any $\nu>4-\alpha$. If we did not do this step first, the obtained fourth derivative would be valid only when \eqref{eq:nu_cond} holds with $\nu>5-\alpha$. We obtain
		\begin{align}
		\phi^{(3)}_{X_2|x}(r) & = \dfrac{-\alpha}{2\pi f_{X_1}(x)} \Bigg[i\alpha x\Big(I_{11}-I_2+I_{62}-2I_{72}\Big)-x^2I_{12}\nonumber\\
		& \hspace{3cm} +\alpha^2\Big(I_3-I_4-2I_{71}+I_{61}\Big)+\alpha(\alpha-1)\Big(I_5-I_{63}+2I_{73}\Big)\Bigg],\label{eq:d4mod}
		\end{align}
		where, in addition to $I_2$, $I_3$, $I_4$ and $I_5$ defined in the Lemma,
		\begin{align*}
		I_{11} & = \Lambda \bigg(\int_{S_2}g_2s_2^3s_1^{-2}\bigg)\bigg(\int_{S_2}g_2s_1\bigg), & I_{12} & = \Lambda \bigg(\int_{S_2}g_2s_2^3s_1^{-2}\bigg),\\
		I_{61} & = \Lambda \bigg(\int_{S_2}g_2s_2^3s_1^{-2}\bigg)\bigg(\int_{S_2}g_2s_1\bigg)^2, & I_{71} & = \Lambda \bigg(\int_{S_2}g_2s_2^2s_1^{-1}\bigg)\bigg(\int_{S_2}g_2s_1\bigg)\bigg(\int_{S_2}g_2s_2\bigg),\\
		I_{62} & = \Lambda \bigg(\int_{S_2}g_2s_2^3s_1^{-2}\bigg)\bigg(\int_{S_2}g_2s_1\bigg), & I_{72} & =  \Lambda \bigg(\int_{S_2}g_2s_2^2s_1^{-1}\bigg)\bigg(\int_{S_2}g_2s_2\bigg),\\
		I_{63} & = \Lambda \bigg(\int_{S_2}g_2s_2^3s_1^{-2}\bigg)\bigg(\int_{S_2}g_3s_1^2\bigg), & I_{73} & = \Lambda \bigg(\int_{S_2}g_2s_2^2s_1^{-1}\bigg)\bigg(\int_{S_2}g_3s_2s_1\bigg).
		\end{align*}
		Both justification and computation of the fourth derivative are obtained by starting from the above representation of the third derivative.
		
		\subsubsection{Justifying inversion: Fourth derivative}
		
		Showing that the dominated convergence theorem holds when differentiating \eqref{eq:d4mod} is the most delicate for the terms: $I_5$, $I_{63}$ and $I_{73}$ -the terms involving the function $g_3$, that is, $|ts_1+rs_2|$ to the power $\alpha-2$. Arguments and bounds that have already been encountered can be used for the other ones. 
		
		Let us show the dominated convergence theorem applies to $I_5$. The cases of $I_{63}$ and $I_{73}$ are similar. We decompose $I_5$ into terms of the form
		\begin{align*}
		& \int_{\mathbb{R}}\int_{S_2}\int_{S_2}\text{trig}\bigg(-tx+a\int_{S_2}(ts_1+rs_2)^{<\alpha>}\Gamma(d\bs)\bigg)e^{-\int_{S_2}|ts_1+rs_2|^{\alpha}\Gamma(d\bs)}\\
		& \hspace{3.5cm} \times|ts_1+rs_2|^{\alpha-1\text{ or }<\alpha-1>}|ts_1'+rs_2'|^{\alpha-2\text{ or }<\alpha-2>}s_2^2s_1^{-1}s_2's_1'\Gamma(d\bs)\Gamma(d\bs')dt.
		\end{align*}
		Consider for instance
		\begin{align*}
		J(r) & := \int_{\mathbb{R}}\int_{S_2}\int_{S_2}\text{cos}\bigg(-tx+a\int_{S_2}(ts_1+rs_2)^{<\alpha>}\Gamma(d\bs)\bigg)e^{-\int_{S_2}|ts_1+rs_2|^{\alpha}\Gamma(d\bs)}\\
		& \hspace{3.5cm} \times|ts_1+rs_2|^{\alpha-1}|ts_1'+rs_2'|^{\alpha-2}s_2^2s_1^{-1}s_2's_1'\Gamma(d\bs)\Gamma(d\bs')dt.
		\end{align*}
		We have
		\begin{align*}
		%
		J'(r) & = \lim_{h\rightarrow0}\dfrac{1}{h}\int_{\mathbb{R}}\int_{S_2}\int_{S_2} \bigg[|ts_1'+(r+h)s_2'|^{\alpha-2}-|ts_1'+rs_2'|^{\alpha-2}\bigg]|ts_1+(r+h)s_2|^{\alpha-1}\\
		& \hspace{4cm} \times\cos\bigg(-tx+a\int_{S_2}(ts_1+(r+h)s_2)^{<\alpha>}\Gamma(d\bs)\bigg)\\
		& \hspace{5cm} \times e^{-\int_{S_2}|ts_1+(r+h)s_2|^{\alpha}\Gamma(d\bs)} s_2^2s_1^{-1}s_2's_1'\Gamma(d\bs)\Gamma(d\bs')dt\\
		%
		& + \lim_{h\rightarrow0}\dfrac{1}{h}\int_{\mathbb{R}}\int_{S_2}\int_{S_2}|ts_1'+rs_2'|^{\alpha-2} \bigg[|ts_1+(r+h)s_2|^{\alpha-1}-|ts_1+rs_2|^{\alpha-1}\bigg]\\
		& \hspace{4cm} \times\cos\bigg(-tx+a\int_{S_2}(ts_1+(r+h)s_2)^{<\alpha>}\Gamma(d\bs)\bigg)\\
		& \hspace{5cm} \times e^{-\int_{S_2}|ts_1+(r+h)s_2|^{\alpha}\Gamma(d\bs)} s_2^2s_1^{-1}s_2's_1'\Gamma(d\bs)\Gamma(d\bs')dt\\
		%
		& + \lim_{h\rightarrow0}\dfrac{1}{h}\int_{\mathbb{R}}\int_{S_2}\int_{S_2}|ts_1'+rs_2'|^{\alpha-2}|ts_1+rs_2|^{\alpha-1}\\
		& \times  \Bigg[\cos\bigg(-tx+a\int_{S_2}(ts_1+(r+h)s_2)^{<\alpha>}\Gamma(d\bs)\bigg)-\cos\bigg(-tx+a\int_{S_2}(ts_1+rs_2)^{<\alpha>}\Gamma(d\bs)\bigg)\Bigg]\\
		& \hspace{5cm} \times e^{-\int_{S_2}|ts_1+(r+h)s_2|^{\alpha}\Gamma(d\bs)} s_2^2s_1^{-1}s_2's_1'\Gamma(d\bs)\Gamma(d\bs')dt\\
		%
		& + \lim_{h\rightarrow0}\dfrac{1}{h}\int_{\mathbb{R}}\int_{S_2}\int_{S_2}|ts_1'+rs_2'|^{\alpha-2}|ts_1+rs_2|^{\alpha-1}\\
		& \hspace{2cm} \times \cos\bigg(-tx+a\int_{S_2}(ts_1+rs_2)^{<\alpha>}\Gamma(d\bs)\bigg) \\
		& \hspace{3cm} \times \Bigg[e^{-\int_{S_2}|ts_1+(r+h)s_2|^{\alpha}\Gamma(d\bs)}-e^{-\int_{S_2}|ts_1+rs_2|^{\alpha}\Gamma(d\bs)}\Bigg] s_2^2s_1^{-1}s_2's_1'\Gamma(d\bs)\Gamma(d\bs')dt\\
		& := K_1 + K_2 + K_3 + K_4
		\end{align*}
		The integrand of $K_4$ can be bounded using inequality \eqref{eq:mvt}, \eqref{eq:e-min} and invoking Lemma \ref{le:310} and \eqref{eq:nu_cond} with $\nu>4-\alpha$. The integrand of $K_3$ can be bounded using \eqref{eq:mvtcrochets} Lemma \ref{le:33}, and concluding with Lemma \ref{le:310} and \eqref{eq:nu_cond} with $\nu>4-\alpha$. Focus now on $K_2$. Using Lemmas \ref{le:33} and \ref{le:34} $(\iota)$, its integrand can be bounded by
		\begin{align*}
		e^{|2r|^\alpha\sigma_2^\alpha} e^{-2^{1-\alpha}\sigma_1^\alpha|t|^\alpha}\Big|t +\dfrac{rs_2'}{s_1'}\Big|^{\alpha-2}\Big|t +\dfrac{rs_2}{s_1}\Big|^{\alpha-2}s_2^3|s_1|^{\alpha-3}|s_1'|^{\alpha-1}|s_2'|.
		\end{align*}
		The later bound does not depend on $h$ and can be shown to be integrable with respect to $t$ using \eqref{eq:nu_cond} with $\nu>4-\alpha$, Lemma \ref{le:cor31} with $\eta=\alpha-2$, $z_2=z_4=0$, $p=0$ and the fact that $\int_\mathbb{R}e^{-c|t|^\alpha}|t|^{2(\alpha-2)}<+\infty$ for $\alpha\in(3/2,2)$. Let us now turn to the term $K_1$ which is more intricate. Appropriate <<integration by parts>> is required. With the change of variable $t=t+\frac{hs_2'}{s_1'}$,
		\begin{align*}
		%
		K_1 & = \lim_{h\rightarrow0}\dfrac{1}{h}\int_{S_2}\int_{S_2}\int_{\mathbb{R}}\Bigg[e^{-\int_{S_2}\bigg|\Big(t-\frac{hs_2'}{s_1'}\Big)s_1+(r+h)s_2\bigg|^\alpha\Gamma(d\bs)}-e^{-\int_{S_2}|ts_1+(r+h)s_2|^\alpha\Gamma(d\bs)}\Bigg]\\
		& \hspace{3cm} \times \cos\Bigg(\Big(t-\dfrac{hs_2'}{s_1'}\Big)x-a\int_{S_2}\bigg(\Big(t-\dfrac{hs_2'}{s_1'}\Big)s_1+(r+h)s_2\bigg)^{<\alpha>}\Gamma(d\bs)\Bigg)\\
		& \hspace{4cm} \times  \bigg|\Big(t-\dfrac{hs_2'}{s_1'}\Big)s_1+(r+h)s_2\bigg|^{\alpha-1}|ts_1'+rs_2'|^{\alpha-2}s_2^2s_1^{-1}s_2's_1'dt\Gamma(d\bs)\Gamma(d\bs')\\
		& + \lim_{h\rightarrow0}\dfrac{1}{h}\int_{S_2}\int_{S_2}\int_{\mathbb{R}}e^{-\int_{S_2}|ts_1+(r+h)s_2|^\alpha\Gamma(d\bs)}\\
		& \hspace{2cm}\times \cos\Bigg(\Big(t-\dfrac{hs_2'}{s_1'}\Big)x-a\int_{S_2}\bigg(\Big(t-\dfrac{hs_2'}{s_1'}\Big)s_1+(r+h)s_2\bigg)^{<\alpha>}\Gamma(d\bs)\Bigg)\\
		& \hspace{3cm}\times  \Bigg[\bigg|\Big(t-\dfrac{hs_2'}{s_1'}\Big)s_1+(r+h)s_2\bigg|^{\alpha-1}-\bigg|ts_1+(r+h)s_2\bigg|^{\alpha-1}\Bigg]\\
		& \hspace{5cm}\times |ts_1'+rs_2'|^{\alpha-2}s_2^2s_1^{-1}s_2's_1'dt\Gamma(d\bs)\Gamma(d\bs')\\
		%
		& + \lim_{h\rightarrow0}\dfrac{1}{h}\int_{S_2}\int_{S_2}\int_{\mathbb{R}}e^{-\int_{S_2}|ts_1+(r+h)s_2|^\alpha\Gamma(d\bs)}\\
		& \hspace{2cm} \times \Bigg[\cos\Bigg(\Big(t-\dfrac{hs_2'}{s_1'}\Big)x-a\int_{S_2}\bigg(\Big(t-\dfrac{hs_2'}{s_1'}\Big)s_1+(r+h)s_2\bigg)^{<\alpha>}\Gamma(d\bs)\Bigg)\\
		& \hspace{6cm} -\cos\Bigg(tx-a\int_{S_2}\bigg(ts_1+(r+h)s_2\bigg)^{<\alpha>}\Gamma(d\bs)\Bigg)\Bigg]\\
		& \hspace{4cm}\times  \big|ts_1+(r+h)s_2\big|^{\alpha-1}|ts_1'+rs_2'|^{\alpha-2}s_2^2s_1^{-1}s_2's_1'dt\Gamma(d\bs)\Gamma(d\bs')\\
		& := K_{11}+K_{12}+K_{13}.
		\end{align*}
		It can be shown that the generalised Lebesgue convergence theorem applies to the terms $K_{11}$ and $K_{12}$ following the proof in \cite{ct98} (p.50-52). Regarding the integrand of $K_{13}$, using the mean value theorem on the cosine, Lemma \ref{le:33} and \eqref{eq:mvtcrochets}, we get for $|h|<|r|$
		\begin{align*}
		& \dfrac{1}{\left|\frac{hs_2'}{s_1'}\right|}e^{|2r|^\alpha\sigma_2^\alpha}e^{-2^{1-\alpha}\sigma_1^\alpha|t|^\alpha} \big|ts_1+(r+h)s_2\big|^{\alpha-1}|ts_1'+rs_2'|^{\alpha-2}s_2^2|s_1|^{-1}|s_2'|^2\\
		& \hspace{2cm} \times \Bigg|\dfrac{hs_2'}{s_1'}x+a\int_{S_2}\bigg(\Big(t-\dfrac{hs_2'}{s_1'}\Big)s_1+(r+h)s_2\bigg)^{<\alpha>}-\bigg(ts_1+(r+h)s_2\bigg)^{<\alpha>}\Gamma(d\bs)\Bigg|\\
		& \le \dfrac{1}{\left|\frac{hs_2'}{s_1'}\right|}e^{|2r|^\alpha\sigma_2^\alpha}e^{-2^{1-\alpha}\sigma_1^\alpha|t|^\alpha} \big|ts_1+(r+h)s_2\big|^{\alpha-1}|ts_1'+rs_2'|^{\alpha-2}s_2^2|s_1|^{-1}|s_2'|^2\\
		& \hspace{2cm} \times \Bigg[\Big|\dfrac{hs_2'}{s_1'}x\Big|+\Big|a\dfrac{hs_2'}{s_1'}\Big|\int_{S_2}|s_1| |ts_1+(r+h)s_2|^{\alpha-1}\Gamma(d\bs)\Bigg]\\
		& \le e^{|2r|^\alpha\sigma_2^\alpha}e^{-2^{1-\alpha}\sigma_1^\alpha|t|^\alpha} \Big|t+\dfrac{rs_2'}{s_1'}\Big|^{\alpha-2}s_2^2|s_1|^{-1}{s_2'}^2|s_1'|^{\alpha-2}\\
		& \hspace{2cm} \times \big(|t|^{\alpha-1}+|2r|^{\alpha-1}\big)\Bigg[|x|+|a|\Gamma(S_2)(|t|^{\alpha-1}+|2r|^{\alpha-1})\Bigg].
		\end{align*}
		The last bound can be shown to be integrable with respect to $t$ using Lemma \ref{le:312} with $\eta=\alpha-2$, $b=0, \alpha-1, 2(\alpha-1)$, $p=0$ and \eqref{eq:nu_cond} with $\nu>4-\alpha$.
		We established that we can invert the derivation and integration signs in all the $K_i$'s, hence in $J'$.
		
		\subsubsection{Lemmas for justifying the inversions in the proof of Lemma  \ref{lem:derivatives}}
		
		\noindent The following elementary lemmas, stated without proof, are used to establish Lemma \ref{lem:derivatives}.
		
		\begin{lem}\label{le:exp}
			For $x,y\in\mathbb{R}$,
			\begin{align}
			|e^{-x}-e^{-y}| \le e^{-\min(x,y)}|x-y|, \label{eq:exp95}\\
			|e^{-x}-e^{-y}| \le e^{-y}e^{|x-y|}|x-y|. \label{eq:exp94}
			\end{align}
		\end{lem}
		
		\begin{lem}\label{le:ti}
			For $\alpha>1$ and $x,y\in\mathbb{R}$,
			\begin{align*}
			\max\Big(2^{1-\alpha}|x|^\alpha-|y|^\alpha,2^{1-\alpha}|y|^\alpha-|x|^\alpha\Big)\le|x+y|^\alpha \le 2^{\alpha-1} \Big(|x|^\alpha+|y|^\alpha\Big).
			\end{align*}
		\end{lem}

		\begin{lem}\label{le:34} 
			For $z\in\mathbb{R}$ and $0<b\le 1$,
			\begin{align*}
			& (\iota) \hspace{1cm} \Big||1+z|^b-1\Big| \le |z|,\\
			& (\iota\iota) \hspace{0.875cm} \Big||1+z|^{<b>}-1\Big| \le 2|z|.
			\end{align*}
		\end{lem}

		\begin{lem}[Lemma 3.3, Cioszek-Georges and Taqqu (1998)]\label{le:33}
			For $\alpha>1$ and $t,r\in\mathbb{R}$,
			\begin{align*}
			\exp\Big\{-\int_{S_2}|ts_1+rs_2|^\alpha\Gamma(d\bs)\Big\} \le \exp\{|r|^\alpha\sigma_2^\alpha\}\exp\{-2^{1-\alpha}\sigma_1^\alpha|t|^\alpha\}.
			\end{align*}
		\end{lem}
		
		\begin{lem}[Lemma 3.1, Cioszek-Georges and Taqqu (1998)]\label{le:310}
			The following inequality holds for $c>0$, $0<\alpha<2$, $-1<\eta<0$ and $-1-\eta<b$:
			\begin{align*}
			\int_{\mathbb{R}} \exp(-c|t|^\alpha) \Big||t+z|^\eta -|t|^\eta\Big||t|^b dt \le const. \, |z|^p
			\end{align*}
			with
			$$
			0\le p <b + \eta + 1 \quad \text{for} \quad -1-\eta < b < 0,
			$$
			and
			$$
			0 \le p < \eta +1 \quad \text{or} \quad b \le p < b + \eta +\eta +1, p \le 1 \quad \text{for} \quad 0 \le b.
			$$
			const. depends only on $c$, $\alpha$, $\eta$, $b$ and $p$.
		\end{lem}
		
		\begin{lem}[Corollary 3.1, Cioszek-Georges and Taqqu (1998)]\label{le:cor31}
			The following inequality holds for $c>0$, $0<\alpha<2$, $-1/2<\eta<0$ and $0\le p <2\eta+1$:
			\begin{align*}
			\int_{\mathbb{R}} \exp(-c|t|^\alpha) \Big||t+z_1|^\eta|t+z_3|^\eta -|t+z_2|^\eta|t+z_4|^\eta\Big| dt \le const. \, (|z_1-z_2|^p+|z_3-z_4|^p),
			\end{align*}
			where const depends only on $c$, $\alpha$, $\eta$ and $p$.
		\end{lem}
		
		\begin{lem}[Lemma 3.12, Cioszek-Georges and Taqqu (1998)]\label{le:312}
			The following inequality holds for $c>0$, $0<\alpha<2$, $-1<\eta<0$, $b\ge0$ and $0\le p <\eta+1$:
			\begin{align*}
			\int_{\mathbb{R}} \exp(-c|t|^\alpha) \Big||t+z_1|^\eta -|t+z_2|^\eta\Big||t|^b dt \le const. \, |z_1-z_2|^p,
			\end{align*}
			where const depends only on $c$, $\alpha$, $\eta$, $b$ and $p$.
		\end{lem}
		
		\subsection{Computation of the derivatives}
		
		We detail the computation of the second order derivative highlighting where appropriate integration by parts intervenes.
		The computations are similar for the third and fourth order derivatives.
		
		
		Note that if $f(x)= |x|^b$, for $x,b\in\mathbb{R}$, $b\ne0$, then for $x\ne0$, $f'(x) = b x^{<b-1>}$ and if $f:x\longmapsto x^{<b>}$, then $f'(x) = b |x|^{b-1}$. This can be shown by distinguishing the cases $x>0$ and $x<0$. 

		\begin{align*}
		\phi^{(2)}_{X_2|x}(r) & = \dfrac{\partial }{\partial r}\phi^{(1)}_{X_2|x}(r)\\
		& = \dfrac{-\alpha}{2\pi f_{X_1}(x)}\lim_{h\rightarrow0}\dfrac{1}{h} \Bigg[\int_{\mathbb{R}}\int_{S_2}e^{-itx}\varphi_{\bX}(t,r+h)g_2(ts_1+(r+h)s_2)s_2\Gamma(d\bs)dt\\
		& \hspace{6cm} -\int_{\mathbb{R}}\int_{S_2}e^{-itx}\varphi_{\bX}(t,r)g_2(ts_1+rs_2)s_2\Gamma(d\bs)dt\Bigg]\\
		& = \dfrac{-\alpha}{2\pi f_{X_1}(x)}\lim_{h\rightarrow0}\dfrac{1}{h} \int_{\mathbb{R}}\int_{S_2}e^{-itx}\bigg[\varphi_{\bX}(t,r+h)-\varphi_{\bX}(t,r)\bigg]g_2(ts_1+(r+h)s_2)s_2\Gamma(d\bs)dt\\
		& \hspace{1cm} + \dfrac{-\alpha}{2\pi f_{X_1}(x)}\lim_{h\rightarrow0}\dfrac{1}{h} \int_{\mathbb{R}}\int_{S_2}e^{-itx}\varphi_{\bX}(t,r)\bigg[g_2(ts_1+(r+h)s_2)-g_2(ts_1+rs_2)\bigg]s_2\Gamma(d\bs)dt\\
		& := A_1 + A_2.
		\end{align*}
		The first limit can be straightforwardly obtained:
		\begin{align*}
		A_1 & = \dfrac{\alpha^2}{2\pi f_{X_1}(x)}\int_{\mathbb{R}} e^{-itx}\varphi_{\bX}(t,r)\bigg(\int_{S_2}g_2(ts_1+rs_2)s_2\Gamma(d\bs)\bigg)^2dt\\
		& = \dfrac{\alpha^2}{2\pi f_{X_1}(x)} \Lambda \bigg(\int_{S_2}g_2s_2\bigg)^2.
		\end{align*}
		The second one requires appropriate integration by parts. With the change of variable $t' = t + \dfrac{hs_2}{s_1}$,
		\begin{align*}
		A_2 & = \dfrac{-\alpha}{2\pi f_{X_1}(x)}\lim_{h\rightarrow0}\dfrac{1}{h} \bigg[\int_{S_2}\int_{\mathbb{R}}e^{-itx}\varphi_{\bX}(t,r)g_2(ts_1+(r+h)s_2)s_2dt\Gamma(d\bs)\\
		& \hspace{6cm}-\int_{S_2}\int_{\mathbb{R}}e^{-itx}\varphi_{\bX}(t,r)g_2(ts_1+rs_2)s_2dt\Gamma(d\bs)\bigg]\\
		& = \dfrac{-\alpha}{2\pi f_{X_1}(x)}\lim_{h\rightarrow0} \dfrac{1}{h}\bigg[\int_{S_2}\int_{\mathbb{R}}e^{-i\Big(t-\dfrac{hs_2}{s_1}\Big)x}\varphi_{\bX}\Big(t-\dfrac{hs_2}{s_1},r\Big)g_2(ts_1+rs_2)s_2dt\Gamma(d\bs)\\
		& \hspace{6cm}-\int_{S_2}\int_{\mathbb{R}}e^{-itx}\varphi_{\bX}(t,r)g_2(ts_1+rs_2)s_2dt\Gamma(d\bs)\bigg]\\
		& = \dfrac{\alpha}{2\pi f_{X_1}(x)}\int_{S_2}\int_{\mathbb{R}}s_2^2s_1^{-1}g_2(ts_1+rs_2)\lim_{h\rightarrow0} \dfrac{1}{-\frac{hs_2}{s_1}}\bigg[e^{-i\Big(t-\dfrac{hs_2}{s_1}\Big)x}\varphi_{\bX}\Big(t-\dfrac{hs_2}{s_1},r\Big)\\
		&  \hspace{11cm}-e^{-itx}\varphi_{\bX}(t,r)\bigg]dt\Gamma(d\bs)
		\end{align*}
		\begin{align*}
		& = \dfrac{\alpha}{2\pi f_{X_1}(x)}\int_{S_2}\int_{\mathbb{R}}s_2^2s_1^{-1}g_2(ts_1+rs_2)\bigg[-ixe^{-itx}\varphi_{\bX}(t,r)+e^{-itx}\dfrac{\partial}{\partial t}\varphi_{\bX}(t,r)\bigg]dt\Gamma(d\bs)\\
		& = \dfrac{-i\alpha x}{2\pi f_{X_1}(x)}\int_{\mathbb{R}}e^{-itx}\varphi_{\bX}(t,r)\bigg(\int_{S_2}s_2^2s_1^{-1}g_2(ts_1+rs_2)\Gamma(d\bs)\bigg)dt\\
		& \hspace{1cm} - \dfrac{\alpha^2}{2\pi f_{X_1}(x)}\int_{\mathbb{R}}e^{-itx}\varphi_{\bX}(t,r)\bigg(\int_{S_2}s_1g_2(ts_1+rs_2)\Gamma(d\bs)\bigg)\bigg(\int_{S_2}s_2^2s_1^{-1}g_2(ts_1+rs_2)\Gamma(d\bs)\bigg)dt\\
		A_2 & = \dfrac{-i\alpha x}{2\pi f_{X_1}(x)}\Lambda\bigg(\int_{S_2}g_2s_2^2s_1^{-1}\bigg) - \dfrac{\alpha^2}{2\pi f_{X_1}(x)}\Lambda\bigg(\int_{S_2}g_2s_2^2s_1^{-1}\bigg)\bigg(\int_{S_2}g_2s_1\bigg)
		\end{align*}
		Combining the expressions obtained for $A_1$ and $A_2$ yields the second derivative.

		\section{Proof of Theorem \ref{theo:cond_moments_stable_general_alpnot1}}
		\label{sec:prooftheo22}
		
		We here finally evaluate the derivatives of Lemma \ref{lem:derivatives} at $r=0$ to obtain the functional forms of the conditinal moments. 
		These proofs yield in particular the expressions of the constants $\btheta_i$, $i=1,\ldots,6$ which intervene in Theorem \ref{theo:cond_moments_stable_general_alpnot1}.
		Lemmas at the end of this section are used to regroup terms and simplify as much as possible the functional forms.
		
		\subsection{Proof of second order conditional moment \eqref{eq:scm_gen} in Theorem \ref{theo:cond_moments_stable_general_alpnot1}}
		\label{sec:tha3}
		
		The second order derivative of the characteristic function of $X_2|X_1=x$ is given by \eqref{eq:d2} in Lemma \ref{lem:derivatives}. Evaluating it at $r=0$ yields
		\begin{align*}
		& \mathbb{E}\Big[X_2^2\Big|X_1=x\Big] \\
		& \hspace{1cm} = -\phi^{(2)}_{X_2|x}(0)\\
		%
		& \hspace{1cm} = \dfrac{\alpha}{2\pi f_{X_1}(x)}\int_{\mathbb{R}} e^{-itx + ia\sigma_1^\alpha\beta_1t^{<\alpha>}} e^{-\sigma_1^\alpha|t|^\alpha}\\
		& \hspace{4cm}\times\Bigg[ ix\sigma_1^\alpha(\kappa_2 t^{<\alpha-1>}-ia\lambda_2|t|^{\alpha-1}) - \alpha\sigma_1^{2\alpha}(\kappa_1t^{<\alpha-1>}-ia\lambda_1|t|^{\alpha-1})^2\\
		%
		& \hspace{5cm}+ \alpha\sigma_1^{2\alpha}(\kappa_2t^{<\alpha-1>}-ia\lambda_2|t|^{\alpha-1})(t^{<\alpha-1>}-ia\beta_1|t|^{\alpha-1}) \Bigg]dt\\
		%
		& \hspace{1cm} = \dfrac{\alpha\sigma_1^{\alpha}}{2\pi f_{X_1}(x)}\int_{\mathbb{R}} e^{-itx + ia\sigma_1^\alpha\beta_1t^{<\alpha>}} e^{-\sigma_1^\alpha|t|^\alpha} \\
		& \hspace{4cm} \times \Bigg[ xa\lambda_2|t|^{\alpha-1} + \alpha\sigma_1^{\alpha}|t|^{2(\alpha-1)} \bigg(\kappa_2-a^2\beta_1\lambda_2-\kappa_1^2 + a^2\lambda_1^2\bigg)\\
		%
		& \hspace{5cm} + ix\kappa_2t^{<\alpha-1>} + i\alpha\sigma_1^{\alpha} t^{<2(\alpha-1)>}\bigg(2a\lambda_1\kappa_1-a(\lambda_2+\beta_1\kappa_2\bigg)\Bigg]dt\\
		& \hspace{1cm} = \dfrac{\alpha\sigma_1^{\alpha}}{\pi f_{X_1}(x)}\Bigg[ax\lambda_2C_1(x) + \kappa_2xS_1(x)\\
		& \hspace{3cm} -\alpha\sigma_1^{\alpha}\Big(\kappa_1^2-a^2\lambda_1^2 + a^2\beta_1\lambda_2 - \kappa_2\Big)C_2(x)-\alpha\sigma_1^{\alpha}\Big(a(\lambda_2+\beta_1\kappa_2)-2a\lambda_1\kappa_1\Big)S_2(x)\Bigg],
		\end{align*}
		where the $\kappa_i$'s and $\lambda_i$'s are given in \eqref{def:sbkl}. 
		Invoking Lemma \ref{le:SCFG} $(\iota\iota\iota)$ yields
		\begin{align*}
		\mathbb{E}\Big[X_2^2\Big|X_1=x\Big] & = \frac{x}{1+(a\beta_1)^2} \Bigg[(a^2\lambda_2\beta_1+\kappa_2)x + a(\lambda_2-\kappa_2\beta_1)\dfrac{1-xH(x)}{\pi f_{X_1}(x)}\Bigg]\\
		& \hspace{0.5cm} - \dfrac{\alpha^2\sigma_1^{2\alpha}}{\pi f_{X_1}(x)}  \mathcal{H}\Big(2(\alpha-1),\btheta_1;x\Big)\\
		& = \kappa_2 x^2  + \frac{ax(\lambda_2-\beta_1\kappa_2)}{1+(a\beta_1)^2} \Bigg[a\beta_1x + \dfrac{1-xH(x)}{\pi f_{X_1}(x)}\Bigg] - \dfrac{\alpha^2\sigma_1^{2\alpha}}{\pi f_{X_1}(x)}  \mathcal{H}\Big(2(\alpha-1),\btheta_1;x\Big),
		\end{align*}
		where $\mathcal{H}$ is given in \eqref{def:hcal} with
		\begin{align*}
		\theta_{11} & = \kappa_1^2-a^2\lambda_1^2 + a^2\beta_1\lambda_2 - \kappa_2,
		& \theta_{12} & = a(\lambda_2+\beta_1\kappa_2)-2a\lambda_1\kappa_1.
		\end{align*}
		
		\subsection{Proof of third order conditional moment \eqref{eq:tcm_gen} in Theorem \ref{theo:cond_moments_stable_general_alpnot1}}
		\label{sec:tha5}
		
		The third order derivative of the characteristic function of $X_2|X_1=x$ is given by \eqref{eq:d3} in Lemma \ref{lem:derivatives}. It can be shown that the $I$'s evaluated at $r=0$ write
		\begin{align*}
		I_1 & = 2\sigma_1^\alpha \mathcal{H}\Big(\alpha-2,\btheta_{1}^I;x\Big), & \btheta_1^I & = \bigg(\kappa_3,-a\lambda_3\bigg),\\
		I_2 & = 2\sigma_1^{2\alpha} \mathcal{H}\Big(2(\alpha-1),\btheta_{2}^I;x\Big), & \btheta_2^I & = \bigg(L,-aK\bigg),\\
		iI_3 & = 2\sigma_1^{3\alpha} \mathcal{H}\Big(3(\alpha-1),\btheta_{3}^I;x\Big), & \btheta_3^I & = \bigg(a\lambda_1(3\kappa_1^2-a^2\lambda_1^2),\kappa_1^3-3a^2\kappa_1\lambda_1^2\bigg),\\
		iI_4 & = 2\sigma_1^{3\alpha} \mathcal{H}\Big(3(\alpha-1),\btheta_{4}^I;x\Big), & \btheta_4^I & = \bigg(a\Big(K+\beta_1L\Big),L-a^2\beta_1K\bigg),\\
		iI_5 = iI_7& = 2\sigma_1^{2\alpha} \mathcal{H}\Big(2\alpha-3,\btheta_{5}^I;x\Big), & \btheta_5^I & = \bigg(aK,L\bigg),\\
		iI_6 & = 2\sigma_1^{2\alpha} \mathcal{H}\Big(2\alpha-3,\btheta_{6}^I;x\Big), & \btheta_6^I & = \bigg(a(\lambda_3+\beta_1\kappa_3),\kappa_3-a^2\beta_1\lambda_3\bigg),
		\end{align*}
		with $K=\kappa_1\lambda_2+\lambda_1\kappa_2$ and $L=\kappa_1\kappa_2-a^2\lambda_1\lambda_2$. Hence,
		\begin{align*}
		\mathbb{E}\Big[X_2^3\Big|X_1=x\Big] = -i\phi^{(3)}_{X_2|x}(0) = \dfrac{\alpha}{\pi f_{X_1}(x)} \Bigg[-x\Big((\alpha-1)K_1-\alpha K_2\Big) + \alpha^2K_3 + \alpha(\alpha-1)K_4\Bigg],
		\end{align*}
		with
		\begin{align*}
		K_1 & = \sigma_1^\alpha \mathcal{H}\Big(\alpha-2,\btheta_{1}^K;x\Big), & \text{with} \quad \btheta_{1}^K & = \btheta_{1}^I,\\
		K_2 & = \sigma_1^{2\alpha} \mathcal{H}\Big(2(\alpha-1),\btheta_{2}^K;x\Big), & \text{with} \quad \btheta_{2}^K & = \btheta_{2}^I,\\
		K_3 & = \sigma_1^{3\alpha} \mathcal{H}\Big(3(\alpha-1),\btheta_{3}^K;x\Big), & \text{with} \quad \btheta_{3}^K & = \btheta_{3}^I - \btheta_{4}^I\\
		K_4 & = \sigma_1^{2\alpha} \mathcal{H}\Big(2\alpha-3,\btheta_{4}^K;x\Big), & \text{with} \quad \btheta_{4}^K & = \btheta_{6}^I - \btheta_{5}^I.
		\end{align*}
		Invoking Lemma \ref{le:SCFG} ($\iota\iota$) for $n=1,2$ and regrouping the terms, we get
		\begin{align*}
		\mathbb{E}\Big[X_2^3\Big|X_1=x\Big] & = \dfrac{\alpha x^2\sigma_1^\alpha}{\pi f_{X_1}(x)} \bigg(\theta_{12}^KC_1(x)-\theta_{11}^KS_1(x)\bigg)\\
		& \hspace{0.5cm} +\dfrac{\alpha}{\pi f_{X_1}(x)} \Bigg[\hspace{0.3cm}\dfrac{\alpha x \sigma_1^{2\alpha}}{2}C_2(x)\bigg(-2\Big(\theta_{11}^K+a\beta_1\theta_{12}^K\Big)+2\theta_{21}^K-\theta_{42}^K\bigg)\\
		& \hspace{2.5cm} + \dfrac{\alpha x \sigma_1^{2\alpha}}{2}S_2(x) \bigg(-2\Big(\theta_{12}^K-a\beta_1\theta_{11}^K\Big)+2\theta_{22}^K+\theta_{41}^K\bigg)\\
		& \hspace{2.5cm} + \dfrac{\alpha^2 \sigma_1^{3\alpha}}{2}C_3(x) \bigg(2\theta_{31}^K + \theta_{41}^K + a\beta_1\theta_{42}^K\bigg)\\
		& \hspace{2.5cm} + \dfrac{\alpha^2 \sigma_1^{3\alpha}}{2}S_3(x) \bigg(2\theta_{32}^K + \theta_{42}^K - a\beta_1\theta_{41}^K\bigg)\Bigg].
		\end{align*}
		Using Lemma \ref{le:SCFG} ($\iota\iota\iota$) yields the conclusion with $\btheta_2=(\theta_{21},\theta_{22})$, $\btheta_3=(\theta_{31},\theta_{32})$ such that 
		\begin{align}
		\theta_{21} & = 3(L+a^2\beta_1\lambda_3-\kappa_3), \label{eq:th2}\\
		\theta_{22}  & = 3a(\lambda_3+\beta_1\kappa_3-K),\\
		\theta_{31} & = a\Big(\lambda_3(1-a^2\beta_1^2) +2\beta_1\kappa_3 + 2\lambda_1(3\kappa_1^2-a^2\lambda_1^2)-3(K+\beta_1L)\Big),\\
		\theta_{32}  & = \kappa_3(1-a^2\beta_1^2)-2a^2\beta_1\lambda_3+2(\kappa_1^3-3a^2\kappa_1\lambda_1^2) + 3(a^2\beta_1K-L),\label{eq:th3}
		\end{align}
		with $K=\kappa_1\lambda_2+\kappa_2\lambda_1$, $L=\kappa_1\kappa_2-a^2\lambda_1\lambda_2$.
		
		%
		
		\subsection{Proof of fourth order conditional moment \eqref{eq:fcm_gen} in Theorem \ref{theo:cond_moments_stable_general_alpnot1}}
		
		The conditional moments are obtained by evaluating the derivatives of the conditional characteristic function at $r=0$. We provide here the proof for the fourth order, which yields the expressions of the vectors $\btheta_4$, $\btheta_5$ and $\btheta_6$ appearing in Equation \eqref{eq:fcm_gen} of Theorem \ref{theo:cond_moments_stable_general_alpnot1}. The fourth order derivative of the characteristic function of $X_2|X_1=x$ is given by \eqref{eq:d4} in Lemma \ref{lem:derivatives}. It can be shown that the $J$'s evaluated at $r=0$ write
		\begin{align*}
		iJ_1 & = 2\sigma_1^{3\alpha} \mathcal{H}\Big(3(\alpha-1),\btheta^J_1;x\Big), & J_{11} = J_{13} & = 2\sigma_1^{3\alpha} \mathcal{H}\Big(3\alpha-4,\btheta^J_{11};x\Big),\\
		iJ_2 & = 2\sigma_1^{3\alpha} \mathcal{H}\Big(3(\alpha-1),\btheta^J_2;x\Big), & J_{14} & = 2\sigma_1^{2\alpha} \mathcal{H}\Big(2\alpha-4,\btheta^J_{14};x\Big),\\
		iJ_3 & = 2\sigma_1^{2\alpha} \mathcal{H}\Big(2\alpha-3,\btheta^J_3;x\Big), & J_{15} & =
		2\sigma_1^{2\alpha} \mathcal{H}\Big(2\alpha-4,\btheta^J_{15};x\Big),\\
		iJ_4 = iJ_5 & = 2\sigma_1^{2\alpha} \mathcal{H}\Big(2\alpha-3,\btheta^J_4;x\Big), & J_{16} & = 2\sigma_1^{2\alpha} \mathcal{H}\Big(2\alpha-4,\btheta^J_{16};x\Big),\\
		J_6 & = 2\sigma_1^{2\alpha} \mathcal{H}\Big(2(\alpha-1),\btheta^J_6;x\Big), & J_{17} & =2\sigma_1^{4\alpha} \mathcal{H}\Big(4(\alpha-1),\btheta^J_{17};x\Big),\\
		J_7 & = 2\sigma_1^{\alpha} \mathcal{H}\Big(\alpha-2,\btheta^J_7;x\Big), & J_{18} & = 2\sigma_1^{4\alpha} \mathcal{H}\Big(4(\alpha-1),\btheta^J_{18};x\Big),\\
		J_8 = J_9 = J_{12} & = 2\sigma_1^{3\alpha} \mathcal{H}\Big(3\alpha-4,\btheta^J_8;x\Big), & J_{19} & = 2\sigma_1^{4\alpha} \mathcal{H}\Big(4(\alpha-1),\btheta^J_{19};x\Big),\\
		J_{10} & = 2\sigma_1^{3\alpha} \mathcal{H}\Big(3\alpha-4,\btheta^J_{10};x\Big), &
		\end{align*}
		where $\btheta^J_i = (\theta^J_{i1},\theta^J_{i2})$, for $i=1,\ldots,19$, 
		\begin{align*}
		\theta^J_{11} & = a\Big(\lambda_2(\kappa_1^2-a^2\lambda_1^2)+2\kappa_1\kappa_2\lambda_1\Big), & \theta^J_{12} & = \kappa_2(\kappa_1^2 - a^2\lambda_1^2) - 2 a^2 \kappa_1 \lambda_1\lambda_2,\\
		\theta^J_{21} & = a\Big(K+\beta_1 L\Big), & \theta^J_{22} & = L - a^2\beta_1 K,\\
		\theta^J_{31} & = a\Big(\beta_1\kappa_4+\lambda_4\Big), & \theta^J_{32} & = \kappa_4 - a^2\beta_1\lambda_4,\\
		\theta^J_{41} & = a K, & \theta^J_{42} & = L,\\
		\theta^J_{61} & = L, & \theta^J_{62} & = - a K,\\
		\theta^J_{71} & = \kappa_4, & \theta^J_{72} & = - a\lambda_4,\\
		\theta^J_{81} & = L - a^2\beta_1K, & \theta^J_{82} & = -a\Big(K + \beta_1L\Big),\\
		\theta^J_{101} & = \kappa_4(1-a^2\beta_1^2) - 2a^2\beta_1\lambda_4, & \theta^J_{102} & = -a \Big(\lambda_4(1-a^2\beta_1^2) + 2\beta_1\kappa_4\Big),\\
		\theta^J_{111} & = \theta^J_{12}, & \theta^J_{112} & = -\theta^J_{11},\\
		\theta^J_{141} & = L, & \theta^J_{142} & = -a K,\\
		\theta^J_{151} & = \kappa_2^2 - a^2\lambda_2^2, & \theta^J_{152} & = -2a\kappa_2\lambda_2,\\
		\theta^J_{161} & = \kappa_4 - a^2\beta_1\lambda_4, & \theta^J_{162} & = -a\Big(\lambda_4+\beta_1\kappa_4\Big),\\
		\theta^J_{171} & = \theta^J_{12} - a\beta_1\theta^J_{11}, & \theta^J_{172} & = -\theta^J_{11} + a \theta^J_{12},\\
		\theta^J_{181} & = \kappa_1^4-6a^2\kappa_1^2\lambda_1^2+a^4\lambda_1^4, & \theta^J_{182} & = -4a\kappa_1\lambda_1(\kappa_1^2-a^2\lambda_1^2),\\
		\theta^J_{191} & = L(1-a^2\beta_1^2) - 2a^2\beta_1K, & \theta_{192} & = -a\Big(K(1-a^2\beta_1^2)+2\beta_1L\Big),
		\end{align*}
		and $K=\kappa_1\lambda_3+\lambda_1\kappa_3$, $L=\kappa_1\kappa_3-a^2\lambda_1\lambda_3$. Hence,
		\begin{align*}
		&\mathbb{E}\Big[X_2^4\Big|X_1=x\Big]  = \phi^{(4)}_{X_2|x}(0)\\
		&   \hspace{0.4cm} = \dfrac{-\alpha}{\pi f_{X_1}(x)} \Bigg[\alpha x \Big(\alpha K_1 + (\alpha-1)K_2\Big) + \alpha x^2 K_6 - (\alpha-1)x^2K_7 + \alpha^2(\alpha-1)K_3 + \alpha(\alpha-1)^2K_4 + \alpha^3 K_5\Bigg],
		\end{align*}
		where
		\begin{align*}
		K_1 & = \sigma_1^{3\alpha}\mathcal{H}\Big(3(\alpha-1),\btheta_{1}^K;x\Big), & \text{with} \quad \btheta_{1}^K & = 3\btheta^J_{1}-2\btheta^J_{2},\\
		K_2 & = \sigma_1^{2\alpha}\mathcal{H}\Big(2\alpha-3,\btheta_{2}^K;x\Big),& \text{with} \quad \btheta_{2}^K & = 2(\btheta^J_3-\btheta^J_4),\\
		K_3 & = \sigma_1^{3\alpha}\mathcal{H}\Big(3\alpha-4,\btheta_{3}^K;x\Big),& \text{with} \quad \btheta_{3}^K & = \btheta^J_{10} - 3\btheta^J_{11} - \btheta^J_{8},\\
		K_4 & = \sigma_1^{2\alpha}\mathcal{H}\Big(2\alpha-4,\btheta_{4}^K;x\Big),& \text{with} \quad \btheta_{4}^K & = 4 \btheta^J_{14} - 3 \btheta^J_{15}-\btheta^J_{16},\\
		K_5 & = \sigma_1^{4\alpha}\mathcal{H}\Big(4(\alpha-1),\btheta_{5}^K;x\Big),& \text{with} \quad \btheta_{5}^K & = 3\btheta^J_{17} - \btheta^J_{18} - \btheta^J_{19},\\
		K_6 & = \sigma_1^{2\alpha} \mathcal{H}\Big(2(\alpha-1),\btheta^K_6;x\Big), & \text{with} \quad \btheta_{6}^K & = \btheta^J_{6}, \\
		K_7 & = \sigma_1^{\alpha} \mathcal{H}\Big(\alpha-2,\btheta^K_7;x\Big), & \text{with} \quad \btheta_{7}^K & = \btheta^J_{7}.
		\end{align*}
		Invoking Lemmas \ref{le:SCFG} ($\iota\iota$) for $n=1,2,3$ and \ref{le:2a-4}, we get
		\begin{align*}
		\mathbb{E}\Big[X_2^4\Big|X_1=x\Big] & = \dfrac{-\alpha}{\pi f_{X_1}(x)} \Bigg[x^3\sigma_1^\alpha\Big(\theta_{72}^K C_1(x)-\theta_{71}^KS_1(x)\Big)\\
		& \hspace{0.5cm} + \dfrac{\alpha x^2 \sigma_1^{2\alpha}}{2}C_2(x)\bigg(-\theta_{22}^K+2\theta_{61}^K-2\Big(\theta_{71}^K+a\beta_1\theta_{72}^K\Big) - \dfrac{\alpha-1}{2\alpha-3}\theta_{41}^K\bigg) \\
		& \hspace{0.5cm} + \dfrac{\alpha x^2 \sigma_1^{2\alpha}}{2}S_2(x)\bigg(\theta_{21}^K+2\theta_{62}^K-2\Big(\theta_{72}^K-a\beta_1\theta_{71}^K\Big) - \dfrac{\alpha-1}{2\alpha-3}\theta_{42}^K\bigg)\Bigg] \\
		& \hspace{0.5cm} + \dfrac{\alpha^2 x \sigma_1^{3\alpha}}{6}C_3(x)\bigg(6\theta_{11}^K+3\Big(\theta_{21}^K+a\beta_1\theta_{22}^K\Big)-2\theta_{32}^K+5\dfrac{\alpha-1}{2\alpha-3}\Big(a\beta_1\theta_{41}^K-\theta_{42}^K\Big)\bigg)\\
		& \hspace{0.5cm} + \dfrac{\alpha^2 x \sigma_1^{3\alpha}}{6}S_3(x)\bigg(6\theta_{12}^K+3\Big(\theta_{22}^K-a\beta_1\theta_{21}^K\Big)+2\theta_{31}^K+5\dfrac{\alpha-1}{2\alpha-3}\Big(\theta_{41}^K+a\beta_1\theta_{42}^K\Big)\bigg)\\
		& \hspace{0.5cm} + \dfrac{\alpha^3 \sigma_1^{4\alpha}}{3}C_4(x)\bigg(\theta_{31}^K+a\beta_1\theta_{32}^K+\dfrac{\alpha-1}{2\alpha-3}\Big(\theta_{41}^K(1-a^2\beta_1^2)+2a\beta_1\theta_{42}^K\Big)+3\theta_{51}^K\bigg)\\
		& \hspace{0.5cm} + \dfrac{\alpha^3 \sigma_1^{4\alpha}}{3}S_4(x)\bigg(\theta_{32}^K-a\beta_1\theta_{31}^K+\dfrac{\alpha-1}{2\alpha-3}\Big(\theta_{42}^K(1-a^2\beta_1^2)-2a\beta_1\theta_{41}^K\Big)+3\theta_{52}^K\bigg)\Bigg].
		\end{align*}
		Using Lemma \ref{le:SCFG} ($\iota\iota\iota$) yields the conclusion. The coefficients $\btheta$'s in the expression \eqref{eq:fcm_gen} are deduced from the $\btheta^K$'s and $\btheta^J$'s as follows:
		
		\vspace*{-0.8cm}
		
		\begin{align}
		\theta_{41} & = -\theta_{22}^K+2\theta_{61}^K-2\Big(\theta_{71}^K+a\beta_1\theta_{72}^K\Big) - \dfrac{\alpha-1}{2\alpha-3}\theta_{41}^K,\label{def:nu_fourthcm1}\\
		\theta_{42} & = \theta_{21}^K+2\theta_{62}^K-2\Big(\theta_{72}^K-a\beta_1\theta_{71}^K\Big) - \dfrac{\alpha-1}{2\alpha-3}\theta_{42}^K,\\
		\theta_{51} & = 6\theta_{11}^K+3\Big(\theta_{21}^K+a\beta_1\theta_{22}^K\Big)-2\theta_{32}^K+5\dfrac{\alpha-1}{2\alpha-3}\Big(a\beta_1\theta_{41}^K-\theta_{42}^K\Big),\\
		\theta_{52} & = 6\theta_{12}^K+3\Big(\theta_{22}^K-a\beta_1\theta_{21}^K\Big)+2\theta_{31}^K+5\dfrac{\alpha-1}{2\alpha-3}\Big(\theta_{41}^K+a\beta_1\theta_{42}^K\Big),\\
		\theta_{61} & = \theta_{31}^K+a\beta_1\theta_{32}^K+\dfrac{\alpha-1}{2\alpha-3}\Big(\theta_{41}^K(1-a^2\beta_1^2)+2a\beta_1\theta_{42}^K\Big)+3\theta_{51}^K,\\
		\theta_{62} & = \theta_{32}^K-a\beta_1\theta_{31}^K+\dfrac{\alpha-1}{2\alpha-3}\Big(\theta_{42}^K(1-a^2\beta_1^2)-2a\beta_1\theta_{41}^K\Big)+3\theta_{52}^K.\label{def:nu_fourthcm2}
		\end{align}
		
		\subsection{Lemmas for the proof of Theorem \ref{theo:cond_moments_stable_general_alpnot1}}
		
		\noindent The following elementary Lemmas, stated without proof, are used to establish Theorem \ref{theo:cond_moments_stable_general_alpnot1}.
		
		\begin{lem}\label{le:SCFG}
			Let $\alpha\in(1,2)$, $b>0$, $c\in\mathbb{R}$. Define for $n\ge1$ and $x\in\mathbb{R}$
			\begin{align*}
			C_n(x) & = \int_{0}^{+\infty} e^{-bt^\alpha}t^{n(\alpha-1)}\cos(tx-ct^\alpha)dt, & F_n(x) & = \int_{0}^{+\infty} e^{-bt^\alpha}t^{n(\alpha-1)-1}\cos(tx-ct^\alpha)dt,\\
			S_n(x) & = \int_{0}^{+\infty} e^{-bt^\alpha}t^{n(\alpha-1)}\sin(tx-ct^\alpha)dt, & G_n(x) & = \int_{0}^{+\infty} e^{-bt^\alpha}t^{n(\alpha-1)-1}\sin(tx-ct^\alpha)dt.
			\end{align*}
			$\iota)$ Then the following hold for any $n\ge1$ and $x\in\mathbb{R}$
			\begin{align*}
			F_{n}(x) & = \dfrac{\alpha\Big(bC_{n+1}(x)-cS_{n+1}(x)\Big)+xS_{n}(x)}{n(\alpha-1)}, &
			G_{n}(x) & = \dfrac{\alpha\Big(cC_{n+1}(x)+bS_{n+1}(x)\Big)-xC_{n}(x)}{n(\alpha-1)}.
			\end{align*}
			$\iota\iota)$ For any $n\ge1$, $\theta_1,\theta_2\in\mathbb{R}$ and $x\in\mathbb{R}$:
			\begin{align*}
			\theta_1F_n(x)+\theta_2G_n(x) & =\dfrac{\alpha\Big[C_{n+1}(x) \Big(b\theta_1+c\theta_2\Big)+S_{n+1}(x) \Big(b\theta_2-c\theta_1\Big)\Big]+x\Big[-\theta_2C_{n}(x) + \theta_1S_{n}(x) \Big]}{n(\alpha-1)}.
			\end{align*}
			$\iota\iota\iota)$ We have for $x\in\mathbb{R}$, $b=\sigma_1^\alpha$ and $c=a\beta_1\sigma_1^\alpha$:
			\begin{align*}
			C_1(x) & = \dfrac{a\beta_1x\pi f_{X_1}(x)+1-xH(x)}{\alpha\sigma_1^\alpha(1+(a\beta_1)^2)}, &
			S_1(x) & = \dfrac{x\pi f_{X_1}(x)-a\beta_1(1-xH(x))}{\alpha\sigma_1^\alpha(1+(a\beta_1)^2)}.
			\end{align*}
		\end{lem}

		\begin{lem}\label{le:2a-4}
			Let $\alpha\in(3/2,2)$, $b>0$, $c\in\mathbb{R}$. Define for $x\in\mathbb{R}$
			\begin{align*}
			h_c(x) & = \int_{0}^{+\infty} e^{-bt^\alpha}t^{2\alpha-4}\cos(tx-ct^\alpha)dt, & h_s(x) & = \int_{0}^{+\infty} e^{-bt^\alpha}t^{2\alpha-4}\sin(tx-ct^\alpha)dt.
			\end{align*}
			Then for any $\theta_1,\theta_2\in\mathbb{R}$ and $x\in\mathbb{R}$,
			\begin{align*}
			\theta_1h_c(x)+\theta_2h_s(x) & = \dfrac{\alpha^2}{3(2\alpha-3)(\alpha-1)}\bigg[C_{4}(x) \Big(\theta_1(b^2-c^2)+2bc\theta_2\Big)+S_{4}(x) \Big(\theta_2(b^2-c^2)-2bc\theta_1\Big)\bigg]\\
			& \hspace{0.5cm} + \dfrac{5\alpha x}{6(2\alpha-3)(\alpha-1)}\bigg[C_{3}(x) \Big(c\theta_1-b\theta_2\Big) + S_{3}(x)\Big(b\theta_1+c\theta_2\Big) \bigg]\\
			& \hspace{0.5cm} - \dfrac{x^2}{2(2\alpha-3)(\alpha-1)}\bigg[\theta_1C_2(x)+\theta_2S_2(x)\bigg].
			\end{align*}
		\end{lem}
		
		\label{page:theo22_pend}

		

		\section{Proof of Theorem \ref{theo:cond_moments_stable_general_alp1_scm}}
		\label{sec:tha4}
		
		\label{page:theo24_pini}
		
		Let $\bX=(X_1,X_2)$ be an $\alpha$-stable vector with $\alpha=1$ and spectral representation $(\Gamma,\boldsymbol{0})$. Its characteristic function, denoted $\varphi_{\bX} (t,r)$ for any $(t,r)\in\mathbb{R}^2$, reads
		\begin{equation}\label{def:char_proof_alp1}
		\varphi_{\bX} (t,r) = \exp\left\{-\int_{S_2} |ts_1+rs_2| + i a(ts_1+rs_2)\ln|ts_1+rs_2|\Gamma(ds)\right\},    
		\end{equation}
		with $a=2/\pi$. The conditional characteristic function of $X_2$ given $X_1=x$, denoted $\phi_{X_2|x}(r)$ for $r\in\mathbb{R}$, is still given by \eqref{def:cond_char_proof}.
		
		\begin{lem}
			\label{lem:alp1}
			Let $(X_1,X_2)$ be an $\alpha$-stable random vector with $\alpha=1$ and spectral representation $(\Gamma,\boldsymbol{0})$. If \eqref{eq:nu_cond} holds with $\nu>0$, the first derivative of $\phi_{X_2|x}$ is given by
			\begin{align*}
			\phi^{(1)}_{X_2|x}(r) & = \dfrac{-1}{2\pi f_{X_1}(x)} \Big(A_1 + iaA_2\Big),
			\end{align*}
			with
			\begin{align}
			A_1 & = \int_{\mathbb{R}}e^{-itx} \varphi_{\bX}(t,r)\bigg(\int_{S_2}s_2(ts_1+rs_2)^{<0>}\Gamma(d\bs)\bigg)dt,\label{dem:j1}\\
			A_2 & = \int_{\mathbb{R}}e^{-itx} \varphi_{\bX}(t,r)\bigg(\int_{S_2}s_2(1+\ln|ts_1+rs_2|)\Gamma(d\bs)\bigg)dt\label{dem:j2}
			\end{align}
			If \eqref{eq:nu_cond} holds with $\nu>1$, the second derivative of $\phi_{X_2|x}$ is given by
			\begin{align}\label{eq:d2_alp1}
			\phi^{(2)}_{X_2|x}(r) & = \dfrac{-1}{2\pi f_{X_1}(x)} \Big(-B_1+ixB_2+B_3\Big),
			\end{align}
			where,
			\begin{align*}
			B_1 & = \int_{\mathbb{R}}e^{-itx}\varphi_{\bX}(t,r)\bigg(\int_{S_2}s_2(ts_1+rs_2)^{<0>}+ias_2(1+\ln|ts_1+rs_2|\Gamma(d\bs)\bigg)^2dt,\\
			B_2 & = \int_{\mathbb{R}}e^{-itx}\varphi_{\bX}(t,r)\bigg(\int_{S_2}\Big((ts_1+rs_2)^{<0>}+ia(1+\ln|ts_1+rs_2|\Big)s_2^2s_1^{-1}\Gamma(d\bs)\bigg)dt,\\
			B_3 & = \int_{\mathbb{R}}e^{-itx}\varphi_{\bX}(t,r)\bigg(\int_{S_2}s_1(ts_1+rs_2)^{<0>}+ias_1(1+\ln|ts_1+rs_2|\Gamma(d\bs)\bigg)\\
			& \hspace{5cm} \times \bigg(\int_{S_2}\Big((ts_1+rs_2)^{<0>}+ia(1+\ln|ts_1+rs_2|\Big)s_2^2s_1^{-1}\Gamma(d\bs)\bigg)dt.
			\end{align*}
		\end{lem}

		\subsection{Justifying inversion of integral and derivative signs}
		
		\noindent \textbf{First derivative}
		
		The terms depending on $r$ in the right-hand side of \eqref{def:char_proof_alp1} are of the form (omitting the factor $1/2\pi f_{X_1}(x)$)
		\begin{align*}
		\int_{\mathbb{R}}e^{-\int_{S_2} |ts_1+rs_2|\Gamma(d\bs)} \text{trig} \bigg(-tx - a\int_{S_2}(ts_1+rs_2)\ln|ts_1+rs_2|\Gamma(d\bs)\bigg)dt.
		\end{align*}
		Consider for instance the term obtained by replacing trig by the cosine function, denoted $I_1$.
		\begin{align*}
		I_1'(r) & = \lim_{h\rightarrow0} \dfrac{1}{h} \int_{\mathbb{R}} \Bigg[e^{-\int_{S_2}|ts_1+(r+h)s_2|\Gamma(d\bs)}-e^{-\int_{S_2}|ts_1+rs_2|\Gamma(d\bs)}\Bigg]\\
		& \hspace{5cm} \times \cos\bigg(tx+ a\int_{S_2}(ts_1+(r+h)s_2)\ln|ts_1+(r+h)s_2|\Gamma(d\bs)\bigg)dt\\
		& + \lim_{h\rightarrow0} \dfrac{1}{h} \int_{\mathbb{R}} e^{-\int_{S_2}|ts_1+rs_2|\Gamma(d\bs)} \Bigg[\cos\bigg(tx+ a\int_{S_2}(ts_1+(r+h)s_2)\ln|ts_1+(r+h)s_2|\Gamma(d\bs)\bigg)\\
		& \hspace{7cm} - \cos\bigg(tx+ a\int_{S_2}(ts_1+rs_2)\ln|ts_1+rs_2|\Gamma(d\bs)\bigg)\Bigg]dt\\
		& := I_{11} + I_{12}
		\end{align*}
		The integrand of $I_{11}$ converges to 
		\begin{align*}
		-e^{-\int_{S_2}|ts_1+rs_2|\Gamma(d\bs)}\cos\bigg(tx+ a\int_{S_2}(ts_1+rs_2)\ln|ts_1+rs_2|\Gamma(d\bs)\bigg)\int_{S_2}s_2(ts_1+rs_2)^{<0>}\Gamma(d\bs).
		\end{align*}
		Using \eqref{eq:exp94} we can bound the integrand of $I_{11}$ by
		\begin{align*}
		\dfrac{1}{|h|} \Bigg|\int_{S_2} |ts_1+(r+h)s_2|-|ts_1+rs_2|\Gamma(d\bs)\Bigg|e^{-\int_{S_2}|ts_1+rs_2|\Gamma(d\bs)} e^{\big|\int_{S_2} |ts_1+(r+h)s_2|-|ts_1+rs_2|\Gamma(d\bs)\big|}.
		\end{align*}
		By Lemma \ref{le:34} $(\iota)$ and the triangle inequality, we can further bound it for $|h|<|r|$ by
		\begin{align*}
		\sigma_2e^{\sigma_2(1+|r|)-\sigma_1|t|},
		\end{align*}
		which does not depend on $h$ and is integrable with respect to $t$ on $\mathbb{R}$. The dominated convergence theorem applies to $I_{11}$. Turning to $I_{12}$, its integrand converges to
		\begin{align*}
		-ae^{-\int_{S_2}|ts_1+rs_2|\Gamma(d\bs)}\sin\bigg(tx+ a\int_{S_2}(ts_1+rs_2)\ln|ts_1+rs_2|\Gamma(d\bs)\bigg)\int_{S_2}s_2(1 + \ln|ts_1+rs_2|)\Gamma(d\bs).
		\end{align*}
		Using the mean value theorem on the cosine, its integrand can be bounded by
		\begin{align}
		& \dfrac{a}{|h|}e^{-\int_{S_2}|ts_1+rs_2|\Gamma(d\bs)} \Bigg|\int_{S_2}(ts_1+(r+h)s_2)\ln|ts_1+(r+h)s_2|-(ts_1+rs_2)\ln|ts_1+rs_2|\Gamma(d\bs)\Bigg|\nonumber\\
		& \le ae^{\sigma_2|r|-\sigma_1|t|}\dfrac{1}{|h|}\int_{S_2}\bigg|(ts_1+(r+h)s_2)\ln|ts_1+(r+h)s_2|-(ts_1+rs_2)\ln|ts_1+rs_2|\bigg|\Gamma(d\bs)\nonumber\\
		& := ae^{\sigma_2|r|-\sigma_1|t|}\Big(Q_1+Q_2\Big),
		\label{eq:alp1_i12}
		\end{align}
		where the two terms $Q_1$ and $Q_2$ involve integrals over $S_2\cap\{\bs:|ts_1+rs_2|\ge2|h|\}$ and $S_2\cap\{\bs:|ts_1+rs_2|<2|h|\}$. Focus on $Q_2$. Introduce the function $f: \mathbb{R}_+ \rightarrow\mathbb{R}_+$ defined for any $z\ge0$ by $f(z) = z|\ln z|$. It is such that $f(0)=0$ and for $z$ small enough ($0<z<e^{-1}$), $f$ is monotone increasing. Since $|ts_1+rs_2|<2|h|$, we also have $|ts_1+(r+h)s_2|<3|h|$. Thus, for $0<|h|<(3e)^{-1}$, the integrand of $Q_2$ can be bounded by
		\begin{align*}
		|h|^{-1} \bigg(\Big|f(|3h|)\Big| + \Big|f(|2h|)\Big|\bigg) \le 2|h|^{-1}\Big|f(|3h|)\Big| \le 6\Big|ln|3h|\Big|
		\end{align*}
		Using Lemma \ref{le:li}, we can bound the later quantity for any $v>0$ by
		\begin{align*}
		6 v^{-1} \Big( 2 + |3h|^v + |3h|^{-v}\Big).
		\end{align*}
		From $|ts_1+rs_2|/2<|h|<(3e)^{-1}$, we deduce that $|3h|^{-v}<\Big(3|ts_1+rs_2|/2\Big)^{-v}$ and
		\begin{align*}
		6 v^{-1} \Big( 2 + |3h|^v + |3h|^{-v}\Big) \le 6 v^{-1} \Big( 2 + e^{-v} + \Big(3|ts_1+rs_2|/2\Big)^{-v}\Big) \le \text{const}_1 + \text{const}_2 |ts_1+rs_2|^{-v},
		\end{align*}
		for some nonnegative constants $\text{const}_1$ and $\text{const}_2$. Hence, the term involving $Q_2$ in \ref{eq:alp1_i12} can be further bounded for any $v>0$ by 
		\begin{align}
		ae^{\sigma_2|r|-\sigma_1|t|}\Big(\text{const}_1 + \text{const}_2 \int_{S_2}\Big|t+\dfrac{rs_2}{s_1}\Big|^{-v}|s_1|^{-v}\Gamma(d\bs)\Big).\label{dem:q2} 
		\end{align}
		The term with $\text{const}_1$ is clearly integrable with respect to $t$ on $\mathbb{R}$. Letting \eqref{eq:nu_cond} hold with $\nu>0$, choose some $v\in(0,\min(\nu,1))$. We show that the second term is bounded by an integrable function of $t$ as we did in Equation \eqref{dem:bound_lem22} using Lemma \ref{le:310} with $\eta=v$, $b=0$, $p=0$, the fact that $\int_{\mathbb{R}}e^{-\sigma_1 |t|} |t|^{-v}dt<+\infty$ and \eqref{eq:nu_cond} with $\nu>v>0$. There remains to be bounded the part involving $Q_1$ in \eqref{eq:alp1_i12}. For this term, we apply the mean value theorem to the function $z\longmapsto z\ln|z|$ and get that
		\begin{align*}
		& |h|^{-1}\bigg|(ts_1+(r+h)s_2)\ln|ts_1+(r+h)s_2|-(ts_1+rs_2)\ln|ts_1+rs_2|\bigg| \\
		& \hspace{2cm} \le |h|^{-1} |hs_2| \Big|1 + \ln|u| \Big| \\
		& \hspace{2cm} \le 1 + \Big|\ln|u|\Big|,
		\end{align*}
		for some $u\in[ts_1+(r+h)s_2\wedge ts_1+rs_2,ts_1+(r+h)s_2\vee ts_1+rs_2]$. Since $Q_1$ is an integral over $S_2\cap\{\bs:|ts_1+rs_2|\ge2|h|\}$, we have $|u|\in\Big[\frac{|ts_1+rs_2|}{2},2|ts_1+rs_2|\Big]$, and because of the quasi-convexity of the function $z\longmapsto\Big|\ln|z|\Big|$, we can bound the above term by
		\begin{align*}
		1 + \Bigg|\ln\bigg|\dfrac{ts_1+rs_2}{2}\bigg|\Bigg| +\Bigg|\ln|2(ts_1+rs_2)|\Bigg| \le \text{const} + 2 \Big|\ln|ts_1+rs_2|\Big|.
		\end{align*}
		Using Lemma \ref{le:li}, we can bound this term for any $v>0$ by
		\begin{align*}
		\text{const} + 2 v^{-1}\Big( 2 + |ts_1+rs_2|^v + |ts_1+rs_2|^{-v} \Big) \le \text{const}_1 + \text{const}_2 |t|^v + \text{const}_3 \Big|t+\dfrac{rs_2}{s_1}\Big|^{-v} |s_1|^{-v}
		\end{align*}
		Hence, the term in \eqref{eq:alp1_i12} involving $Q_1$ can be bounded for any $v>0$ by
		\begin{align}
		ae^{\sigma_2|r|-\sigma_1|t|}\Big(\text{const}_1 + \text{const}_2 |t|^v + \text{const}_3 \int_{S_2}\Big|t+\dfrac{rs_2}{s_1}\Big|^{-v}|s_1|^{-v}\Gamma(d\bs)\Big).   \label{dem:q1} 
		\end{align}
		which can be shown to be integrable with respect to $t$ on $\mathbb{R}$ as we did above for the term with $Q_2$. The dominated convergence theorem applies to $I_{12}$ and thus to $I_1$. We can derivate $\phi_{X_2|x}$ under the integral sign.

		\noindent \textbf{Second derivative}\\
		Let us start with $A_2$, which is the most delicate. It is composed of terms of the form
		\begin{align*}
		& \int_{\mathbb{R}} e^{-\int_{S_2}|ts_1+rs_2|\Gamma(d\bs)}\text{trig}\bigg(-tx-a\int_{S_2}(ts_1+rs_2)\ln|ts_1+rs_2|\Gamma(d\bs)\bigg)\\
		& \hspace{8cm} \times \bigg(\int_{S_2}s_2(1+\ln|ts_1+rs_2|)\Gamma(d\bs)\bigg)dt,
		\end{align*}
		where <<trig>> stands for sine or cosine. Denoting the one with cosine as $K_2$, we have
		\begin{align*}
		K_2 & = \lim_{h\rightarrow0}\dfrac{1}{h}\int_{\mathbb{R}}\Bigg[e^{-\int_{S_2}|ts_1+(r+h)s_2|\Gamma(d\bs)}-e^{-\int_{S_2}|ts_1+rs_2|\Gamma(d\bs)}\Bigg]\\
		&  \hspace{3cm} \times \cos\bigg(tx+a\int_{S_2}(ts_1+(r+h)s_2)\ln|ts_1+(r+h)s_2|\Gamma(d\bs)\bigg)\\
		& \hspace{8cm} \times \bigg(\int_{S_2}s_2(1+\ln|ts_1+(r+h)s_2|)\Gamma(d\bs)\bigg)dt\\
		%
		& + \lim_{h\rightarrow0}\dfrac{1}{h}\int_{\mathbb{R}}e^{-\int_{S_2}|ts_1+rs_2|\Gamma(d\bs)}\Bigg[\cos\bigg(tx+a\int_{S_2}(ts_1+(r+h)s_2)\ln|ts_1+(r+h)s_2|\Gamma(d\bs)\bigg)\\
		& \hspace{6cm} - \cos\bigg(tx+a\int_{S_2}(ts_1+rs_2)\ln|ts_1+rs_2|\Gamma(d\bs)\bigg)\Bigg]\\
		& \hspace{8cm}\times\bigg(\int_{S_2}s_2(1+\ln|ts_1+(r+h)s_2|)\Gamma(d\bs)\bigg)dt\\
		%
		& + \lim_{h\rightarrow0}\dfrac{1}{h}\int_{\mathbb{R}}e^{-\int_{S_2}|ts_1+rs_2|\Gamma(d\bs)}\cos\bigg(tx+a\int_{S_2}(ts_1+rs_2)\ln|ts_1+rs_2|\Gamma(d\bs)\bigg)\\
		& \hspace{4cm} \times\Bigg[\int_{S_2}s_2\ln|ts_1+(r+h)s_2|-s_2\ln|ts_1+rs_2|\Gamma(d\bs)\Bigg]dt\\
		& := K_{21} + K_{22} + K_{23}.
		\end{align*}
		The integrand of $K_{21}$ converges to 
		\begin{align*}
		& -e^{-\int_{S_2}|ts_1+rs_2|\Gamma(d\bs)} \cos\bigg(tx+a\int_{S_2}(ts_1+rs_2)\ln|ts_1+rs_2|\Gamma(d\bs)\bigg)\\
		& \hspace{4cm}\times \bigg(\int_{S_2}s_2(ts_1+rs_2)^{<0>}\Gamma(d\bs)\bigg)\bigg(\int_{S_2}s_2(1+\ln|ts_1+rs_2|)\Gamma(d\bs)\bigg).
		\end{align*}
		Using \eqref{eq:exp94}, the triangle inequality and \eqref{le:33}, it can be bounded by
		\begin{align}
		\sigma_2e^{\sigma_2(1+|r|) - \sigma_1|t|} \int_{S_2}|s_2|\Big|1+\ln|ts_1+(r+h)s_2|\Big|\Gamma(d\bs).\label{dem:k11}
		\end{align}
		The integrand of the above expression can be bounded using Lemma \ref{le:li} for any $v>0$ by
		\begin{align*}
		& 1 + v^{-1}\Big(2 + |ts_1+(r+h)s_2|^v + |ts_1+(r+h)s_2|^{-v}\Big) \\
		& \hspace{2cm} \le \text{const}_1 + \text{const}_2 |t|^v + \text{const}_3 \Big|t+\dfrac{(r+h)s_2}{s_1}\Big|^{-v} |s_1|^{-v},
		\end{align*}
		hence, \eqref{dem:k11} is bounded by
		\begin{align*}
		\sigma_2e^{\sigma_2(1+|r|) - \sigma_1|t|} \Big(\text{const}_1 + \text{const}_2 |t|^v + \text{const}_3\int_{S_2} \Big|t+\dfrac{(r+h)s_2}{s_1}\Big|^{-v} |s_1|^{-v}\Gamma(d\bs)\Big).
		\end{align*}
		The terms involving $\text{const}_1$ and $\text{const}_2$ are clearly integrable with respect to $t$. The last term is more intricate as it still depends on $h$. We will show that the generalised Lebesgue dominated convergence theorem (Theorem 19, p.89 in \cite{roy10}) applies. Denoting
		\begin{align*}
		T(h) & = e^{- \sigma_1|t|}\Big|t+\dfrac{(r+h)s_2}{s_1}\Big|^{-v} |s_1|^{-v},
		\end{align*}
		it can be shown that $T(0)$ is integrable with respect to $t$ on $\mathbb{R}$ and $\Gamma$ on $S_2$ invoking the usual arguments. Also, choosing some $v\in(0,1)$, with have by Lemma \ref{le:312} with $\eta=-v$, $b=0$ and $0<p<1-v$,
		\begin{align*}
		\Bigg|\int T(h) - T(0) \Bigg| & \le \int_{S_2} |s_1|^{-v} \int_{\mathbb{R}} e^{- \sigma_1|t|}\Bigg|\Big|t+\dfrac{(r+h)s_2}{s_1}\Big|^{-v}-\Big|t+\dfrac{rs_2}{s_1}\Big|^{-v}\Bigg| dt \Gamma(d\bs)\\
		& \le \text{const }\int_{S_2} |s_1|^{-v} \Big|\dfrac{hs_2}{s_1}\Big|^p \Gamma(d\bs)\\
		& \le \text{const }|h|^p \int_{S_2} |s_1|^{-v-p}\Gamma(d\bs) \underset{h\rightarrow0}{\longrightarrow} 0,
		\end{align*}
		because \eqref{eq:nu_cond} holds with $\nu>1$ and $v+p<v+1-v<1$. Since $T(0)$ is integrable and $\lim_{h\rightarrow0}\int T(h)=\int T(0)$, the generalised dominated convergence theorem applies to $K_{21}$. We turn to $K_{22}$. Its integrand converges to 
		\begin{align*}
		& -ae^{-\int_{S_2}|ts_1+rs_2|\Gamma(d\bs)} \sin\bigg(tx+a\int_{S_2}(ts_1+rs_2)\ln|ts_1+rs_2|\Gamma(d\bs)\bigg)\\
		& \hspace{8cm}\times\bigg(\int_{S_2}s_2(1+\ln|ts_1+rs_2|)\Gamma(d\bs)\bigg)^2.
		\end{align*}
		With the usual inequalities and Lemma \ref{le:li}, it can be bounded for any $v>0$ by
		\begin{align*}
		& \dfrac{a}{|h|}e^{\sigma_2|r|-\sigma_1|t|} \Bigg|\int_{S_2}(ts_1+(r+h)s_2)\ln|ts_1+(r+h)s_2|-(ts_1+rs_2)\ln|ts_1+rs_2|\Gamma(d\bs)\Bigg|\\
		& \hspace{9cm}\times\Bigg|\int_{S_2}s_2(1+\ln|ts_1+(r+h)s_2|)\Gamma(d\bs)\Bigg|\\
		& \le ae^{\sigma_2|r|-\sigma_1|t|}\Big(Q_1+Q_2\Big)\Big(\sigma_2 + \int_{S_2}\Big|\ln|ts_1+(r+h)s_2|\Big|\Gamma(d\bs) \Big)\\
		& \le ae^{\sigma_2|r|-\sigma_1|t|}\Big(Q_1+Q_2\Big)\Big(\text{const}_1 + \text{const}_2|t|^v + \text{const}_3 \int_{S_2}\Big|t+\dfrac{(r+h)s_2}{s_1}\Big|^{-v} |s_1|^{-v}\Gamma(d\bs) \Big),
		\end{align*}
		where, similarly to \eqref{eq:alp1_i12}, the two terms $Q_1$ and $Q_2$ involve integrals over $S_2\cap\{\bs:|ts_1+rs_2|\ge2|h|\}$ and $S_2\cap\{\bs:|ts_1+rs_2|<2|h|\}$. After expansion, the terms with $\text{const}_1$ and $\text{const}_2$ are readily dealt with by following the method developed for \eqref{eq:alp1_i12}. Focus on the remaining term 
		$$
		a\int_{S_2}e^{\sigma_2|r|-\sigma_1|t|}(Q_1+Q_2)\Big|t+\dfrac{(r+h)s_2}{s_1}\Big||s_1|^{-v}\Gamma(d\bs).
		$$
		In view of the bounds \eqref{dem:q2} and \eqref{dem:q1}, the integrand can be bounded (up to a multiplicative constant) by 
		\begin{align*}
		U(h) = e^{-\sigma_1|t|} \Big|t+\dfrac{rs_2}{s_1}\Big|^{-v}  \Big|t+\dfrac{(r+h)s_2'}{s_1'}\Big|^{-v} |s_1|^{-v}|s_1'|^{-v}.
		\end{align*}
		Choosing some $v\in(0,1/2)$, we can invoke Lemma \eqref{le:cor31} with $\eta=-v$, $p=0$ and the fact that $\int_{\mathbb{R}}e^{-\sigma_1|t|}|t|^{-2v}dt<+\infty$ to show that $U(0)$ is integrable on the one hand. On the other hand we can again invoke Lemma \eqref{le:cor31}, this time with $\eta=-v$, $0<p<1-2v$, and the fact that \eqref{eq:nu_cond} holds with $\nu>1>v+1-2v>v+p$ to show that $\int U(h) \rightarrow \int U(0)$. The generalised dominated convergence theorem applies to $K_{12}$.\\
		We turn to $K_{23}$ for which <<appropriate integration by parts>> is required. After obvious manipulations,
		\begin{align*}
		%
		K_{23} & = \lim_{h\rightarrow0}\dfrac{1}{h} \int_{\mathbb{R}} \int_{S_2} s_2'\ln|ts_1'+rs_2'|\Bigg[e^{-\int_{S_2}\Big|\Big(t-\frac{hs_2'}{s_1'}\Big)s_1+rs_2\Big|\Gamma(d\bs)}-e^{-\int_{S_2}|ts_1+rs_2|\Gamma(d\bs)}\Bigg]\\
		& \hspace{1cm} \times \cos\Bigg(\Big(t-\dfrac{hs_2'}{s_1'}\Big)x+a\int_{S_2}\bigg(\Big(t-\dfrac{hs_2'}{s_1'}\Big)s_1+rs_2\bigg)\ln\bigg|\Big(t-\dfrac{hs_2'}{s_1'}\Big)s_1+rs_2\bigg|\Gamma(d\bs)\Bigg)\Gamma(d\bs')\\
		%
		& + \lim_{h\rightarrow0}\dfrac{1}{h}\int_{\mathbb{R}}\int_{S_2} s_2'\ln|ts_1'+rs_2'|e^{-\int_{S_2}|ts_1+rs_2|\Gamma(d\bs)}\\
		& \hspace{2cm}\times \Bigg[\cos\Bigg(\Big(t-\dfrac{hs_2'}{s_1'}\Big)x+a\int_{S_2}\bigg(\Big(t-\dfrac{hs_2'}{s_1'}\Big)s_1+rs_2\bigg)\ln\bigg|\Big(t-\dfrac{hs_2'}{s_1'}\Big)s_1+rs_2\bigg|\Gamma(d\bs)\Bigg)\\
		& \hspace{4cm} - \cos\Bigg(tx+a\int_{S_2}(ts_1+rs_2)\ln|ts_1+rs_2|\Gamma(d\bs)\Bigg)\Bigg]\Gamma(d\bs')\\
		& := L_1 + L_2.
		\end{align*}
		Starting with $L_1$, its integrand converges to
		\begin{align*}
		& e^{-\int_{S_2}|ts_1+rs_2|\Gamma(d\bs)} \cos\bigg(tx+a\int_{S_2}(ts_1+rs_2)\ln|ts_1+rs_2|\Gamma(d\bs)\bigg)\\
		& \hspace{5cm} \times \bigg(\int_{S_2}s_1(ts_1+rs_2)^{<0>} \Gamma(d\bs)\bigg) \bigg(\int_{S_2}\ln|ts_1+rs_2| {s_2}^2{s_1}^{-1}\Gamma(d\bs)\bigg)
		\end{align*}
		It can be bounded using \eqref{eq:exp95} and Lemma \ref{le:34} $(\iota)$ by
		\begin{align*}
		& \bigg|\dfrac{s_2'\ln|ts_1'+rs_2'|}{h}\bigg|\exp\Bigg\{-\min\bigg(\int_{S_2}\Big|\Big(t-\dfrac{hs_2'}{s_1'}\Big)s_1+rs_2\Big|\Gamma(d\bs),\int_{S_2}|ts_1+rs_2|\Gamma(d\bs)\bigg)\Bigg\}\\
		& \hspace{5cm} \times \Bigg|\int_{S_2} \bigg|\Big(t-\dfrac{hs_2'}{s_1'}\Big)s_1+rs_2\bigg| - |ts_1+rs_2|\Gamma(d\bs)\Bigg|\\
		& \le  e^{\sigma_2|r|} \exp\Bigg\{-\sigma_1\min\bigg(\Big|t-\dfrac{hs_2'}{s_1'}\Big|,|t|\bigg)\Bigg\}\Big|s_2'\ln|ts_1'+rs_2'|\Big| \dfrac{1}{|h|} \int_{S_2}\Big|\dfrac{hs_2'}{s_1'}s_1\Big|\Gamma(d\bs)\\
		& \le \sigma_1 e^{\sigma_2|r|} \exp\Bigg\{-\sigma_1\min\bigg(\Big|t-\dfrac{hs_2'}{s_1'}\Big|,|t|\bigg)\Bigg\} \Big|\ln|ts_1'+rs_2'|\Big| |s_2'|^2 |s_1'|^{-1}\\
		& := V(h).
		\end{align*}
		We follow a similar procedure as the one used in \cite{ct98} (p.51) to deal with the min inside the exponential. Focus on the case $\dfrac{hs_2}{s_1}>0$ (the converse case is similar). We have
		\begin{align*}
		\min\bigg(\Big|t-\dfrac{hs_2'}{s_1'}\Big|,|t|\bigg) = \left\{
		\begin{array}{ll}
		\Big|t-\dfrac{hs_2'}{s_1'}\Big|, \hspace{1cm} \text{if} \quad t \ge hs_2'/2s_1',\\
		\hspace{0.5cm} |t|, \hspace{1.5cm} \text{if} \quad t < hs_2'/2s_1'.  \\
		\end{array}
		\right.
		\end{align*}
		Thus, up to a multiplicative constant,
		\begin{align*}
		\int_{\mathbb{R}} V(h) dt & = \int_{\frac{hs_2}{2s_1}}^{+\infty}e^{-\sigma_1|t-\frac{hs_2}{s_1}|} \Big|\ln|ts_1+rs_2|\Big| |s_2|^2 |s_1|^{-1} dt + \int_{-\infty}^{-\frac{hs_2}{2s_1}} e^{-\sigma_1|t|} \Big|\ln|ts_1+rs_2|\Big| |s_2|^2 |s_1|^{-1} dt\\
		& = \int_{-\frac{hs_2}{2s_1}}^{+\infty} e^{-\sigma_1|t|}\bigg|\ln\Big|ts_1+rs_2+\dfrac{hs_2}{s_1}\Big|\bigg||s_2|^2 |s_1|^{-1}dt + \int_{-\infty}^{-\frac{hs_2}{2s_1}} e^{-\sigma_1|t|} \Big|\ln|ts_1+rs_2|\Big| |s_2|^2 |s_1|^{-1} dt\\
		& = \int_{\mathbb{R}} e^{-\sigma_1|t|}\bigg[\Big|\ln|ts_1+(r+h)s_2|\Big|\mathds{1}_{\{t\ge -hs_2/2s_1\}} + \Big|\ln|ts_1+rs_2|\Big|\mathds{1}_{\{t \le -hs_2/2s_1\}}\bigg]|s_2|^2 |s_1|^{-1} dt.
		\end{align*}
		Thus, using Lemma \ref{le:li}, we can bound the integrand for any $v>0$ and $|h|<|r|$ by
		\begin{align*}
		& e^{-\sigma_1|t|}\bigg[\Big|\ln|ts_1+(r+h)s_2|\Big| + \Big|\ln|ts_1+rs_2|\Big|\bigg]|s_2|^2 |s_1|^{-1}\\
		& \hspace{1cm} \le v^{-1}e^{-\sigma_1|t|}\bigg[\text{const}_1 + \text{const}_2 |t|^v \\
		& \hspace{4cm} + \text{const}_3 \Big|t +\dfrac{rs_2}{s_1}\Big|^{-v} |s_1|^{-v} + \text{const}_4 \Big|t +\dfrac{(r+h)s_2}{s_1}\Big|^{-v}|s_1|^{-v} \bigg]|s_2|^2 |s_1|^{-1}.
		\end{align*}
		Clearly, the terms involving $\text{const}_1$ and $\text{const}_2$ are integrable with respect to $t$ and $\Gamma$. Denoting the last term as $V_4(h) := e^{-\sigma_1|t|} \Big|t +\dfrac{(r+h)s_2}{s_1}\Big|^{-v} |s_2|^2 |s_1|^{-1-v}$, we show that the generalised dominated convergence theorem applies. As \eqref{eq:nu_cond} holds for some $\nu>1$, choose $v = \dfrac{\nu-1}{2}>0$ if $\nu<2$, and some $v\in(0,1)$ if $\nu\ge2$. The integrability of $V_4(0)$ (and at the same time, of the term involving $\text{const}_3$) is obtained from Lemma \ref{le:310} with $\eta=-v$, $b=0$, $p=0$ and the fact that $\int_{\mathbb{R}}e^{-\sigma_1|t|}|t|^{-v}dt<+\infty$. Doing so indeed yields
		\begin{align*}
		& \Bigg|\int_{S_2}|s_2|^2 |s_1|^{-1-v}\int_{\mathbb{R}} e^{-\sigma_1|t|} \Big|t +\dfrac{rs_2}{s_1}\Big|^{-v} - |t|^{-v}|s_2|^2 |s_1|^{-1-v} dt \Bigg| \Gamma(d\bs) \\
		& \hspace{5cm} \le \int_{S_2}\int_{\mathbb{R}} e^{-\sigma_1|t|} \Bigg|\Big|t +\dfrac{rs_2}{s_1}\Big|^{-v} - |t|^{-v}\Bigg|dt    \Gamma(d\bs) \\
		& \hspace{5cm} \le \text{const} \int_{S_2}|s_1|^{-\nu}|s_1|^{\nu-1-v}\Gamma(d\bs)\\
		& \hspace{5cm} \le \text{const} \int_{S_2}|s_1|^{-\nu}\Gamma(d\bs)\\
		& \hspace{5cm} < +\infty,
		\end{align*}
		since $\nu-1-v = \dfrac{\nu-1}{2}>0$ if $\nu\in(1,2)$ and $\nu-1-v>\nu-2>0$ if $\nu\ge2$. The convergence $\int V_4(h) \rightarrow \int V_4(0)$ can be obtained from Lemma \ref{le:312} with $\eta=-v$, $b=0$ and $0<p<v$. The generalised dominated convergence hence applies to $L_1$.\\
		We turn to $L_2$. Its integrand converges to
		\begin{align*}
		& e^{-\int_{S_2}|ts_1+rs_2|\Gamma(d\bs)} \sin\bigg(tx+a\int_{S_2}(ts_1+rs_2)\ln|ts_1+rs_2|\Gamma(d\bs)\bigg)\\
		& \hspace{4cm} \times \bigg(x+a\int_{S_2}s_1(1+\ln|ts_1+rs_2|)\Gamma(d\bs)\bigg)\ln|ts_1'+rs_2'| {s_2'}^2{s_1'}^{-1}.
		\end{align*}
		Applying the mean value theorem to the cosine function and the usual bounds, we can bound it by
		\begin{align}
		& e^{\sigma_2|r|-\sigma_1|t|}  \Big|s_2'^2 s_1'^{-1}\ln|ts_1'+rs_2'|\Big| \nonumber\\
		& \hspace{0.25cm} \dfrac{1}{\left|\frac{hs_2'}{s_1'}\right|} \Bigg|-\dfrac{hs_2'}{s_1'}x+a\int_{S_2}\bigg(\Big(t-\dfrac{hs_2'}{s_1'}\Big)s_1+rs_2\bigg)\ln\Big|\Big(t-\dfrac{hs_2'}{s_1'}\Big)s_1+rs_2\Big|-(ts_1+rs_2)\ln|ts_1+rs_2|\Gamma(d\bs)\Bigg|\nonumber\\
		& \le e^{\sigma_2|r|-\sigma_1|t|} \Big|s_2'^2 s_1'^{-1}\ln|ts_1'+rs_2'|\Big|\nonumber\\
		& \hspace{0.5cm} \Bigg(|x|+\dfrac{a}{\left|\frac{hs_2'}{s_1'}\right|} \int_{S_2}\bigg|\bigg(\Big(t-\dfrac{hs_2'}{s_1'}\Big)s_1+rs_2\bigg)\ln\Big|\Big(t-\dfrac{hs_2'}{s_1'}\Big)s_1+rs_2\Big|-(ts_1+rs_2)\ln|ts_1+rs_2|\bigg|\Gamma(d\bs)\Bigg).\label{dem:l2}
		\end{align}
		The term involving $|x|$ can be treated using the usual arguments. The one with the integral is of course the most delicate. Let us split this integral into two parts as:
		\begin{align*}
		& \int_{S_2}\dfrac{1}{\left|\frac{hs_2'}{s_1'}\right|}\bigg|\bigg(\Big(t-\dfrac{hs_2'}{s_1'}\Big)s_1+rs_2\bigg)\ln\Big|\Big(t-\dfrac{hs_2'}{s_1'}\Big)s_1+rs_2\Big|-(ts_1+rs_2)\ln|ts_1+rs_2|\bigg|\Gamma(d\bs)\\
		& \hspace{2cm} := Q_1 + Q_2,
		\end{align*}
		where $Q_1$ and $Q_2$ involve integrals over $S_2\cap\{\bs:|ts_1+rs_2|\ge2|hs_2'/s_1'|\}$ and $S_2\cap\{\bs:|ts_1+rs_2|<2|hs_2'/s_1'|\}$ respectively. We will first majorise $Q_1$ and $Q_2$, and then use these bounds in inequality \eqref{dem:l2}. Consider $Q_2$ and define the function $g$ such that for any $z>0$
		\begin{equation*}
		g(z)  = \left\{\begin{array}{ccc}
		f(z)=z|\ln z|, & \text{if} & 0<z<e^{-1},\\
		z(2 + \ln z), & \text{if} & z\ge e^{-1}.
		\end{array}
		\right.
		\end{equation*}
		It is easily checked that $g$ is continuous, strictly increasing and such that for any $z>0$, $0\le f(z)\le g(z)$. The integrand of $Q_2$ can be bounded as
		\begin{align*}
		\dfrac{1}{\left|\frac{hs_2'}{s_1'}\right|} \Bigg(\bigg|f\bigg(\Big(t-\dfrac{hs_2'}{s_1'}\Big)s_1+rs_2\bigg)\bigg| + \Big|f\Big(ts_1+rs_2\Big)\Big|\Bigg) & \le \dfrac{1}{\left|\frac{hs_2'}{s_1'}\right|} \Bigg(\bigg|g\bigg(\Big(t-\dfrac{hs_2'}{s_1'}\Big)s_1+rs_2\bigg)\bigg| + \Big|g\Big(ts_1+rs_2\Big)\Big|\Bigg)\\
		& \le \dfrac{1}{\left|\frac{hs_2'}{s_1'}\right|} \Bigg(\bigg|g\bigg(\Big|\dfrac{3hs_2'}{s_1'}\Big|\bigg)\bigg| + \bigg|g\bigg(\Big|\dfrac{2hs_2'}{s_1}\Big|\bigg)\bigg|\Bigg) \\
		& \le \dfrac{2}{\left|\frac{hs_2'}{s_1'}\right|} g\Big(\dfrac{3hs_2'}{s_1'}\Big).
		\end{align*}
		By Lemma \eqref{le:li}, with bound further the right-hand side for any $v>0$ by
		\begin{align*}
		\dfrac{2}{\left|\frac{hs_2'}{s_1'}\right|} g\Big(\dfrac{3hs_2'}{s_1'}\Big) \le \text{const}_1 + \text{const}_2 \Big|\dfrac{3hs_2'}{s_1'}\Big|^{v} + \text{const}_3\Big|\dfrac{3hs_2'}{s_1'}\Big|^{-v}. 
		\end{align*}
		On the one hand if $\Big|\dfrac{3hs_2'}{s_1'}\Big|<e^{-1}$, given that $(3|ts_1+rs_2|/2)^{-v}>(3hs_2'/s_1')^{-v}$,
		\begin{align*}
		\text{const}_1 + \text{const}_2 \Big|\dfrac{3hs_2'}{s_1'}\Big|^{v} + \text{const}_3\Big|\dfrac{3hs_2'}{s_1'}\Big|^{-v} & \le \text{const}_1 + \text{const}_2 \Big|t+\dfrac{rs_2}{s_1}\Big|^{-v}|s_1|^{-v}.
		\end{align*}
		On the other hand if  $\Big|\dfrac{3hs_2'}{s_1'}\Big|\ge e^{-1}$, then for $|h|<|r|$,
		\begin{align}
		\text{const}_1 + \text{const}_2 \Big|\dfrac{3hs_2'}{s_1'}\Big|^{v} + \text{const}_3\Big|\dfrac{3hs_2'}{s_1'}\Big|^{-v} & \le \text{const}_1 + \text{const}_2|s_1'|^{-v}.\label{dem:l2_q2}
		\end{align}
		Focusing now on $Q_1$, we can use the mean value theorem to bound its integrand by
		\begin{align*}
		|s_1| \Big|1+ \ln|u|\Big|,
		\end{align*}
		for some $u\in\Big[ts_1+rs_2-hs_2's_1/s_1' \wedge ts_1+rs_2,ts_1+rs_2-hs_2's_1/s_1' \vee ts_1+rs_2\Big]$. Given that $|ts_1+rs_2|\ge2|hs_2'/s_1'|$, we have $|u|\in\Big[\frac{|ts_1+rs_2|}{2},2|ts_1+rs_2|\Big]$ and thus, we further bound the above inequality using Lemma \ref{le:li} for any $v>0$ by
		\begin{align}
		& |s_1|\Big(\text{const}_1 + \text{const}_2 |ts_1+rs_2|^v + \text{const}_3 |ts_1+rs_2|^{-v}  \Big)  \nonumber\\
		& \hspace{4cm} \le \text{const}_1 + \text{const}_2 |t|^v + \text{const}_3\Big|t+\dfrac{rs_2}{s_1}\Big|^{-v} |s_1|^{1-v}. \label{dem:l2_q1}
		\end{align}
		Hence, using \eqref{dem:l2_q2} and \eqref{dem:l2_q1} in \eqref{dem:l2}, and making use again of Lemma \eqref{le:li} to bound $\Big|\ln|ts_1'+rs_2'|\Big|$, we can bound integrand of $L_2$ for any $v>0$ by
		\begin{align*}
		& e^{-\sigma_1|t|} \bigg(\text{const}_1 + \text{const}_2 |t|^v + \text{const}_3 \Big|t+\dfrac{rs_2'}{s_1'}\Big|^{-v} \bigg) |s_1'|^{-1-v} \\
		& \hspace{2.5cm} \times \bigg(|x| + \text{const}_4 + \text{const}_5 |t|^{v} + \text{const}_6 |s_1'|^{-v} + \text{const}_7 \Big|t+\dfrac{rs_2}{s_1}\Big|^{-v} |s_1|^{1-v} \bigg)
		\end{align*}
		It can be shown that all the terms obtained after expansion can be bounded by functions integrable with respect to $t$ and $\Gamma$ using the usual combinations of either Lemma \ref{le:310} or Lemma \ref{le:cor31} with $\eta=-v$, $b=0$, $p=0$, 
		the fact that $\int_{\mathbb{R}}e^{-\sigma_1|t|}|t|^{-v}<+\infty$, $\int_{\mathbb{R}}e^{-\sigma_1|t|}|t|^{-2v}<+\infty$ for appropriately chosen values $v>0$, 
		and \eqref{eq:nu_cond} with $\nu>1$. The detail we have to pay attention to is precisely to chose an appropriate exponent $v>0$ so that it satisfies the constraint $\eqref{eq:nu_cond}$ 
		and ensures the finiteness of the two integrals in $t$. 
		The later imposes us to have $v\in(0,1/2)$. 
		Regarding the former, we identify that the most negative power of which $|s_1|$ appears in the above bound after expansion is $-1-2v$. 
		We need $\nu-1-2v>0$. 
		Choosing  $v=(\nu-1)/4$ if $1< \nu < 3$ and any $v\in(0,1/2)$ if $\nu\ge3$ enables to satisfy both constraints, validating the use of the dominated convergence theorem for $L_2$, and finally, for $B_2$ in \eqref{dem:j2}.\\
		\indent The proof is essentially similar, somewhat easier, for $B_1$ in \eqref{dem:j1} for which the only difficulty is to perform the <<appropriate integration by parts>> when it comes to differentiating the term involving $(ts_1+rs_2)^{<0>}$.

		\subsection{Evaluating at $r=0$}

		Since $\mathbb{E}\Big[X_2^2\Big|X_1=x\Big]=-\phi^{(2)}_{X_2|x}(0)$, we evaluate \eqref{eq:d2_alp1} at $r=0$ and get
		\begin{align*}
		\varphi_{\bX}(t,0) & = \exp\{-\sigma_1|t|-ia\sigma_1\beta_1t\ln|t| + it\mu_1\},\\
		A_1/2 & = \sigma_1^2 \Big((\kappa_1^2-a^2q_0^2)H_c(0)+2a\kappa_1q_0H_s(0)\Big)\\
		& \hspace{3cm} + 2a\lambda_1\sigma_1^2\Big(-aq_0H_c(1)+\kappa_1H_s(1)\Big) - a^2\lambda_1^2\sigma_1^2H_c(2),\\
		iA_2/2 & = \sigma_1\Big(-ak_1H_c(0)+\kappa_2H_s(0)\Big) - a\lambda_2\sigma_1H_c(1),\\
		A_3/2 & = \sigma_1 \Big((\sigma_1\kappa_2+a\mu_1k_1)H_c(0) + (\sigma_1ak_1-\mu_1\kappa_2)H_s(0)\Big)\\
		& \hspace{1cm} + a\sigma_1\Big((\lambda_2\mu_1-a\sigma_1\beta_1k_1)H_c(1) + \sigma_1(\lambda_2+\beta_1\kappa_2)H_s(1)\Big) - a^2\sigma_1^2\beta_1\lambda_2H_c(2),
		\end{align*}
		where $k_1=\sigma_1^{-1}\int_{S_2}(s_2/s_1)^{2}s_1\ln|s_1|\Gamma(d\bs)$, and the $H_c$'s and $H_s$'s are defined at Lemma \ref{le:HcHs}. Using the result of the same Lemma under $\beta_1\ne0$ and $\beta_1=0$, and regrouping the terms allows to retrieve the two formulae of Theorem \ref{theo:cond_moments_stable_general_alp1_scm} with
		\begin{align}
		& \hspace*{2cm} U(x)  = \int_{0}^{+\infty} e^{-\sigma_1 t} \sin\Big(t(x-\mu_1) + a\sigma_1\beta_1t\ln t \Big)dt, \label{def:uvw1} \\
		& \hspace*{2cm} V(x)  = \int_{0}^{+\infty} e^{-\sigma_1 t} (1+\ln t) \cos\Big(t(x-\mu_1) + a\sigma_1\beta_1t\ln t \Big)dt,\label{def:uvw2} \\
		& \hspace*{1.9cm} W(x)  = \int_{0}^{+\infty} e^{-\sigma_1 t} (1+\ln t)^2 \cos\Big(t(x-\mu_1) + a\sigma_1\beta_1t\ln t \Big)dt.\label{def:uvw3}
		\end{align}
		
		\subsection{Lemmas for the proof of Theorem \ref{theo:cond_moments_stable_general_alp1_scm}}
		
		\begin{lem}\label{le:li}
			For any $x>0$ and $v>0$
			\begin{align*}
			|\ln x| \le \dfrac{1}{v}\Big(2 + x^v + x^{-v}\Big).
			\end{align*}
		\end{lem}
		
		We provide here two Lemmas which are used in the proof of Theorem \ref{theo:cond_moments_stable_general_alp1_scm}.
		
		\begin{lem}
			\label{le:HcHs}
			Let for any $n\ge0$,
			\begin{align*}
			H_c(n) & = \int_{0}^{+\infty} e^{-\sigma_1 t} (1+\ln t)^n \cos\Big(t(x-\mu_1)+a\sigma_1\beta_1t\ln t\Big)dt,\\
			H_s(n) & = \int_{0}^{+\infty} e^{-\sigma_1 t} (1+\ln t)^n \sin\Big(t(x-\mu_1)+a\sigma_1\beta_1t\ln t\Big)dt.
			\end{align*}
			Then, if $\beta_1\ne0$,
			\begin{align*}
			H_c(1) & = \dfrac{1}{a\sigma_1\beta_1}\Big(\sigma_1H_s(0)-(x-\mu_1)H_c(0)\Big), & H_s(1) & = \dfrac{1}{a\sigma_1\beta_1}\Big(1 - \sigma_1H_c(0)-(x-\mu_1)H_s(0)\Big).
			\end{align*}
			If $\beta_1=0$,
			\begin{align*}
			H_c(0) & = \pi f_{X_1}(x),\\
			H_s(0) & = \dfrac{x-\mu_1}{\sigma_1}\pi f_{X_1}(x),\\
			H_s(1) - \dfrac{x-\mu_1}{\sigma_1}H_c(1)& = \dfrac{\pi F_{X_1}(x)}{\sigma_1}.
			\end{align*}
		\end{lem}
		
		\textit{Proof.}
		\noindent The equalities of Lemmas \ref{le:SCFG}-\ref{le:HcHs} can be obtained by integrating by parts. We provide details for the last equality of Lemma \ref{le:HcHs} when $\beta_1=0$. Integrating the exponential by parts, we obtain
		\begin{align*}
		H_s(1) = \dfrac{1}{\sigma_1} \int_{0}^{+\infty} e^{-\sigma_1 t} t^{-1} \sin\Big(t(x-\mu_1)\Big)dt + \dfrac{x-\mu_1}{\sigma_1} H_c(1)
		\end{align*}
		Denote $A(x) = \int_{0}^{+\infty} e^{-\sigma_1 t} t^{-1} \sin\Big(t(x-\mu_1)\Big)dt$ for $x\in\mathbb{R}$ ($A$ is well defined since $e^{-\sigma_1 t} t^{-1} \sin\Big(t(x-\mu_1)\Big) \rightarrow x-\mu_1$ as $t\rightarrow0$). It can be shown that we can derivate $A$ under the integral sign and get
		\begin{align*}
		A'(x) & = \int_{0}^{+\infty} e^{-\sigma_1 t} \cos\Big(t(x-\mu_1)\Big)dt = \pi f_{X_1}(x),
		\end{align*}
		Since $X_1$ is Cauchy distributed when $\alpha=1$ and $\beta_1=0$,
		\begin{align*}
		A(x) & = \pi F_{X_1}(x) + \text{const} = \text{Arctg}\Big(\dfrac{x-\mu_1}{\sigma_1}\Big) + \dfrac{\pi}{2} + \text{const},
		\end{align*}
		and evaluating the integral form of $A$ at $\mu_1$, we deduce that $\text{const}=-\pi/2$. Thus, $A(x) = \pi \Big(F_{X_1}(x)-1/2\Big)$.
		
		\label{page:theo24_pend}
		
		\section{Proof of Proposition \ref{prop:equiv}}
		\label{sec:prop_asymp}
		
		\subsection{Case $\boldsymbol{\alpha\ne1}$}
		\label{sec:propa1_alpnot1}
		
		First assume that $|\beta_1|\ne1$. We will focus on the case $x\rightarrow+\infty$. The case $x\rightarrow-\infty$ can be obtained by considering the vector $(X_1,X_2)$, whose parameter are $\beta_1^\ast=-\beta_1$, $\kappa_1^\ast=-\kappa_1$ and $\lambda_1^\ast=\lambda_1$ and noticing that $\mathbb{E}\Big[X_2^p\Big|X_1=x\Big] = \mathbb{E}\Big[X_2^p\Big|-X_1=-x\Big]$. For $p=1$, the result is already known (see \cite{har91}). For $p=2,3,4$, we have from the proofs of \eqref{eq:scm_gen}-\eqref{eq:fcm_gen}, that
		\begin{align*}
		\mathbb{E}\Big[X_2^p\Big|X_1=x\Big] & = \dfrac{\alpha \sigma_1^\alpha}{\pi f_{X_1}(x)} \bigg[x^{p-1}\mathcal{H}\Big(\alpha-1,(a\lambda_p,\kappa_p);x\Big) + \sum_{i=2}^pb_{i,p}x^{p-i}\mathcal{H}\Big(i(\alpha-1),\bnu_{i};x\Big)\bigg],
		\end{align*}
		for some coefficients $b$'s. From the proof of Corollary 3.2 in \cite{har91}, we deduce the following limit:
		\begin{align*}
		x^{\alpha}\mathcal{H}\Big(\alpha-1,(a\lambda_p,\kappa_p);x\Big) \underset{x\rightarrow+\infty}{\longrightarrow} \Big(\kappa_p   + \lambda_p\Big)\sin\Big(\frac{\pi \alpha}{2}\Big)\Gamma(\alpha).
		\end{align*}
		
		\vspace*{-0.4cm}
		
		\noindent We also have
		
		\vspace*{-1cm}
		
		\begin{align}\label{eq:f_decay}
		x^{\alpha+1}f_{X_1}(x) \underset{x\rightarrow+\infty}{\longrightarrow} \dfrac{1}{\pi}\sigma_1^\alpha(1+\beta_1)\sin\Big(\frac{\pi \alpha}{2}\Big)\Gamma(1+\alpha).
		\end{align}
		
		\vspace*{-0.4cm}
		
		\noindent Hence,
		
		\vspace*{-1cm}
		
		\begin{align*}
		x^{-p}\dfrac{\alpha \sigma_1^\alpha x^{p-1}}{\pi f_{X_1}(x)}\mathcal{H}\Big(\alpha-1,(a\lambda_p,\kappa_p);x\Big) \longrightarrow \dfrac{\kappa_p + \lambda_p}{1+\beta_1} ,
		\end{align*}
		as $x\rightarrow+\infty$. It remains to be shown that $\dfrac{\sum_{i=2}^pb_{i,p}x^{p-i}\mathcal{H}\Big(i(\alpha-1),\bnu_{i};x\Big)}{x^{p-1}\mathcal{H}\Big(\alpha-1,(a\lambda_p,\kappa_p);x\Big)} \underset{x\rightarrow+\infty}{\longrightarrow} 0$. 
		By Theorem 127 in \cite{tit48}, for $i=2,3,4$,
		
		\vspace*{-1cm}
		
		\begin{align*}
		\mathcal{H}\Big(i(\alpha-1),\bnu_{i};x\Big) \underset{x\rightarrow+\infty}{=} O\Big(x^{-i(\alpha-1)-1}\Big).
		\end{align*}
		
		\vspace*{-0.4cm}
		
		\noindent Hence, 
		
		\vspace*{-1cm}
		
		\begin{align*}
		\Bigg|\dfrac{x^{p-i}\mathcal{H}\Big(i(\alpha-1),\bnu_{i};x\Big)}{x^{p-1}\mathcal{H}\Big(\alpha-1,(a\lambda_p,\kappa_p);x\Big)} \Bigg| \underset{x\rightarrow+\infty}{=} O\Big(x^{\alpha(1-i)}\Big) \longrightarrow 0.
		\end{align*}
		
		\noindent Now assume that $|\beta_1|=1$. For instance if $\beta_1=1$, the distribution of $X_1$ is \textit{totally skewed to the right}. On the one hand, we have $\lambda_p=\beta_1\kappa_p$. On the other hand, the right tail of $f_{X_1}$ still decays as \eqref{eq:f_decay}, 
		yielding the conclusion.
		
		\subsection{Case $\boldsymbol{\alpha=1}$}
		\label{sec:propa1_alp1}
		
		The form of the conditional second order moment when $\alpha=1$ requires to distinguish the cases $\beta_1\ne0$ and $\beta_1=0$.
		
		\paragraph{Case $\boldsymbol{\beta_1\ne0}$}

		We only consider $|\beta_1|<1$ and $x\longrightarrow+\infty$, the other cases being similar.
		Since $|x|\rightarrow+\infty$, we have $x-\mu_1 \sim x$ and we may assume that $\mu_1=0$. From \cite{har91}, 
		we know that $U(x)\sim x^{-1}$. Notice that
		\begin{align*}
		W(x) & = \int_{0}^{+\infty} e^{-\sigma_1t}(1+\ln t)^2\cos(a\sigma_1\beta_1t\ln t)\cos(tx)dt \\
		& \hspace{0.5cm} - \int_{0}^{+\infty} e^{-\sigma_1t}(1+\ln t)^2\sin(a\sigma_1\beta_1t\ln t)\sin(tx)dt.
		\end{align*}
		Because the factors of $\cos(tx)$ and $\sin(tx)$ are integrable, we have by the Riemann-Lebesgue Lemma that $W(x) \underset{x\rightarrow+\infty}{\longrightarrow} 0$. Having also
		\begin{align*}
		f_{X_1}(x) & \sim \dfrac{\sigma_1(1+\beta_1)}{\pi}x^{-2},
		\end{align*}
		we deduce the following limits
		\begin{align*}
		\Big(2a\sigma_1q_0(\lambda_1-\beta_1\kappa_1)+2(\kappa_1\lambda_1-\lambda_2)x\Big)\dfrac{\sigma_1U(x)}{\beta_1 \pi f_{X_1}(x)}x^{-2} & \underset{x\rightarrow+\infty}{\longrightarrow } \dfrac{2(\kappa_1\lambda_1-\lambda_2)}{(1+\beta_1)\beta_1},\\
		\Big(\lambda_2+\beta_1\kappa_2-2\kappa_1\lambda_1 + a^2\sigma_1\beta_1(\lambda_1^2-\beta_1\lambda_2)W(x)\Big) \dfrac{\sigma_1x^{-2}}{\pi f_{X_1}(x)} & \underset{x\longrightarrow+\infty}{\longrightarrow } \dfrac{\lambda_2 + \beta_1\kappa_2-2\kappa_1\lambda_1}{(1+\beta_1)\beta_1}.
		\end{align*}
		Hence,
		\begin{align*}
		x^{-2}\mathbb{E}\Big[X_2^2\Big|X_1=x\Big] & \underset{x\rightarrow+\infty}{\longrightarrow} \dfrac{\lambda_2}{\beta_1} + \dfrac{2(\kappa_1\lambda_1-\lambda_2)}{(1+\beta_1)\beta_1} + \dfrac{\lambda_2 + \beta_1\kappa_2-2\kappa_1\lambda_1}{(1+\beta_1)\beta_1} = \dfrac{\kappa_2+\lambda_2}{1+\beta_1}
		\end{align*}
		
		\paragraph{Case $\boldsymbol{\beta_1=0}$}
		
		From \cite{har91},
		\begin{align*}
		V(x) \longrightarrow -\dfrac{\pi}{2x},
		\end{align*}
		hence,
		\begin{align*}
		2a\sigma_1\lambda_1\Big(a\sigma_1q_0-\kappa_1(x-\mu_1)\Big)\dfrac{V(x)}{\pi f_{X_1}(x)}x^{-2} \longrightarrow a\pi\lambda_1\kappa_1.
		\end{align*} 
		Moreover,
		\begin{align*}
		a\sigma_1\dfrac{F_{X_1}(x) -1/2}{f_{X_1}(x)} x^{-2}\longrightarrow \dfrac{1}{2} a\pi (\lambda_2-2\kappa_1\lambda_1).
		\end{align*}
		It can be shown that $W(x)\longrightarrow0$. Therefore,
		\begin{align*}
		x^{-2}\mathbb{E}\Big[X_2^2\Big|X_1=x\Big] & \underset{x\rightarrow+\infty}{\longrightarrow} \kappa_2 + \dfrac{1}{2}a\pi(\lambda_2-2\kappa_1\lambda_1)  + a\pi\kappa_1\lambda_1 = \kappa_2+\lambda_2
		\end{align*}
		\label{page:prop21_pend}

		\section{Proof of Lemma \ref{prop:multistable}}
		\label{sec:lem_spec_rep_MA}
		
		The characteristic function of $\boldsymbol{X_t}$ reads, for any $\bu=(u_1,u_2)\in\mathbb{R}^{2}$:
		\begin{align*}
		\varphi_{\boldsymbol{X_t}}(\bu) = \mathbb{E}\left(\exp\left\{i \sum_{j=1}^{2}u_jX_{j,t}\right\}\right)
		& = \prod_{k\in\mathbb{Z}}\mathbb{E}\left[i\left(\sum_{j=1}^{2}u_ja_{k,j} \right)\varepsilon_{t+k}\right].
		\end{align*}
		We obtain for $\alpha\ne1$,
		\begin{align}
		\varphi_{\boldsymbol{X_t}}(\bu) & =  \exp\left\{-\sum_{k\in\mathbb{Z}}\sigma^\alpha|\sum_{j=1}^{2}u_ja_{k,j}|^\alpha \left(1-i\beta \sign\Big(\sum_{j=1}^{2}u_ja_{k,j}\Big)\text{tg}\bigg(\dfrac{\pi\alpha}{2}\bigg)\right) +i\sum_{j=1}^{2}u_j \sum_{k\in\mathbb{Z}}a_{k,j}\mu\right\}.\label{dem:alp}
		\end{align}
		And for $\alpha=1$,
		\begin{align}
		\varphi_{\boldsymbol{X_t}}(\bu) & = \exp\left\{-\sum_{k\in\mathbb{Z}}\sigma|\sum_{j=1}^{2}u_ja_{k,j}|\left(1+i\beta\frac{2}{\pi}\sign\Big(\sum_{j=1}^{2}u_ja_{k,j}\Big)\ln\Big|\sum_{j=1}^{2}u_ja_{k,j}\Big|\right) + i  \sum_{j=1}^{2}u_j\sum_{k\in\mathbb{Z}}a_{k,j}\mu\right\}. \label{dem:alp1}
		\end{align}
		Replacing \eqref{def:spectral} in \eqref{def:char_fun_vec}, we retrieve the two above formulae.
		\label{page:lem31_pend}
		
		
		
		\section{Proof of the asymptotic moments in Section \ref{ex:ar1}}
		\label{sec:ex_ar1}
		
		The results in Section \ref{ex:ar1} follow from Proposition \ref{prop:2sidedma_condmom} applied to $X_t =\sum_{k\in\mathbb{Z}}\rho^{k}\mathds{1}_{\{k\ge0\}} \varepsilon_{t+k}$.
		Regarding the asymptotic behaviours of moments, we give the proof for the excess kurtosis. 
		The other limits and equivalents are obtained in a similar manner. 
		Letting $\alpha\in(3/2,2)$ ensures the existence of the fourth order moment. 
		Since we assume $\rho>0$, it follows that $\lambda_p = \beta_1 \kappa_p$ for $p=1,2,3,4$. 
		Using Proposition \ref{prop:equiv}, one can show that as $x$ tends to infinity
		\begin{align*}
		\gamma_2(x,h) & \longrightarrow \dfrac{\kappa_4 -4\kappa_1\kappa_3+6\kappa_1^2\kappa_2-3\kappa_1^4}{\Big(\kappa_2-\kappa_1^2\Big)^2} - 3.
		\end{align*}
		Substituting the $\kappa_p$'s by $\rho^{h(\alpha-p)}$ and rearranging terms yields the conclusion.

		\section{Proof of Proposition \ref{prop:ar1_bernoulli}}
		\label{sec:prop_ar1_bernoulli}
		
		From Proposition 5.2.4 p.110 and Equation (15.3.9) p.438 in \cite{ks20} applied to the noncausal AR(1) process $(X_t)$ with MA($\infty$) representation $X_t=\sum_{k\in\mathbb{Z}}\rho^{k}\mathds{1}_{\{k\ge0\}}\varepsilon_{t+k}$, we have that,
    $\mathbb{P}\bigg(\dfrac{X_{t+h}}{|X_t|}\in A\bigg||X_t|>x\bigg)\longrightarrow \mathbb{P}\Big(\Theta_h\in A\Big)$,
    as $x\rightarrow+\infty$, for any continuity set $A\subset \mathbb{R}$, and $\Theta_h$ is defined by $\Theta_h:= S\rho^{-h}\mathds{1}_{\{h+\tau\le0\}}$, with $S$ a discrete random variable such that $\mathbb{P}(S=1)=1-\mathbb{P}(S=-1)=c:=\dfrac{1+\beta}{2}$, and $\tau$ a discrete random variable independent from $S$ such that $\mathbb{P}(\tau=k) = \rho^{-\alpha k}(1-\rho^{\alpha})\mathds{1}_{\{k\le 0\}}$, for all $k\in\mathbb{Z}$.
    For $A=[\rho^{-h}-\delta,\rho^{-h}+\delta]$ with $\delta\in(0,\rho^{-h})$, we have $[1-\rho^{h}\delta,1+\rho^{h}\delta]\cap\{0,1\}=\{1\}$, $[-1-\rho^{h}\delta,-1+\rho^{h}\delta]\cap\{0,1\}=\emptyset$ and
    \begin{align*}
    \mathbb{P}\Big(\Theta_h\in A\Big) & = \mathbb{P}\Big(S\rho^{-h}\mathds{1}_{\{h+\tau\le0\}}\in [\rho^{-h}-\delta,\rho^{-h}+\delta]\Big)\\
    & = \mathbb{P}\Big(S\mathds{1}_{\{\tau\le-h\}}\in [1-\delta\rho^{h},1+\delta\rho^{h}]\Big)\\
    & = \mathbb{P}\Big(\mathds{1}_{\{\tau\le-h\}}\in [1-\delta\rho^{h},1+\delta\rho^{h}]\Big)\mathbb{P}\Big(S=1\Big) + \mathbb{P}\Big(\mathds{1}_{\{\tau\le-h\}}\in [-1-\delta\rho^{h},-1+\delta\rho^{h}]\Big)\mathbb{P}\Big(S=-1\Big)\\
    & = \mathbb{P}\Big(\tau\le-h\Big) c\\
    & = c \sum_{k\le -h} \rho^{-\alpha k}(1-\rho^{\alpha})\mathds{1}_{\{k\le 0\}}\\
    & = c \rho^{\alpha h}.
    \end{align*}
    Thus, $\mathbb{P}\bigg(\dfrac{X_{t+h}}{|X_t|}\in [\rho^{-h}-\delta,\rho^{-h}+\delta]\bigg||X_t|>x\bigg)\longrightarrow c \rho^{\alpha h}$ as $x\rightarrow+\infty$.
    Similarly, it can be shown that $\mathbb{P}\bigg(\dfrac{X_{t+h}}{|X_t|}\in [-\rho^{-h}-\delta,-\rho^{-h}+\delta]\bigg||X_t|>x\bigg)\longrightarrow (1-c) \rho^{\alpha h}$ as $x\rightarrow+\infty$. Hence, for $s\in\{-1,+1\}$, provided $\mathbb{P}(S=s)>0$,
    \begin{align*}
    \mathbb{P}\bigg(\dfrac{X_{t+h}}{X_t}\in [\rho^{-h}-\delta,\rho^{-h}+\delta]\bigg|sX_t>x\bigg) & = \mathbb{P}\bigg(\dfrac{X_{t+h}}{|X_t|}\in [s\rho^{-h}-\delta,s\rho^{-h}+\delta]\bigg||X_t|>x,sX_t>x\bigg) \\
    & = \dfrac{\mathbb{P}\bigg(\dfrac{X_{t+h}}{|X_t|}\in [s\rho^{-h}-\delta,s\rho^{-h}+\delta]\bigg||X_t|>x\bigg)}{\mathbb{P}\Big(sX_t>x\bigg||X_t|>x\Big)}\\
    & \longrightarrow \dfrac{\rho^{\alpha h} \mathbb{P}(S=s) }{\mathbb{P}(S=s)} = \rho^{\alpha h}.
    \end{align*}
    The proof for the limit of $\mathbb{P}\bigg(\dfrac{X_{t+h}}{X_t}\in[-\delta,\delta]\bigg|sX_t >x\bigg)$ is similar.
    
    \label{page:prop41_pend}

    \section{Proof of Proposition \ref{prop:mar1q_bernoulli}}
    \label{sec:prop_mar1q_bernoulli}
    
    Assume Assumption \ref{asu:not_exp_colin} holds for some $\epsilon>0$. Let us first show that $d>0$.
    \textit{Ad absurdum}, assume that $d=0$. 
    Then there exists $k,\ell\in\{0,\ldots,h\}$, $k\ne\ell$, say $k<\ell$, such that $A_k = A_\ell$.
    This entails 
    $$
    (\underbrace{\rho^{-1},\ldots,\rho^{-k}}_{k },\underbrace{\rho^{-k}a_{-1},\ldots,\rho^{-k}a_{-(h-k)}}_{h-k}) = (\underbrace{\rho^{-1},\ldots,\rho^{-\ell}}_{\ell},\underbrace{\rho^{-\ell}a_{-1},\ldots,\rho^{-\ell}a_{-(h-\ell)}}_{h-\ell}).
    $$
    Since $k<\ell$, the above equality implies that $\rho^{-k}a_{-1}=\rho^{-k-1}$, and hence $a_{-1} = \rho^{-1}$.
    But then for $k=-1$ and $\lambda=\rho$,
    \begin{align*}
    \lambda(a_{k+m},\ldots,a_{k+2},a_{k+1},a_k) & = \rho(\rho^{m-1},\ldots,\rho,1,\rho^{-1}) = \brho, 
    \end{align*}
    which violates Assumption \ref{asu:not_exp_colin}.
    Hence, $d>0$.\\
    
    Let us now establish the main result of Proposition \ref{prop:mar1q_bernoulli}. For $u>0$ and $\bx_0\in\mathbb{R}^n$, denote generically $\mathcal{B}_u(\bx_0)=\left\{\bx\in\mathbb{R}^{n}: \hspace{0.3cm} \big|\big|\bx-\bx_0\big|\big|<u\right\}$.
    By Point $\iota\iota)$ Proposition 5.2.4 p.110 and Equation (15.3.9) p.438 in \cite{ks20},
    \begin{align}
    \mathbb{P}\Bigg(\dfrac{\boldsymbol{\overline{X}}_{t+h}}{|X_t|} \in \mathcal{B}_{\delta}(s\bA_k) \Bigg||X_t|>x,\dfrac{\boldsymbol{\underline{X}}_{t}}{|X_t|}\in\mathcal{B}_{\eta}\big(s\brho\big)\Bigg) & = \dfrac{\mathbb{P}\Bigg(\dfrac{(\boldsymbol{\underline{X}}_{t},\boldsymbol{\overline{X}}_{t+h})}{|X_t|} \in \mathcal{B}_{\eta}\big(s\brho\big)\times\mathcal{B}_{\delta}(s\bA_k) \Bigg||X_t|>x\Bigg)}{\mathbb{P}\Bigg(\dfrac{\boldsymbol{\underline{X}}_{t}}{|X_t|} \in \mathcal{B}_{\eta}\big(s\brho\big) \Bigg||X_t|>x\Bigg)} \nonumber \\
    & \underset{x\rightarrow+\infty}{\longrightarrow} \dfrac{\mathbb{P}\bigg((\Theta_{-m},\ldots,\Theta_{0},\Theta_1,\ldots,\Theta_h)\in\mathcal{B}_{\eta}\big(s\brho\big)\times\mathcal{B}_{\delta}(s\bA_k)\bigg)}{\mathbb{P}\bigg((\Theta_{-m},\ldots,\Theta_{0})\in\mathcal{B}_{\eta}\big(s\brho\big)\bigg)},\label{eq:limit_tail_proc}
    \end{align}
    where for all $k$, $\Theta_k = S\dfrac{a_{-N-k}}{|a_{-N}|}$, with $\mathbb{P}(S=1)=1-\mathbb{P}(S=-1)=c:=\dfrac{1+\beta}{2}$, $\mathbb{P}(N=j) = \dfrac{|a_{-j}|^\alpha}{\sum_{\ell\in\mathbb{Z}}|a_\ell|^\alpha}$ for all $j\in\mathbb{Z}$, and $S$ and $N$ are independent.
    There is no issue with division by zero since $a_{-N}=0$ implies $\mathbb{P}(N=j)=0$.
    Let us first focus on the denominator in \eqref{eq:limit_tail_proc}.
    We have
    \begin{align}
    \mathbb{P}\bigg((\Theta_{-m},\ldots,\Theta_{0})\in\mathcal{B}_{\eta}\big(s\brho\big)\bigg) & = \mathbb{P}\bigg(S(a_{m-N},\ldots,a_{1-N},a_{-N})/|a_{-N}|\in\mathcal{B}_{\eta}\big(s\brho\big)\bigg)\nonumber\\
    & = \sum_{\zeta\in\{-1,+1\}}\mathbb{P}\bigg(\zeta(a_{m-N},\ldots,a_{1-N},a_{-N})/|a_{-N}|\in\mathcal{B}_{\eta}\big(s\brho\big)\bigg)\mathbb{P}(S=\zeta) \nonumber\\
    & = \sum_{\zeta\in\{-1,+1\}}\sum_{j\in\mathbb{Z}}\mathds{1}_{\left\{\zeta(a_{m-j},\ldots,a_{1-j},a_{-j})/|a_{-j}|\in\mathcal{B}_{\eta}\big(s\brho\big)\right\}}\mathbb{P}(N=j)\mathbb{P}(S=\zeta) \nonumber \\
    & = \sum_{\zeta\in\{-1,+1\}}\sum_{j\in\mathcal{J}_\zeta}\dfrac{|a_{j}|^\alpha}{\sum_{\ell\in\mathbb{Z}}|a_\ell|^\alpha}\mathbb{P}(S=\zeta),\label{eq:proba_Js}
    \end{align}
    where for $\zeta\in\{-1,1\}$
    $$
    \mathcal{J}_\zeta := \left\{j\in\mathbb{Z}: \hspace{0.3cm} \zeta\dfrac{(a_{m+j},\ldots,a_{1+j},a_{j})}{|a_{j}|}\in\mathcal{B}_\eta(s\brho)\right\}.
    $$
    By convention, if $a_j=0$, then the index $j$ drops from the sum above and thus from $\mathcal{J}_\zeta$. Let us show that $\mathcal{J}_s = \{j\ge0\}$ and $\mathcal{J}_{-s} = \emptyset$. Notice first that for any $\zeta,s\in\{-1,1\}$, Assumption \ref{asu:not_exp_colin} guarantees that
    $$
    \bigg|\bigg|\zeta\dfrac{(a_{m+j},\ldots,a_{1+j},a_{j})}{|a_{j}|}-s\brho\bigg|\bigg| = \bigg|\bigg|\dfrac{s\zeta}{|a_j|}(a_{m+j},\ldots,a_{1+j},a_{j})-\brho\bigg|\bigg| >\epsilon>\eta,
    $$
    for all $j\le-1$ such that $a_j\ne0$. Thus, $\zeta(a_{m+j},\ldots,a_{1+j},a_{j})/|a_{j}|\not\in\mathcal{B}_\eta(s\brho)$ for $j\le-1$, and $\zeta,s\in\{-1,1\}$. Hence, $\mathcal{J}_\zeta\subset\{j\ge0\}$, for $\zeta\in\{-1,1\}$.
    Now, for $j\ge0$, since $a_j=\rho^{j}$
    $$
    \zeta\dfrac{(a_{m+j},\ldots,a_{1+j},a_{j})}{|a_{j}|} = \zeta\dfrac{(\rho^{m+j},\ldots,\rho^{1+j},\rho^{j})}{|\rho^{j}|} = \zeta\brho \in \mathcal{B}_\eta(\zeta\brho).
    $$
    Hence, $\{j\ge0\}\subset\mathcal{J}_s$  which shows that $\mathcal{J}_s = \{j\ge0\}$.
    However, for $j\ge0$, $-s(a_{m+j},\ldots,a_{1+j},a_{j})/|a_{j}| = -s\brho$, and $\|(-s\brho)-s\brho\|=2\|\brho\|>2\epsilon>\eta$ by Assumption \ref{asu:not_exp_colin} with $\lambda=0$, and $-s\brho\not\in\mathcal{B}_\eta(s\brho)$. 
    Thus, $\mathcal{J}_{-s} = \emptyset$. \\
    Therefore, \eqref{eq:proba_Js} yields
    \begin{align*}
    \sum_{\zeta\in\{-1,+1\}}\sum_{j\in\mathcal{J}_\zeta}\dfrac{|a_{j}|^\alpha}{\sum_{\ell\in\mathbb{Z}}|a_\ell|^\alpha}\mathbb{P}(S=\zeta) & = \sum_{j\ge0}\dfrac{|a_{j}|^\alpha}{\sum_{\ell\in\mathbb{Z}}|a_\ell|^\alpha}\mathbb{P}(S=s)\\
    & = \dfrac{\mathbb{P}(S=s)}{\sum_{\ell\in\mathbb{Z}}|a_\ell|^\alpha}\sum_{j\ge0} \rho^{j \alpha}\\
    & = \dfrac{\mathbb{P}(S=s)}{(1-\rho^\alpha)\sum_{\ell\in\mathbb{Z}}|a_\ell|^\alpha},
    \end{align*}
    which shows that:
    \begin{align}\label{eq:theta_denom}
    \mathbb{P}\bigg((\Theta_{-m},\ldots,\Theta_{0})\in\mathcal{B}_{\eta}\big(s\brho\big)\bigg) = \dfrac{\mathbb{P}(S=s)}{(1-\rho^\alpha)\sum_{\ell\in\mathbb{Z}}|a_\ell|^\alpha}.    
    \end{align}
    Let us now turn to the numerator in \eqref{eq:limit_tail_proc}.
    Proceeding as above, we obtain that
    \begin{align}
    \mathbb{P}\bigg((\Theta_{-m},\ldots,\Theta_{0},\Theta_1,\ldots,\Theta_h)\in\mathcal{B}_{\eta}\big(s\brho\big)\times\mathcal{B}_{\delta}(s\bA_k)\bigg) & = \sum_{\zeta\in\{-1,+1\}}\sum_{j\in\mathcal{J}_\zeta}\dfrac{|a_{j}|^\alpha}{\sum_{\ell\in\mathbb{Z}}|a_\ell|^\alpha}\mathbb{P}(S=\zeta),\label{eq:proba_Js_numerator}
    \end{align}
    where for $\zeta\in\{-1,1\}$
    $$
    \mathcal{J}_\zeta := \left\{j\in\mathbb{Z}: \hspace{0.3cm} \zeta\dfrac{(a_{m+j},\ldots,a_{1+j},a_{j},a_{j-1},\ldots,a_{j-h})}{|a_{j}|}\in\mathcal{B}_\eta(s\brho)\times \mathcal{B}_\delta(s\bA_k)\right\}.
    $$
    With similar considerations as above regarding the part $\zeta\dfrac{(a_{m+j},\ldots,a_{1+j},a_{j})}{|a_{j}|}\in\mathcal{B}_\eta(s\brho)$, we obtain that $\mathcal{J}_s\subset\{j\ge0\}$ whereas $\mathcal{J}_{-s}=\emptyset$.
    Also,
    \begin{align*}
  \mathcal{J}_s & = \left\{j\ge0: \hspace{0.3cm} s\dfrac{(a_{m+j},\ldots,a_{1+j},a_{j},a_{j},a_{j-1},\ldots,a_{j-h})}{|a_{j}|}\in \mathcal{B}_\eta(s\brho)\times\mathcal{B}_\delta(s\bA_k)\right\}\\
  & = \left\{j\ge0: \hspace{0.3cm} \bA_j\in \mathcal{B}_\delta(\bA_k)\right\},
    \end{align*}
    where $\bA_j=\bA_h$ for all $j\ge h$.
    Since $\delta<d$, we have by definition of $d$ that $\bA_j\in\mathcal{B}_\delta(\bA_k)$ if and only if $j = k$ in the case $0\le k \le h-1$, and $\bA_j\in\mathcal{B}_\delta(\bA_h)$ if and only if $j\ge h$, that is
    $$
    \mathcal{J}_s = \left\{\begin{array}{cc}
         \{k\}, & \text{ if } k\in\{0,\ldots,h-1\},  \\
         \{j\ge h\}, & \text{ if } k=h.
    \end{array}\right.
    $$
    Therefore, \eqref{eq:proba_Js_numerator} yields
    \begin{align*}
    \mathbb{P}\bigg((\Theta_{-m},\ldots,\Theta_{0},\Theta_1,\ldots,\Theta_h)\in\mathcal{B}_{\eta}\big(s\brho\big)\times\mathcal{B}_{\delta}(s\bA_k)\bigg) & = \left\{\begin{array}{cc}
         \dfrac{\mathbb{P}(S=s) \rho^{\alpha k}}{ \sum_{\ell\in\mathbb{Z}}|a_\ell|^\alpha}, & \text{ if } k\in\{0,\ldots,h-1\},  \\
         & \\
         \dfrac{\mathbb{P}(S=s) \rho^{\alpha h}}{(1-\rho^\alpha) \sum_{\ell\in\mathbb{Z}}|a_\ell|^\alpha}, & \text{ if } k=h.
    \end{array}\right.
    \end{align*}
    Finally, combining this with \eqref{eq:limit_tail_proc} and \eqref{eq:theta_denom}, and provided $\mathbb{P}(S=s)>0$, we obtain that 
    $$
    \mathbb{P}\Bigg(\dfrac{\boldsymbol{\overline{X}}_{t+h}}{|X_t|} \in \mathcal{B}_{\delta}(s\bA_k) \Bigg||X_t|>x,\dfrac{\boldsymbol{\underline{X}}_{t}}{|X_t|}\in\mathcal{B}_{\eta}\big(s\brho\big)\Bigg) = \left\{\begin{array}{cc}
         \rho^{\alpha k}(1-\rho^\alpha), & \text{ if } k\in\{0,\ldots,h-1\},  \\
         & \\
         \rho^{\alpha h}, & \text{ if } k=h.
    \end{array}\right.
    $$
    which concludes the proof.
    \label{page:prop42_pend}
    

\section{Proof of Proposition \ref{prop:crashodds_noncausal_ma}}

\label{page:prop_noncausMA}
    
        We first have by Bayes formula that for any $h\ge1$, $\delta\in(0,\epsilon^h)$, $s\in\{-1,+1\}$:
        \begin{align*}
        \mathbb{P}\bigg(\dfrac{X_{t+h}}{X_t}\in[-\delta,\delta]\bigg|sX_t>x\bigg) & = \mathbb{P}\bigg(\dfrac{X_{t+h}}{X_t}\in[-\delta,\delta]\bigg||X_t|>x,\sign(X_t)=s\bigg)\\
        & = \dfrac{\mathbb{P}\bigg(\dfrac{(X_t,X_{t+h})}{|X_t|}\in\{s\}\times[-\delta,\delta]\bigg||X_t|>x\bigg)}{\mathbb{P}\Big(sX_t>x\Big||X_t|>x\Big)}.
        \end{align*}
        Letting $S$ a random variable such that $\mathbb{P}(S=1)=1-\mathbb{P}(S=-1)=c:=\dfrac{1+\beta}{2}$, then \eqref{eq:reg_var_errors} implies that
        $$
        \mathbb{P}\Big(sX_t>x\Big||X_t|>x\Big) \underset{x\rightarrow+\infty}{\longrightarrow} \mathbb{P}(S=s).
        $$
        Now, by Proposition 5.2.4 p.110 and Equation (15.3.9) p.438 in \cite{ks20} applied to $X_t=\sum_{k\in\mathbb{Z}}a_k \varepsilon_{t+k}$, we have that:
        \begin{align*}
        \mathbb{P}\bigg(\dfrac{(X_t,X_{t+h})}{|X_t|}\in\{s\}\times[-\delta,\delta]\bigg||X_t|>x\bigg) \underset{x\rightarrow+\infty}{\longrightarrow} \mathbb{P}\Big((\Theta_0,\Theta_h)\in\{s\}\times[-\delta,\delta]\Big),
        \end{align*}
        where the $\Theta_j$'s are random variables such that $\Theta_j := S \dfrac{a_{-N-j}}{|a_{-N}|}$ for all $j\in\mathbb{Z}$, with $N$ the random variable such that $\mathbb{P}(N=j) = |a_{-j}|^{\alpha} / \sum_{k\in\mathbb{Z}}|a_k|^\alpha$, and with $S$ and $N$ furthermore independent.
        There is no issue of division by zero since for $a_{-j}=0$ we have $\mathbb{P}(N=j)=0$.
        Now, since $|a_k|=a_k$ for all $k\in\mathbb{Z}$,
        \begin{align*}
        \mathbb{P}\Big((\Theta_0,\Theta_h)\in\{s\}\times[-\delta,\delta]\Big) & = \mathbb{P}\bigg(S\Big(1,\dfrac{a_{-N-h}}{a_{-N}}\Big)\in\{s\}\times[-\delta,\delta]\bigg)\\
        & = \mathbb{P}(S=s)\mathbb{P}\bigg(\dfrac{a_{-N-h}}{a_{-N}}\in[-\delta,\delta]\bigg)\\
        & = \mathbb{P}(S=s)\sum_{j\in\mathcal{J}_{h}} \dfrac{|a_{j}|^\alpha}{\sum_{k\in\mathbb{Z}}|a_k|^\alpha},
        \end{align*}
		where $\mathcal{J}_{h}$ is the set of indexes defined by
		$$
		\mathcal{J}_{h} := \bigg\{j\in\mathbb{Z}: \hspace{0.15cm} a_j\ne 0 \text{ and } \dfrac{a_{j-h}}{a_j}\in[-\delta,\delta]\bigg\}.
		$$
		By assumption, $a_j=0$ for all $j<0$, which implies that $\mathcal{J}_{h}\subset \{j\ge0\}$.
		Also by assumption, we have that for any $j\ge0$, $a_{j}/a_{j+1} > \epsilon > 0$, which implies that for all $j \ge h$
		$$
		\dfrac{a_{j-h}}{a_{j}} > \epsilon^h > \delta.
		$$
		Hence, $a_{j-h}/a_{j} \not\in[-\delta,\delta]$ for all $j\ge h$.
		Last, for $j\in\{0,\ldots,h-1\}$, we have that $a_{j-h}=0$ and thus $a_{j-h}/a_{j} \in[-\delta,\delta]$.
		We deduce that $\mathcal{J}_{h}=\{0,\ldots,h-1\}$, and therefore
		$$
		\sum_{j\in\mathcal{J}_{h}} \dfrac{|a_{j}|^\alpha}{\sum_{k\in\mathbb{Z}}|a_k|^\alpha} = \dfrac{\sum_{j=0}^{h-1}|a_j|^\alpha}{\sum_{k\in\mathbb{Z}}|a_k|^\alpha}.
		$$
		Finally, we conclude that provided $\mathbb{P}(S=s)>0$
		\begin{align*}
		   \mathbb{P}\bigg(\dfrac{X_{t+h}}{X_t}\in[-\delta,\delta]\bigg|sX_t>x\bigg) & = \dfrac{\mathbb{P}\bigg(\dfrac{(X_t,X_{t+h})}{|X_t|}\in\{s\}\times[-\delta,\delta]\bigg||X_t|>x\bigg)}{\mathbb{P}\Big(sX_t>x\Big||X_t|>x\Big)}\\
		   & \underset{x\rightarrow+\infty}{\longrightarrow} \dfrac{\mathbb{P}\Big((\Theta_0,\Theta_h)\in\{s\}\times[-\delta,\delta]\Big)}{\mathbb{P}(S=s)}\\
		   & = \dfrac{ \mathbb{P}(S=s)\dfrac{\sum_{j=0}^{h-1}|a_j|^\alpha}{\sum_{k\in\mathbb{Z}}|a_k|^\alpha}}{\mathbb{P}(S=s)}\\
		   & =  \dfrac{\sum_{j=0}^{h-1}|a_j|^\alpha}{\sum_{k\in\mathbb{Z}}|a_k|^\alpha}.
		\end{align*}
		Which concludes the proof.
		\label{page:prop43_pend}

		\begin{center}
			{\sc Additional References}
		\end{center}
		

	\end{appendices}

\end{document}